\documentclass[11pt]{article}
\usepackage{amsmath} 
\usepackage{amssymb,latexsym,times}
\usepackage{mathptmx}
\usepackage{amscd}
\usepackage{stmaryrd}
\usepackage[all]{xy}
\catcode`\@=11

\catcode`\@=11
\renewcommand\subsubsection{\@startsection{subsubsection}{2}{\z@}%
                                     {-3.25ex\@plus -1ex \@minus -.2ex}%
                                     {-0.01 mm}
                                     {\normalfont\bfseries}}
\newtheorem{Thm}{Theorem}[section]
\newtheorem{Lem}[Thm]{Lemma}
\newtheorem{Cor}[Thm]{Corollary}
\newtheorem{Prop}[Thm]{Proposition}

\newtheorem{example}[Thm]{Example}
\newtheorem{remark}[Thm]{Remark}

\newcommand{\A}{\mathcal{A}}

\newcommand{\Z}{\mathbb{Z}}
\newcommand{\Q}{\mathbb{Q}}
\newcommand{\RR}{\mathbb{R}}

\newcommand{\N}{\mathbb{N}}
\newcommand{\C}{\mathbb{C}}
\newcommand{\R}{\mathcal{R}}

\newcommand{\M}{{\mathcal M}}

\newcommand{\rk}{\operatorname{rk}}

\newcommand{\md}{\mathrm{mod}}

\newcommand{\Hom}{\operatorname{Hom}} 
\newcommand{\Ext}{\operatorname{Ext}}
\newcommand{\Ex}{\operatorname{Ex}}

\newcommand{\Aut}{\operatorname{Aut}}

\newcommand{\dimv}{\underline{\dim}}

\newcommand{\bsm}{\begin{smallmatrix}}
\newcommand{\esm}{\end{smallmatrix}}

\def\CC{{\mathcal C}}
\def\resp{{\em resp.\ }}

\def\<{\langle\,}
\def\>{\,\rangle}

\def\eg{{\em e.g. }}

\def\T{{\mathcal T}}
\def\TT{{\mathbf T}}

\def\B{{\mathcal B}}
\def\U{{\mathcal U}}

\def\bE{{\mathbf E}}
\def\g{\mathfrak g}

\def\<{\langle}
\def\>{\rangle}
\def\n{{\mathfrak n}}
\def\si{\sigma}

\def\1{\mathbf 1}

\def\aa{{\mathbf a}}
\def\cc{{\mathbf c}}
\def\ii{{\mathbf i}}

\def\f{{\mathfrak f}}

\def\a{\alpha} 
\def\Y{\mathcal{Y}}
\def\YY{\mathbf{Y}}

\def\L{\Lambda}
\def\mod{{\rm mod}\,}

\def\Ext{{\rm Ext}}

\def\Hom{{\rm Hom}}

\def\K{\mathcal{K}}
\def\de{\delta}
\def\De{\Delta}
\def\Sl{\mathfrak{sl}}
\def\la{\lambda}
\def\AA{\mathcal{A}}

\def\hL{\widehat{\Lambda}}

\def\M{{\mathfrak M}}
\def\O{{\mathcal O}}
\def\1{{\mathbf 1}}
\def\ie{{\em i.e. }}

\def\h{\mathfrak{h}}
\def\hI{\widehat{I}\,}
\def\hJ{\widehat{J}\,}

\def\tG{\widetilde{\Gamma}}

\def\b{\beta}
\def\ga{\gamma}

\def\bx{\overline{x}}
\def\bb{\mathbf b}

\def\bB{\mathbf{B}}
\def\bi{\mathbf{i}}

\def\bd{\mathbf{d}}

\def\vpi{\varpi}
\def\cN{\mathcal{N}}
\def\bD{\widetilde{D}}
\def\ds{\displaystyle}
\def\proof{\medskip\noindent {\it Proof --- \ }}
\def\cqfd{\hfill $\Box$ \bigskip}
\def\H{\mathcal{H}}
\def\bchi{\widetilde{\chi}}
\def\hQ{\widehat{Q}}
\def\hDe{\widehat{\De}}
\def\Ga{\Gamma}
\def\tx{\widetilde{x}}

\def\tL{\widetilde{\Lambda}}
\def\tG{\widetilde{\Gamma}}
\def\veps{\varepsilon}
\def\I{\mathcal{I}}

\begin{document}

\title{\bf Quantum Grothendieck rings and derived Hall algebras}
\author{D. Hernandez, B. Leclerc}

\date{}

\maketitle

\begin{abstract}
Let $\g = \n\oplus\h\oplus\n_-$ be a simple Lie algebra over $\C$ of type $A, D, E$, and let
$U_q(L\g)$ be the associated quantum loop algebra.
Following Nakajima \cite{N}, Varagnolo-Vasserot \cite{VV}, and the
first author \cite{H1},
we study a $t$-deformation $\K_t$ of the Grothendieck ring of
a tensor category $\CC_\Z$ of finite-dimensional $U_q(L\g)$-modules. 
We obtain a presentation of $\K_t$ by generators and relations. 

Let $Q$ be a Dynkin quiver of the same type as $\g$.
Let $DH(Q)$ be the derived Hall algebra of the bounded derived category
$D^b(\md(FQ))$ over a finite field $F$, 
introduced by To\"en \cite{T}. Our presentation shows that
the specialization of $\K_t$ at $t=\sqrt{|F|}$ is isomorphic to 
$DH(Q)$. Under this isomorphism, the classes of fundamental $U_q(L\g)$-modules 
are mapped to scalar multiples of the classes of indecomposable objects in $DH(Q)$.

Our presentation of $\K_t$ is deduced from the preliminary study of a tensor subcategory
$\CC_{Q}$ of $\CC_\Z$ analogous to the heart $\mod(FQ)$ of the 
triangulated category $D^b(\md(FQ))$. We show that the $t$-deformed
Grothendieck ring $\K_{t,\,Q}$ of $\CC_{Q}$ is isomorphic to the positive
part of the quantum enveloping algebra of $\g$, and that the basis 
of classes of simple objects of $\K_{t,\,Q}$ corresponds to the dual of 
Lusztig's canonical basis. The proof relies on the algebraic characterizations
of these bases, but we also give a geometric approach in the last section.

It follows that for every orientation $Q$ of the Dynkin diagram,
the category $\CC_{Q}$ gives a new categorification of the 
coordinate ring $\C[N]$ of a unipotent group $N$ 
with Lie algebra $\n$, 
together with its dual canonical basis.

\end{abstract}

\bigskip
\setcounter{tocdepth}{1}
\tableofcontents


\section{Introduction}\label{sect1}

{\bf 1.1}\quad Let $\g$ be a simple Lie algebra of type $A, D, E$ over $\C$.
We denote by $\g = \n\oplus \h\oplus\n_-$ a triangular decomposition of $\g$.
Let $v$ be an indeterminate, and let 
\[
U_v(\g) = U_v(\n)\otimes U_v(\h)\otimes U_v(\n_-)
\]
be the corresponding Drinfeld-Jimbo quantum enveloping algebra over $\C(v)$,
defined via a $v$-analogue of the Chevalley-Serre presentation of $U(\g)$.
Using a geometric realization of $U_v(\n)$ in terms of perverse sheaves on
varieties of representations of a quiver $Q$ of the same Dynkin type as $\g$, 
Lusztig \cite{Lu} has defined a canonical basis $\bB$ of $U_v(\n)$ with favorable positivity 
properties. This was inspired by a seminal work of Ringel \cite{Ri}, showing that 
the twisted Hall algebra of the category $\mod(FQ)$ of representations of $Q$ over a finite field $F$,
is isomorphic to the specialization of $U_v(\n)$ at $v=\sqrt{|F|}$.

\bigskip\noindent
{\bf 1.2}\quad One can associate with $\g$ another quantum algebra. 
Let $L\g = \C[t,t^{-1}]\otimes \g$ be the loop algebra of $\g$.
Let $q$ be a nonzero complex number,
which is not a root of unity.
Via a $q$-analogue of the loop presentation of 
$U(L\g)$, Drinfeld \cite{D} has defined the quantum loop algebra $U_q(L\g)$,
an algebra over $\C$. 
The finite-dimensional representations of $U_q(L\g)$ have attracted a lot
of attention, because of their connection with the trigonometric solutions 
of the quantum Yang-Baxter equation with spectral parameter.
In this paper we focus on a tensor subcategory $\CC_\Z$ of the category
of finite-dimensional $U_q(L\g)$-modules, whose simple objects are parametrized
by a discrete set (for the precise definition of $\CC_\Z$ see \cite[\S 3.7]{HL} or
\S\ref{subsectCZ} below).  
Denote by $\R$ the complexified Grothendieck ring of $\CC_\Z$.
Let $t$ be another indeterminate. 
By works of Nakajima \cite{N} and Varagnolo-Vasserot \cite{VV}, the $\C$-algebra $\R$ has an interesting
$t$-deformation $\R_t$ over $\C(t)$.
The first author \cite{H1} has introduced a slightly different 
deformation $\K_t$.
These $t$-deformations are important because they contain for 
every simple object $L$ of $\CC_\Z$ a ``class'' $[L]_t$
which can be characterized by axioms similar to those of Lusztig for the canonical
basis $\bB$. As a consequence, Nakajima \cite{N} has shown that one can  
calculate algorithmically the character of $L$.

\bigskip\noindent
{\bf 1.3}\quad
Surprisingly, these deformed Grothendieck rings have not been much studied from
the ring theoretic point of view,
and for instance, to the best of our knowledge, there is no available presentation
by generators and relations in the literature. One of the main results of this paper (Theorem~\ref{presentation})
is a presentation of $\K_t$, with a similar flavor as the familiar Drinfeld-Jimbo
presentation of $U_v(\n)$. More precisely, this presentation shows that
$\K_t$ is obtained by taking an infinite number of copies of $U_t(\n)$
labelled by $m\in\Z$, and then imposing $t$-boson relations between
generators of copies sitting at adjacent integers, and $t$-commutation relations
between generators of non-adjacent copies.  

\bigskip\noindent
{\bf 1.4}\quad
Let $D^b(\md(FQ))$ be the bounded derived category of $\md(FQ)$.
To\"en \cite{T} has attached to this triangulated category an associative
algebra called the derived Hall algebra of 
$D^b(\md(FQ))$ (see also \cite{XX}). 
Let $DH(Q)$ denote the twisted derived Hall algebra obtained by twisting
To\"en's multiplication by means of the Ringel form, as in \cite{S}. 
It follows from our presentation of $\K_t$ that:
\begin{Thm}\label{Th_DHall}
\begin{itemize}
\item[(a)] The specialization of $\K_t$ at $t=\sqrt{|F|}$ is isomorphic to 
$DH(Q)$. 

\item[(b)]Under this isomorphism, the classes of fundamental $U_q(L\g)$-modules 
are mapped to scalar multiples of the classes of indecomposable stalk complexes in $DH(Q)$,
and the basis of classes of standard $U_q(L\g)$-modules is mapped to a
rescaling of the natural basis of $DH(Q)$ indexed by isoclasses of 
objects of $D^b(\md(FQ))$.
\end{itemize}
\end{Thm}
There is a similar result for the $t$-deformed Grothendieck ring $\R_t$
of \cite{N,VV}, but the twisted derived Hall algebra should be replaced
by a non-twisted one (Remark~\ref{Remark_untwistedDH}).

\bigskip\noindent
{\bf 1.5}\quad
To obtain our presentation of $\K_t$ we first consider a tensor subcategory
$\CC_{Q}$ of $\CC_\Z$ which ``looks like $\md(FQ)$ inside $D^b(\md(FQ))$''.
Recall that in \cite{HL} we have introduced an increasing sequence $(\CC_\ell)_{\ell>0}$
of subcategories of $\CC_\Z$.
When $Q$ is a bipartite orientation of the Dynkin diagram and 
the Coxeter number $h$ is even, 
$\CC_Q$ is just the subcategory $\CC_\ell$ with $\ell=h/2-1$.
The general definition of $\CC_Q$ for an arbitrary orientation $Q$
will be given in \S\ref{sectCh'} below.
Let $\K_{t,\,Q}$ be the subalgebra of $\K_t$ spanned by the elements
$[L]_t$ associated with the simple objects $L$ of~$\CC_Q$.
Note that $\K_t$ and $\K_{t,\,Q}$ are algebras over $\C(t^{1/2})$,
where $t^{1/2}$ is a square root of $t$. 

The quantum algebra $U_v(\n)$ is endowed with a distinguished scalar product.
Let $\bB^*$ be the basis of $U_v(\n)$ adjoint to the canonical basis $\bB$ with
respect to this scalar product.
Let~$v^{1/2}$ be a square root of $v$, and set 
$\U_{v}(\n):=\C(v^{1/2})\otimes U_v(\n)$. The main step for obtaining the presentation
of $\K_t$ is:
\begin{Thm}\label{mainTh}
\begin{itemize}
\item[(a)] There is a $\C$-algebra isomorphism $\Phi \colon \K_{t,\,Q} \overset{\sim}{\to} \U_{v}(\n)$
with $\Phi(t^{1/2}) = v^{1/2}$.
\item[(b)] For every simple object $L$ of $\CC_{Q}$, the image $\Phi([L]_t)$
belongs to $\bB^*$ (up to some half-integral power of $v$).
\end{itemize}
\end{Thm}

Nakajima obtained in \cite{Ncl} similar results for the first subcategory $\CC_1$ of \cite{HL}.
Namely, he showed that the classes $[L]_t$ of simple objects
of $\CC_1$ can be identified with a subset of the basis $\widetilde{\bB}^*$
of $U_v(\widetilde{\n})$.
Here~$\widetilde{\n}$ denotes the  
positive part of the Kac-Moody algebra
of rank $2\rk(\g)$ attached to the decorated Dynkin diagram of $\g$,
and $\widetilde{\bB}$ is Lusztig's canonical basis of $U_v(\widetilde{\n})$.
For example, if $\g$ has type $A_3$, $\widetilde{\g}$ has type $E_6$.
Note that in Theorem~\ref{mainTh}, we do not use $\widetilde{\n}$, but only $\n$.

\bigskip\noindent
{\bf 1.6}\quad
Let $A_v(\n)$ be the graded dual of the vector space $U_v(\n)$. It can be 
endowed with a multiplication coming from the comultiplication of $U_v(\g)$,
and regarded as the quantum
coordinate ring of the unipotent group $N$ with Lie algebra $\n$ (see \eg \cite{GLS}).
The basis $\bB^*$ can be identified with a basis of $A_v(\n)$ called the dual canonical basis. 
It specializes when $v\mapsto 1$ to a basis $\B$ of the coordinate ring $\C[N]$.

By specializing $v^{1/2}$ and $t^{1/2}$ to $1$ in Theorem~\ref{mainTh}, 
we see that the complexified Grothendieck
ring $\R_{Q}$ of $\CC_{Q}$ can be identified with the coordinate ring 
$\C[N]$ in such a way that the basis of $\R_{Q}$ consisting of the classes
of simple objects becomes Lusztig's dual canonical basis $\B$ of $\C[N]$.
We can therefore state:
\begin{Thm}
The tensor category $\CC_{Q}$ is a categorification of the ring $\C[N]$
and its dual canonical basis $\B$. 
\end{Thm}

Note that, by work of Khovanov-Lauda \cite{KL}, Rouquier \cite{R}, and Varagnolo-Vasserot \cite{VV2},
$(\C[N], \B)$ has another categorification in terms of KLR-algebras. 
In type $A_n$, KLR-algebras are isomorphic to blocks of affine Hecke algebras,
and the category $\CC_Q$ for an equi-oriented quiver~$Q$ is related to
a category of representations of affine Hecke algebras 
through the quantum affine Schur-Weyl duality. 
It would be interesting to find for other Dynkin quivers $Q$
similar functors between $\CC_{Q}$ and the module categories of the corresponding KLR-algebras. 

\bigskip\noindent
{\bf 1.7}\quad
The first author \cite{H3} has shown that tensor products of simple objects
of $\CC_\Z$ have the following remarkable property: a tensor product
$L_1\otimes \cdots \otimes L_k$ of simple objects $L_i$ is simple if and only if for every
pair $1\le i < j \le k$ the tensor product $L_i\otimes L_j$ is simple.
Using Theorem~\ref{mainTh} this yields the following: 
\begin{Cor}\label{corollary_binary}
A product $b_1\cdots b_k$ of elements $b_i$ of the dual canonical basis 
$\bB^*$ of $U_v(\n)$ belongs to $\bB^*$ up to a power of $v$ if and
only if for every
pair $1\le i < j \le k$ the product $b_i b_j$ belongs to $\bB^*$ up to a power of $v$. 
\end{Cor}
Corollary~\ref{corollary_binary} was expected in relation with the program of
Berenstein and Zelevinsky \cite{BZ, BZ2} of describing $\bB^*$ in terms of quantum 
cluster algebras.
But it was only known in a few low rank cases.

\bigskip\noindent
{\bf 1.8}\quad
Theorem~\ref{mainTh} also gives new supporting evidence for some conjectures 
formulated in \cite{GLS} and \cite{HL}.
It was conjectured in \cite[\S13]{HL} that for every $\ell \in \N$,
the Grothendieck ring $\R_\ell$ of $\CC_{\ell}$
has a particular cluster algebra structure for which
all cluster monomials are classes of simple objects. 
In \cite{GLS}, it is shown that $U_v(\n)$ has a quantum cluster algebra
structure, and it is conjectured that all quantum cluster monomials
belong to $\bB^*$.
Suppose that $Q$ is bipartite and $h$ is even. Set $h'=h/2-1$.
By comparing initial seeds, one sees that 
the quantum cluster structure of $\K_{t,\,Q}$ obtained by transporting
via $\Phi^{-1}$ the quantum cluster structure of
$U_v(\n)$ is a $t$-analogue of the cluster structure of $\R_{h'}$ conjectured
in \cite{HL}. 
Thus, by Theorem~\ref{mainTh}, the two conjectures of \cite{GLS} for $U_v(\n)$ and of \cite{HL}
for $\R_{h'}$ are essentially equivalent.

In \cite{HL} and \cite{Ncl}, the conjecture for $\R_\ell$ was proved in the
first non trivial case $\ell = 1$.
(In \cite{HL} some combinatorial steps of the proof were only verified 
for $\g$ of type $A_n$ and $D_4$; the proof of \cite{Ncl} is general and 
uses geometric representation theory.)
Since $\CC_1$ is a tensor subcategory of $\CC_Q$ (for every $\g$ except $\Sl_2$ and $\Sl_3$),
$\K_{t,\,Q}$ contains a subring $\K_{t,1}$ corresponding to $\CC_1$.
It is easy to see that $\Phi(\K_{t,1})$ is equal to the subalgebra
$\U_{v}(\n(w))$ of \cite{GLS} where $w = c^2$ is the square of the Coxeter
element of the Weyl group of $\g$ corresponding to the bipartite quiver $Q$.
This is a quantum cluster algebra 
of finite cluster type, equal to the Dynkin type of $\g$ in the classification
of Fomin and Zelevinsky. Thus, using \cite{HL,Ncl,Qin}, Theorem~\ref{mainTh} readily implies:
\begin{Cor}\label{Cor_cluster}
Let $w = c^2$ be as above. Then $\bB^* \cap \,U_v(\n(w))$ is equal
to the set of quantum cluster monomials of $U_v(\n(w))$. 
\end{Cor}

For $\g$ of type $A_n$, Lampe \cite{La} has given a direct proof 
of the fact that 
the quantum cluster \emph{variables} of $U_v(\n(w))$ belong to $\bB^*$.

\bigskip\noindent
{\bf 1.9}\quad
Since the bases $\bB^*$ and $\{[L]_t\}$ have geometric origin,
it is natural to ask for a geometric explanation of Theorem~\ref{mainTh}~(b).
In the final part of the paper, we show (Theorem~\ref{thmgeo})
that the quiver representation spaces $E_\bd$ 
used by Lusztig to define the canonical basis of $U_v(\n)$ 
are isomorphic to some 
particular graded quiver varieties $\M^\bullet_0(W^\bd)$
used by Nakajima for describing the classes $[L]_t$ of the simple objects $L$ of $\CC_Q$.
Moreover the intersection cohomology sheaves of closures of $G_\bd$-orbits in $E_\bd$
can be identified with the intersection cohomology sheaves of closures of strata
in $\M^\bullet_0(W^\bd)$.  
This is inspired by a similar result of Nakajima \cite{Ncl} for the 
category~$\CC_1$. 

\bigskip\noindent
{\bf 1.10}\quad
We now give an overview of the structure
of the paper. 
In Section~\ref{sectCartanAuslander}, we set up our notation and 
introduce an important bijection $\varphi$ between the set of
fundamental modules of $\CC_\Z$ and the vertices of the 
Auslander-Reiten quiver of $D^b(KQ)$. We use this bijection
to express the entries of 
the inverse of the quantum Cartan matrix of $\g$ 
in terms of the Ringel form of $Q$, or in terms of the scalar product of
the weight lattice of $\g$
(Proposition~\ref{formula_inverse}).
By construction, the quantum Grothendieck ring $\K_{t}$ is a
subring of a quantum torus $\Y_{t}$ over $\C(t^{1/2})$.
The $t$-commutation relations between
generators of $\Y_{t}$ are expressed in terms of entries of 
the inverse of the quantum Cartan matrix of $\g$ \cite{H1},
hence by Proposition~\ref{formula_inverse}, in terms of
scalar products of weights of $\g$.
The quantum Grothendieck ring $\K_{t,\,Q}$ is a
subring of a subtorus $\Y_{t,\,Q}$ of $\Y_{t}$,
of rank $r$ equal to the number of positive roots of $\g$. 
On the other hand, by \cite{GLS}, $\U_{v}(\n)$ has
an explicit embedding into a quantum torus $\T_{v,\,Q}$ of rank $r$ over $\C(v^{1/2})$, 
whose generators are certain unipotent quantum flag minors. 
The explicit $v$-commutation relations 
between these minors involve scalar products of roots and weights of $\g$.
Comparing these two presentations, 
we show that there is an isomorphism 
$\Phi\colon \Y_{t,\,Q} \to \T_{v,\,Q}$ mapping $t^{1/2}$ to $v^{1/2}$
(Proposition~\ref{main_prop}).

The proof that $\Phi$ restricts to an isomorphism from $\K_{t,\,Q}$ to $\U_{v}(\n)$
is based on some explicit systems of algebraic identities satisfied by the generators of both 
algebras. In Section~\ref{sect4}, we recall from \cite{GLS} a system of quantum determinantal 
identities occuring in $U_{v}(\n)$, and in Section~\ref{sect5} we derive a quantum 
$T$-system for the $(q,t)$-characters of the Kirillov-Reshetikhin modules.
(In~\cite[\S4]{N-KR}, a quantum $T$-system was already obtained for the 
$t$-deformed product used in \cite{N,VV}.
A~quantum cluster algebra related to the quantum $T$-system of type $A_1$ is also studied in \cite{DFK}.) 
Comparing these two systems we obtain that 
$\Phi$ maps the classes of the Kirillov-Reshetikhin modules of $\CC_{Q}$
to certain quantum minors of $\U_{v}(\n)$ (multiplied by explicit
half-integral powers of~$v$). 
In particular, $\Phi$ maps 
the classes of the fundamental modules of $\CC_{Q}$ in $\K_{t,\,Q}$
to the generators of the dual PBW-basis of $\U_{v}(\n)$ associated with
$Q$ (up to powers of $v^{1/2}$). 
This proves the first part of Theorem~\ref{mainTh}.
The second part is deduced from the algebraic characterizations
of $\bB^*$ and of the classes $[L]_t$  
(Section~\ref{sect6}). 
After some examples, we give the proof of Corollary~\ref{corollary_binary}.

The above-mentioned presentation of $\K_t$ (Theorem~\ref{presentation})
is deduced from Theorem~\ref{mainTh} in Section~\ref{sect_pres}, 
and in Section~\ref{sect_Hall} we prove the isomorphism
with the derived Hall algebra $DH(Q)$ stated in Theorem~\ref{Th_DHall}.
Finally, in Section~\ref{sect7}, we explain our geometric approach
to Theorem~\ref{mainTh}~(b) (Theorem~\ref{thmgeo}).

\section{Cartan matrices and Auslander-Reiten quivers}\label{sectCartanAuslander}

\subsection{Cartan matrix}\label{notation}
Let $\g$ be a simple Lie algebra of type $A, D, E$.
We denote by $I$ the set of vertices of its Dynkin diagram,
and we put $n = |I|$. 
The \emph{Cartan matrix} of $\g$ is the $I\times I$ matrix $C$
with entries
\[
 C_{ij} = 
\left\{
\begin{array}{cl}
2 & \mbox{if } i=j,\cr
-1 & \mbox{if $i$ and $j$ are adjacent vertices of the Dynkin diagram},\cr
0 &\mbox{otherwise.}
\end{array}
\right. 
\]
We shall often use the shorthand
notation $i\sim j$ to say that $C_{ij}=-1$.

We denote by $P$ the weight lattice  of $\g$, and by 
$\vpi_i\ (i\in I)$ its basis of fundamental weights.
The simple roots are defined by
\[
\a_i = \sum_{j\in I} C_{ij}\vpi_j,\qquad (i\in I). 
\]
The set of simple roots is denoted by $\Pi := \{\a_i\mid i\in I\}$.
We denote by $(\cdot, \cdot)$ the scalar product of $P$
defined by $(\a_i,\vpi_j) = \de_{ij}$.
Equivalently $(\a_i,\a_j) = C_{ij}$.
The Weyl group $W$ is generated by the reflexions $s_i$
acting on $P$ by
\[
s_i(\la) = \la -(\la,\a_i)\a_i,\qquad (\la\in P,\ i\in I). 
\]
The root system of $\g$ is $\De:= W\Pi$. It decomposes as
$\De = \De_+ \sqcup \De_-$, where $\De_+ = \De \cap (\oplus_{i\in I} \N \a_i)$
and $\De_- = -\De_+$. We write $r:=|\De_+|$.

A \emph{Coxeter element} of $W$ is a product of the form 
$c = s_{i_1}\cdots s_{i_n}$ where $(i_1,\ldots,i_n)$ is an
arbitrary ordering of $I$. All Coxeter elements are conjugate
in $W$. Their common order is called the \emph{Coxeter number} and
denoted by $h$. We have $hn = 2r$.
 
\subsection{Quivers}\label{subsect_quivers}

Let $Q$ be an orientation of the Dynkin diagram of $\g$. 
In other words, $Q$ is a Dynkin quiver of the same Dynkin type as $\g$.

For $i\in I$, we denote by $s_i(Q)$ the quiver obtained from
$Q$ by changing the orientation of every arrow with source $i$
or target $i$. Let $w = s_{i_1}\cdots s_{i_k}\in W$ be a reduced
decomposition. We say that $\bi=(i_1,\ldots,i_k)$  
is \emph{adapted} to $Q$ if $i_1$ is a source of $Q$,
$i_{2}$ is a source of $s_{i_1}(Q)$, \ldots, $i_k$ is a source
of $s_{i_{k-1}}\cdots s_{i_1}(Q)$.
There is a unique Coxeter element having reduced expressions adapted to $Q$.
We shall denote it by $\tau$.

\begin{figure}[t]
\[
\xymatrix@-1.0pc{
&&&&4\ar[ld]
\\
&&& \ar[lld]\ar[ld]3 &
\\
&1 &2&&
}
\qquad\qquad
\xymatrix@-1.0pc{
&{}\save[]+<0cm,2ex>*{\vdots}\restore&{}\save[]+<0cm,2ex>*{\vdots}\restore& (3,3) &{}\save[]+<0cm,2ex>*{\vdots}\restore
\\
&(1,2) \ar[rru]&(2,2)\ar[ru]&&\ar[lu](4,2)
\\
&&&\ar[llu]\ar[lu](3,1)\ar[ru]&
\\
&(1,0)\ar[rru] &(2,0)\ar[ru]&&\ar[lu](4,0)
\\
&&& \ar[llu]\ar[lu](3,-1)\ar[ru]&
\\
&(1,-2)\ar[rru] &(2,-2)\ar[ru]& & \ar[lu](4,-2) 
\\
& &&  \ar[llu]\ar[lu](3,-3)\ar[ru]& \\
&{}\save[]+<0cm,2ex>*{\vdots}\restore&{}\save[]+<0cm,2ex>*{\vdots}\restore&&{}\save[]+<0cm,2ex>*{\vdots}\restore
}
\]
\caption{\label{FigD4_1} {\it A quiver $Q$ of type $D_4$ and its repetition quiver $\hQ$}.}
\end{figure}
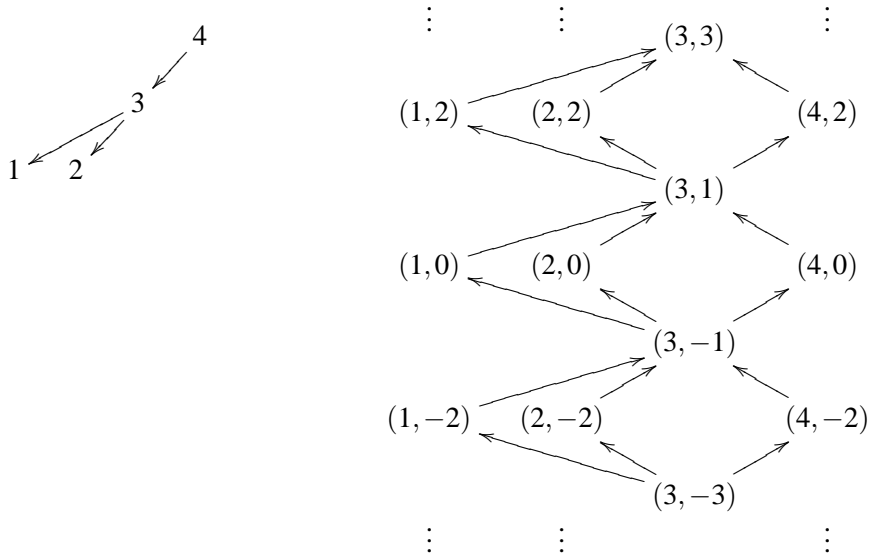

We denote by $Q_1$ the set of arrows of $Q$. 
A \emph{height function} $\xi\colon I \to \Z$ on $Q$ is a function
satisfying 
\[
\xi_j = \xi_i-1\quad \mbox{if}\quad  i\to j \in Q_1.
\]
Since $Q$ is connected, two height functions differ by a constant. 
We fix such a function $\xi$. Define
\[
\hI := \{(i,p)\in I\times\Z \mid p-\xi_i \in 2\Z\}.
\]
We attach to $Q$ the infinite \emph{repetition quiver} $\hQ$, defined as
the oriented graph with vertex set $\hI$ and two types of arrows:
\begin{itemize}
 \item[(i)] if there is an arrow $i\to j$ in $Q$ we have 
arrows $(i,p)\to (j,p+1)$ in $\hQ$ for all $(i,p)\in\hI$;
 \item[(ii)] if there is an arrow $i\to j$ in $Q$ we have 
arrows $(j,q)\to (i,q+1)$ in $\hQ$ for all $(j,q)\in\hI$.
\end{itemize} 
Note that $\hQ$ depends only on the Dynkin diagram, and not
on the choice of orientation $Q$. In fact, it is well known that
$\hQ$ is the quiver of a $\Z$-covering of the preprojective algebra
associated with~$Q$. In the literature, this quiver is often denoted
by $\Z Q$.
An example is shown in
Figure~\ref{FigD4_1}, where the height function is 
$\xi_1=\xi_2=0,\ \xi_3=1,\ \xi_4=2$.

Let $\hDe := \De_+\times\Z$. We now describe a natural labelling of 
the vertices of $\hQ$ by $\hDe$.
For $i\in I$, let $B(i)$ be the subset of $I$ consisting of all $j$'s
such that there is a path from $j$ to $i$ in $Q$.
Define
\[
 \ga_i := \sum_{j\in B(i)} \a_j,\qquad (i\in I).
\]
We have $\ga_i\in \De_+$. There is a unique bijection $\varphi\colon \hI \to \hDe$
defined inductively as follows:
\begin{itemize}
 \item[(a)] $\varphi(i,\xi_i) = (\ga_i,0)$ for $i\in I$;
 \item[(b)] suppose that $\varphi(i,p) = (\b,m)$; then 
\begin{itemize}
 \item[$\bullet$] $\varphi(i,p-2) = (\tau(\b),m)$ if $\tau(\b)\in \De_+$;
 \item[$\bullet$] $\varphi(i,p-2) = (-\tau(\b),m-1)$ if $\tau(\b)\in \De_-$;
 \item[$\bullet$] $\varphi(i,p+2) = (\tau^{-1}(\b),m)$ if $\tau^{-1}(\b)\in \De_+$;
 \item[$\bullet$] $\varphi(i,p+2) = (-\tau^{-1}(\b),m+1)$ if $\tau^{-1}(\b)\in \De_-$.
\end{itemize}
\end{itemize}
Note that this second labelling of $\hQ$ depends on $Q$.
This is illustrated in
Figure~\ref{FigD4_2}.
\begin{figure}[t]
\[
\xymatrix@-1.0pc{
&&&&4\ar[ld]
\\
&&& \ar[lld]\ar[ld]3 &
\\
&1 &2&&
}
\quad
\def\objectstyle{\scriptstyle}
\xymatrix@-1.0pc{
&{}\save[]+<0cm,2ex>*{\vdots}\restore&{}\save[]+<0cm,2ex>*{\vdots}\restore& (\a_1+\a_2+\a_3,1) &{}\save[]+<0cm,2ex>*{\vdots}\restore
\\
&(\a_1,1) \ar[rru]&(\a_2,1)\ar[ru]&&\ar[lu](\a_4,0)
\\
&&&\ar[llu]\ar[lu](\a_3+\a_4,0)\ar[ru]&
\\
&(\a_1+\a_3+\a_4,0)\ar[rru] &(\a_2+\a_3+\a_4,0)\ar[ru]&&\ar[lu](\a_3,0)
\\
&&& \ar[llu]\ar[lu](\a_1+\a_2+2\a_3+\a_4,0)\ar[ru]&
\\
&(\a_2+\a_3,0)\ar[rru] &(\a_1+\a_3,0)\ar[ru]& & \ar[lu](\a_1+\a_2+\a_3+\a_4,0) 
\\
& &&  \ar[llu]\ar[lu](\a_1+\a_2+\a_3,0)\ar[ru]& 
\\
&(\a_1,0)\ar[rru] &(\a_2,0)\ar[ru]& & \ar[lu](\a_4,-1) 
\\
& &&  \ar[llu]\ar[lu](\a_3+\a_4,-1)\ar[ru]& \\
&{}\save[]+<0cm,2ex>*{\vdots}\restore&{}\save[]+<0cm,2ex>*{\vdots}\restore&&{}\save[]+<0cm,2ex>*{\vdots}\restore
}
\]
\caption{\label{FigD4_2} {\it The labelling of $\hQ$ by $\hDe$ for $Q$ of type $D_4$}.}
\end{figure}

\subsection{Auslander-Reiten theory}\label{subsect_AR}

The quiver $\hQ$ with its labelling by $\hDe$ arises in the representation
theory of the path algebra $KQ$ of~$Q$ over a field $K$, as we shall
now recall.
We refer the reader to \cite{ARS,ASS,GaR,Ri0} for background 
on quiver representations and Auslander-Reiten theory.

Let $\md(KQ)$ be the abelian category of representations of $Q$ over $K$.
For an object $X$ of $\md(KQ)$ we write $\dimv(X)$ for its dimension vector.
We define the Ringel bilinear form
\[
\< X, Y\> := \dim(\Hom(X,Y)) - \dim(\Ext^1(X,Y)), \qquad (X,Y\in\md(KQ)),
\]
and the symmetric form $(X,Y):= \<X,Y\> + \<Y,X\>$.
It is known that these forms depend only on the dimension vectors 
$\dimv(X)$ and $\dimv(Y)$. Moreover, if we identify in the standard
way  $\dimv(X)$ and $\dimv(Y)$  with elements of the root lattice of $\g$,
then $(X,Y)$ coincides with the natural scalar product $(\dimv(X),\dimv(Y))$. 
In this picture, $\a_i$ is the
dimension vector of the simple $KQ$-module~$S_i$ supported on vertex $i$,
and $\ga_i$ is the dimension vector of its injective envelope $I_i$.
Recall that, by Gabriel's theorem, the isoclasses of indecomposable 
$KQ$-modules are in natural bijection with~$\De_+$.
They form the vertices of the Auslander-Reiten quiver $\Ga_Q$ of $\md(KQ)$.
The map $\b \mapsto (\b,0)$ identifies $\Ga_Q$ with the full subgraph of $\hQ$ with set
of vertices $\De_+\times \{0\}$. 
The map $\tau$ restricted to the dimension vectors in $\De_+$ of non projective $KQ$-modules 
is the Auslander-Reiten translation of $\md(KQ)$ \cite{ARS}.

Let $D^b(\md(KQ))$ be the bounded derived category of $KQ$.
Its indecomposable objects are the stalk complexes $X[i]$, 
consisting of an indecomposable object $X$ of $\md(KQ)$ sitting
in degree $i\in\Z$, and zero objects in all other degrees. 
Thus, the isoclasses of indecomposable objects of $D^b(\md(KQ))$
are naturally labelled by $\hDe$. 
Using this labelling,
the quiver $\hQ$ is identified with the Auslander-Reiten quiver of the triangulated
category $D^b(\md(KQ))$ \cite{Ha}.

\subsection{Quantum Cartan matrix}\label{cartan}

Let $z$ be an indeterminate, and
let $C(z)$ be the matrix with entries
\[
 C_{ij}(z) = 
\left\{
\begin{array}{cl}
z+z^{-1} & \mbox{if } i=j,\cr
-1 & \mbox{if $i\sim j$},\cr
0 &\mbox{otherwise.}
\end{array}
\right. 
\]
Thus $C(1)$ is just the Cartan matrix $C$ of $\g$.
Since $\det(C) \not = 0$, $\det(C(z))\not = 0$.  
We denote by $\widetilde{C}(z)$ the inverse of the matrix $C(z)$.
This is a matrix with entries $\widetilde{C}_{ij}(z)\in\Q(z)$.
Denoting by $A$ the adjacency matrix of the Dynkin diagram we have
\[
C(z) = (z+z^{-1})I - A, 
\]
therefore
\[
\widetilde{C}(z) = \sum_{k\ge 0} (z+z^{-1})^{-k-1}A^k. 
\]
Hence the entries of $\widetilde{C}(z)$ have power series expansions in $z$
of the form
\[
\widetilde{C}_{ij}(z) = \sum_{m \ge 1} \widetilde{C}_{ij}(m)\, z^{m},
\]
where $\widetilde{C}_{ij}(m)\in\Z$.
Note that since $C(z)$ is a symmetric matrix, we have 
$\widetilde{C}_{ij}(m) = \widetilde{C}_{ji}(m)$.

\subsection{Formula for $\widetilde{C}_{ij}(m)$}\label{formula_inverse}

We will now give several equivalent expressions for the coefficients $\widetilde{C}_{ij}(m)$.
For other expressions of $\widetilde{C}_{ij}(z)$ in type $A_n$ and $D_n$, see \cite[Appendix C]{FR1}. 

Fix an orientation $Q$ of the Dynkin diagram, and recall from \S\ref{subsect_quivers}
and \S\ref{subsect_AR}
the associated notation $\xi_i$, $\ga_i$, the Coxeter transformation $\tau$, and
the Ringel form $\<\cdot,\cdot\>$.
\begin{Prop}\label{prop1.1}
Let $m\ge 1$.
If $m+\xi_i-\xi_j-1$ is odd then $\widetilde{C}_{ij}(m) = 0$. Otherwise
\begin{equation}\label{eq1}
\widetilde{C}_{ij}(m) = \left(\tau^{(m+\xi_i-\xi_j-1)/2}(\ga_i), \vpi_j\right).  
\end{equation}
Equivalently, 
\begin{equation}\label{eq2}
\widetilde{C}_{ij}(m) = \left\< \tau^{(m+\xi_i-\xi_j-1)/2}(I_i), I_j\right\>.
\end{equation}
\end{Prop}
\proof
Let us denote temporarily by $D_{ij}(m)$ the value of $\widetilde{C}_{ij}(m)$
predicted by the proposition. We want to show that
\[
\sum_{k\in I,\ m\ge 1} C_{ik}(z) D_{kj}(m) z^m = \de_{ij},\qquad (i,j\in I). 
\]
Using the definition of $C_{ik}(z)$, 
this is equivalent to show that
\begin{equation}\label{eqtoprov}
\sum_{m\ge 1}\left(\left(z+z^{-1}\right)D_{ij}(m) - \sum_{k \sim i} D_{kj}(m)\right)z^m = \de_{ij},
\qquad (i,j\in I).
\end{equation}

The coefficient of $z^0$ in the left-hand side is equal to $D_{ij}(1)$.
If $\xi_i-\xi_j$ is odd then by definition $D_{ij}(1) = 0$. Otherwise,
if $\xi_i-\xi_j = 2l$, then $D_{ij}(1) = (\tau^l(\ga_i), \vpi_j)$ is the coefficient
of $\a_j$ in $\tau^l(\ga_i)$. Let $(\b,m)$ be the vertex of $\hQ$ in the column
of $(\ga_i,0)$ and at the same height as $(\ga_j,0)$.
Such a vertex exists because $\xi_i-\xi_j$ is even, and clearly $\b = \pm \tau^l(\ga_i)$.
Now it is a well-known fact from the combinatorics of Auslander-Reiten quivers that
for all vertices $(\ga, s)$ of $\hQ$ at the same height as $(\ga_j, 0)$ the 
coefficient of $\a_j$ in $\ga$ is 0, except if $(\ga, s) = (\ga_j,0)$ in which
case it is equal to 1. Hence we have $D_{ij}(1) = \de_{ij}$.

Consider now the coefficient of $z^m\ (m\ge 1)$ in (\ref{eqtoprov}). 
We need to show that
\begin{equation}\label{eqtoprov2}
D_{ij}(m+1)+D_{ij}(m-1) - \sum_{k:\,k \sim i} D_{kj}(m) = 0,
\qquad (i,j\in I, \ m\ge 1). 
\end{equation}
Note that for $k\sim i$ we have $\xi_k = \xi_i \pm 1$, hence
if $m+\xi_i-\xi_j$ is odd, all summands of the left-hand side are zero.
Otherwise, writing $m+\xi_i-\xi_j = 2l$, the left-hand side of
(\ref{eqtoprov2}) is 
\[
\left(\tau^l(\ga_i) + \tau^{l-1}(\ga_i) - \sum_{k\sim i} \tau^{l+(\xi_k-\xi_i+1)/2}(\ga_k)\ ,\ \vpi_j\right). 
\]
Now it is again a familiar fact from the combinatorics of Auslander-Reiten quivers that
\[
\tau^l(\ga_i) + \tau^{l-1}(\ga_i) =  \sum_{k\sim i} \tau^{l+(\xi_k-\xi_i+1)/2}(\ga_k),
\]
since the roots $\tau^{l}(\ga_i)$, $\tau^{l-1}(\ga_i)$, and $\tau^{l+(\xi_k-\xi_i+1)/2}(\ga_k)$
with $k\sim i$,
form a mesh. This proves (\ref{eq1}).

Finally, if $\b = \dimv X$ then $(\b,\vpi_j)$ is equal to the coefficient of
$\a_j$ in $\b$, hence 
\[
(\b,\vpi_j) = \dim(\Hom(X, I_j)) = \<X,\, I_j\>, 
\]
because $I_j$ is injective. This proves (\ref{eq2}). \cqfd

\begin{example}
{\rm
Take $\g$ of type $A_4$.
One has for instance
\[
\begin{array}{ccl}
\widetilde{C}_{11}(z)&=&z^{1}-z^{9}+z^{11}-z^{19}+\cdots \\[2mm]
\widetilde{C}_{12}(z)&=&z^{2}-z^{8}+z^{12}-z^{18}+\cdots \\[2mm]
\widetilde{C}_{13}(z)&=&z^{3}-z^{7}+z^{13}-z^{17}+\cdots \\[2mm]
\widetilde{C}_{14}(z)&=&z^{4}-z^{6}+z^{14}-z^{16}+\cdots \\[2mm]
\widetilde{C}_{21}(z)&=&z^{2}-z^{8}+z^{12}-z^{18}+\cdots \\[2mm]
\widetilde{C}_{22}(z)&=&z^{1}+z^{3}-z^{7}-z^{9}+z^{11}+z^{13}-z^{17}-z^{19}+\cdots \\[2mm]
\widetilde{C}_{23}(z)&=&z^{2}+z^{4}-z^{6}-z^{8}+z^{12}+z^{14}-z^{16}-z^{18}\cdots \\[2mm]
\widetilde{C}_{24}(z)&=&z^{3}-z^{7}+z^{13}-z^{17}+\cdots \\[2mm]
\end{array}
\]
Let us choose the sink-source orientation $Q$ with height function
$\xi_1=0, \xi_2=1, \xi_3=0, \xi_4=1$.
Then $\tau = s_2s_4s_1s_3$, and since $\tau^5 = 1$, the roots $\tau^l(\ga_i)$ are all determined by:
\[
\begin{array}{llll}
\ga_1 = \a_1+\a_2, & \ga_2 = \a_2, &\ga_3 = \a_2+\a_3+\a_4, &\ga_4 = \a_4,\\[2mm]
\tau(\ga_1) = \a_3+\a_4, & \tau(\ga_2) = \a_1+\a_2+\a_3+\a_4, &\tau(\ga_3) = \a_1+\a_2+\a_3, 
&\tau(\ga_4) = \a_2+\a_3,\\[2mm] 
\tau^2(\ga_1) = -\a_4, & \tau^2(\ga_2) = \a_3, &\tau^2(\ga_3) = -\a_2, &\tau^2(\ga_4) = \a_1,\\[2mm]
\tau^3(\ga_1) = -\a_2-\a_3, & \tau^3(\ga_2) = -\a_2-\a_3-\a_4, &\tau^3(\ga_3) = -\a_1-\a_2-\a_3-\a_4, 
&\tau^3(\ga_4) = -\a_1-\a_2,\\[2mm]
\tau^4(\ga_1) = -\a_1, & \tau^4(\ga_2) = -\a_1-\a_2-\a_3, &\tau^4(\ga_3) = -\a_3, &\tau^4(\ga_4) = -\a_3-\a_4.
\end{array} 
\]
For instance by Proposition~\ref{prop1.1}, $\widetilde{C}_{23}(6)$ is equal to the 
coefficient of $\a_3$ in $\tau^3(\ga_2) = -\a_2-\a_3-\a_4$, namely to $-1$.
} 
\end{example}

\begin{Cor}
For $i,j\in I$ and $m \ge 1$ we have
\[
\widetilde{C}_{ij}(m+2h) = \widetilde{C}_{ij}(m). 
\]
\end{Cor}
\proof
Since $\tau^h = 1$, this follows immediately from Proposition~\ref{prop1.1}.
\cqfd


\section{Quantum tori} \label{sect3}


\subsection{The quantum torus $\YY_t$}\label{secY}

Recall from \S\ref{subsect_quivers} the labelling set $\hI$ of $\hQ$.
Define 
\[
\Y := \C\left[Y_{i,p}^{\pm 1}\mid (i,p) \in \hI\right]
\] 
to be the Laurent polynomial ring
generated by a collection of commutative variables $Y_{i,p}$ labelled by $\hI$.
This ring is related to a tensor subcategory $\CC_\Z$ of the category of
finite-dimensional $U_q(L\g)$-modules considered in \cite{HL} (see below \S\ref{subsectCZ}).

Let $t$ be an indeterminate. 
Following \cite{H1} we introduce a $t$-deformed version $(\YY_t,*)$ of $\Y$,
with noncommutative multiplication denoted by $*$. 
This is the $\C(t)$-algebra generated by  
variables still denoted by $Y_{i,p}$, subject to the 
$t$-commutation relations
\begin{equation}\label{tcom}
Y_{i,p}*Y_{j,s} := t^{\cN(i,p;j,s)} Y_{j,s}*Y_{i,p},\qquad ((i,p),\,(j,s)\in \hI), 
\end{equation}
where
\begin{equation}\label{defcN}
\cN(i,p;j,s) := \widetilde{C}_{ij}(p-s-1) - \widetilde{C}_{ij}(p-s+1)
-\widetilde{C}_{ij}(s-p-1) + \widetilde{C}_{ij}(s-p+1).  
\end{equation}
Here we have extended the definition of $\widetilde{C}_{ij}(m)$ to
every $m\in\Z$ by setting $\widetilde{C}_{ij}(m) = 0$ if $m\le 0$.
Note that, since $\widetilde{C}(z)$ is symmetric, we have
\begin{equation}
\cN(i,p;j,s) = - \cN(j,s;i,p),\qquad (i,j\in I,\ p,s\in\Z). 
\end{equation}
If $p=s$ then $\cN(i,p;j,s) = 0$. 
Otherwise, without loss of generality 
we can assume that $p<s$. Then, (\ref{defcN}) simplifies as 
\begin{equation}\label{defcN2}
\cN(i,p;j,s) = 
\widetilde{C}_{ij}(s-p+1) - \widetilde{C}_{ij}(s-p-1),\qquad (p<s).  
\end{equation}
We regard the noncommutative ring $(\YY_t,*)$ as a quantum torus
of infinite rank.

\begin{remark}\label{Remark_differ}
{\rm
In \cite{VV} and \cite{N}, the construction of a $t$-deformed
Grothendieck ring is based on a slightly different quantum torus.
Namely, in these papers the product is defined by:
\begin{equation}\label{tcom'}
Y_{i,p}\,\cdot\,Y_{j,s} := t^{\cN'(i,p;j,s)} 
Y_{j,s}\,\cdot\,Y_{i,p},\qquad ((i,p),\,(j,s)\in \hI), 
\end{equation}
where instead of (\ref{defcN}) the following exponent is used:
\begin{equation}\label{defcN'}
\cN'(i,p;j,s) := -2\left(\widetilde{C}_{ij}(p-s-1) -\widetilde{C}_{ij}(s-p-1)\right). 
\end{equation}
For instance, in type $A_3$, we have
\[
Y_{1,0} * Y_{2,1} = t Y_{2,1} * Y_{1,0},
\qquad
\mbox{whereas}
\qquad
Y_{1,0} \,\cdot\, Y_{2,1} = Y_{2,1} \,\cdot\, Y_{1,0}. 
\]
In \cite{N,VV}, the definition of the product comes from 
a convolution operation for certain perverse sheaves on quiver varieties,
and the deformation parameter $t$ encodes the natural grading of complexes of sheaves.
Our product $*$ comes from \cite{H1} and  
the original construction of $q$-characters.
Indeed in \cite{FR}, the variables $Y_{i,p}\in\Y$ are defined as
formal power series in elements of $U_q(L\g)$, and they pairwise commute.
In \cite{H1}, these formal power series are replaced by certain infinite sums~$\tilde{Y}_{i,p}$
in elements of the quantum affine algebra $U_q(\widehat{\g})$ (with non trivial central charge~$c$), 
which can be seen as vertex operators. 
The original variable $Y_{i,p}$
is just one factor of the complete variable $\tilde{Y}_{i,p}$. The 
relations of the quantum affine algebra then give rise 
to $t$-commutation relations between the
$\tilde{Y}_{i,p}$, where the parameter $t$ appears as a formal power series with coefficients in
$\C[c^{\pm 1}]$ \cite[Theorem 3.11]{H1}. The defining relations (\ref{tcom}) (\ref{defcN}) of $*$ are obtained by replacing
$t$ by $t^{-1}$ in those $t$-commutation relations. 
}
\end{remark}

Recall from \S\ref{subsect_quivers} the bijection $\varphi\colon \hI \to \hDe$.
\begin{Prop}\label{expression_N}
Let $(i,p)$ and $(j,s)$ be elements of $\hI$ with $p<s$.
There holds
\[
\cN(i,p;j,s) = \left(\tau^{(s-p+\xi_i-\xi_j)/2}(\ga_i),\ \ga_j\right). 
\]
Moreover, if $\varphi(i,p) = (\b,m)$ and $\varphi(j,s) = (\de,l)$, then
\[
\cN(i,p;j,s) = (-1)^{l-m}\,(\b,\de). 
\]
\end{Prop}
\proof
First note that the definition of $\hI$ implies
that $s-p+\xi_i-\xi_j\in 2\Z$.
By Proposition~\ref{prop1.1}, we have
\begin{eqnarray*}
\cN(i,p;j,s) &=& \widetilde{C}_{ij}(s-p+1) - \widetilde{C}_{ij}(s-p-1)\\
&=& \left\<\tau^{(s-p+\xi_i-\xi_j)/2}(I_i),\ I_j\right\>
- \left\<\tau^{(s-p+\xi_i-\xi_j)/2-1}(I_i),\ I_j\right\>.
\end{eqnarray*}
Now recall the classical formula
\[
\<\tau^{-1}(X),Y\> = - \<Y,X\>,\qquad (X,Y\in\md(FQ)). 
\]
It follows that
\begin{eqnarray*}
\cN(i,p;j,s) &=&  \left\<\tau^{(s-p+\xi_i-\xi_j)/2}(I_i),\ I_j\right\> 
+ \left\<I_j, \tau^{(s-p+\xi_i-\xi_j)/2}(I_i)\right\> \\
&=& \left(\tau^{(s-p+\xi_i-\xi_j)/2}(I_i),\ I_j\right) \\
&=& \left(\tau^{(s-p+\xi_i-\xi_j)/2}(\ga_i),\ \ga_j\right).
\end{eqnarray*}
This proves the first equality.
The second equality is immediately deduced from the first one if we
note that, by definition of the bijection $\varphi$,
\[
\tau^{(\xi_i - p)/2}(\ga_i) = (-1)^m\b,
\qquad
\tau^{(\xi_j - s)/2}(\ga_i) = (-1)^l\de.
\]
\cqfd

\begin{remark}\label{Remark_differ2}
{\rm
For the product of \cite{VV} and \cite{N}, we have, for $p<s$, a similar expression
\[
\cN'(i,p;j,s) = 2\left\<\tau^{(s-p+\xi_i-\xi_j)/2-1}(\ga_i),\ \ga_j\right\>
=-2\left<\ga_j,\tau^{(s-p+\xi_i-\xi_j)/2}(\ga_i)\right\>
=(-1)^{l-m+1}2\<\de,\b\>,
\]
in which the symmetric scalar product $(\cdot,\cdot)$ is replaced by the
non-symmetric Ringel form $\<\cdot,\cdot\>$.
}
\end{remark}

\subsection{Commutative monomials}\label{commut_prod}
Let us adjoin a square root $t^{1/2}$ of $t$ and extend the quantum torus $(\YY_t,*)$
to
\[
(\Y_t,*):= \C(t^{1/2})\otimes_{\C(t)} (\YY_t,*).  
\]
We notice that the expression
\[
t^{\frac{1}{2}\cN(j,s;\,i,p)}\, Y_{i,p}*Y_{j,s} =
t^{\frac{1}{2}\cN(i,p;\,j,s)}\, Y_{j,s}*Y_{i,p} \in \Y_{t}
\]
is invariant under permutation of $(i,p)$ and $(j,s)$.
We can then denote it as a \emph{commutative} monomial $Y_{i,p}Y_{j,s}=Y_{j,s}Y_{i,p}$,
and write
\[
Y_{i,p}*Y_{j,s} = t^{\frac{1}{2}\cN(i,p;\,j,s)}\, Y_{i,p}Y_{j,s}. 
\]
More generally, for a family $(u_{i,p}\mid (i,p)\in \hI)$ of integers with finitely many nonzero
components, the expression
\[
t^{\frac{1}{2}\sum_{(i,p)<(j,s)} u_{i,p} u_{j,s} \cN(j,s;\,i,p)}
\hskip-0.4cm
\stackrel{\longrightarrow}{\mbox{\Large\quad *}_{\hskip-0.4cm (i,p)\in\hI}}
Y_{i,p}^{u_{i,p}} 
\]
does not depend on the chosen ordering of $\hI$ used to define it.
We will denote it as a commutative monomial $\prod_{(i,p)\in\hI} Y_{i,p}^{u_{i,p}}$,
and write
\[
\underset{(i,p)\in\hI}{\overset{\longrightarrow}{\mbox{\Large *}}} 
\,Y_{i,p}^{u_{i,p}} = 
t^{\frac{1}{2}\sum_{(i,p)<(j,s)} u_{i,p} u_{j,s} \cN(i,p;\,j,s)}
\prod_{(i,p)\in\hI} Y_{i,p}^{u_{i,p}}.
\] 
The commutative monomials form a basis of the $\C(t^{1/2})$-vector space
$\Y_{t}$.
It will be convenient to denote commutative monomials by
\[
m = \prod_{(i,p)\in \hI} Y_{i,p}^{u_{i,p}(m)}. 
\]
A commutative monomial $m$ is said to be \emph{dominant} if
$u_{i,p}(m) \ge 0$ for every $(i,p)\in \hI$.

The noncommutative product of two commutative monomials $m_1$ and $m_2$
is given by 
\begin{equation}
\label{mult_mono}
m_1 * m_2 = 
t^{\frac{1}{2}D(m_1,\,m_2)} m_1m_2
= t^{D(m_1,\,m_2)} m_2*m_1, 
\end{equation}
where 
\[
D(m_1,m_2) = \sum_{(i,p),(j,s)\in \hI} u_{i,p}(m_1)u_{j,s}(m_2)
\cN(i,p;\,j,s),
\]
and 
\[
m_1m_2 = \prod_{(i,p)\in \hI} Y_{i,p}^{u_{i,p}(m_1)+u_{i,p}(m_2)}, 
\]
denotes the commutative product.

\subsection{The quantum torus $\YY_{t,\,Q}$} \label{torusT}

Recall the bijection $\varphi \colon \hI \to \hDe$ of \S\ref{subsect_quivers}.
Define
\[
\hI_Q := \varphi^{-1}\left(\De_+\times \{0\}\right) \subset \hI, 
\]
and let $\YY_{t,\,Q}$ be the the $\C(t)$-subalgebra
of $(\YY_t,*)$ generated by the variables $Y_{i,p}\ ((i,p)\in \hI_Q)$.
This is a quantum torus of rank $r = |\De_+|$.
We will also use the extended torus 
\[
(\Y_{t,\,Q},*):= \C(t^{1/2})\otimes_{\C(t)} (\YY_{t,\,Q},*).
\]
\begin{example}\label{exampleYQ}
{\rm 
We take $\g$ of type $D_4$ and choose $Q$ as in Figure~\ref{FigD4_1} and
Figure~\ref{FigD4_2}. Comparing the two figures we see that $\YY_{t,\,Q}$
is generated by
\[
Y_{1,0}^{\pm1},\ Y_{1,-2}^{\pm1},\ Y_{1,-4}^{\pm1},\
Y_{2,0}^{\pm1},\ Y_{2,-2}^{\pm1},\ Y_{2,-4}^{\pm1},\
Y_{3,1}^{\pm1},\ Y_{3,-1}^{\pm1},\ Y_{3,-3}^{\pm1},\
Y_{4,2}^{\pm1},\ Y_{4,0}^{\pm1},\ Y_{4,-2}^{\pm1}. 
\]
}
\end{example}

\subsection{The quantum torus $\TT_{v,\,Q}$}\label{qtorusA}

Let $w_0$ be the longest element of $W$. Let $\bi = (i_1,\ldots,i_r)$
be a reduced expression of $w_0$ adapted to~$Q$ (see \S\ref{subsect_quivers}).
Following \cite[\S11]{GLS}, we introduce a quantum torus $\TT_{v,\,Q}$ of rank $r$ over $\C(v)$.
(The indeterminate $v$ is denoted by $q$ in \cite{GLS}). 
Its generators are certain unipotent quantum minors
\[
D_{\vpi_{i_k},\,\la_k},\qquad (1\le k\le r) 
\]
in the quantum coordinate ring $A_v(\n)$.
Here $\la_k$ is the weight given by
\[
\la_k = s_{i_1}\cdots s_{i_k}(\vpi_{i_k}),\qquad (1\le k\le r). 
\]
The definition of $A_v(\n)$ will be recalled in \S\ref{defAqn} below.
At this stage we only need to know the explicit $v$-commutation 
relations satisfied by these minors. 
It is shown in \cite[Lemma 11.2]{GLS} that for $k<l$ there holds
\begin{equation}\label{commutD}
D_{\vpi_{i_k},\,\la_k} D_{\vpi_{i_l},\,\la_l} = 
v^{(\vpi_{i_k}-\la_k,\ \vpi_{i_l}+\la_l)}
D_{\vpi_{i_l},\,\la_l} D_{\vpi_{i_k},\,\la_k}. 
\end{equation}
 
For $1\le k\le r$, set $k^-:= \max(\{s<k \mid i_s = i_k\} \cup \{0\}).$
Define 
\begin{equation}\label{defEk}
Z_k :=
D_{\vpi_{i_{k}},\,\la_{k}} \left(D_{\vpi_{i_{k^-}},\,\la_{k^-}}\right)^{-1}, 
\end{equation}
where if $k^-=0$ we understand $D_{\vpi_{i_{k^-}},\,\la_{k^-}} = 1$.
Clearly, $Z_k\ (1\le k\le r)$ is another set of generators of $\TT_{v,\,Q}$.
Let
\begin{equation}\label{eq_beta}
\b_k = s_{i_1}\cdots s_{i_{k-1}}(\a_{i_k}),\qquad (1\le k\le r). 
\end{equation}
Note that we have
\begin{equation}\label{relation-n}
\la_{k} = \la_{k^-} - \b_{k}, \qquad (1\le k \le r), 
\end{equation}
where if $k^-=0$ we use the convention $\la_{k^-}=\varpi_{i_k}$.

\begin{Prop}\label{lemcruc}
For $1\le k<l \le r$, we have:
\begin{equation}\label{toprove}
Z_k Z_l = v^{-(\b_k,\b_l)} Z_l Z_k. 
\end{equation}
\end{Prop}
\proof
Let us introduce the integers $\mu_{kl}$ and $\nu_{kl}$ such that
\[
D_{\vpi_{i_k},\la_k}D_{\vpi_{i_l},\la_l} = v^{\mu_{kl}} D_{\vpi_{i_l},\la_l} D_{\vpi_{i_k},\la_k},
\qquad
Z_k Z_l = v^{\nu_{kl}} Z_l Z_k,
\qquad (1\le k,\ l\le r).
\]
By definition of $Z_k$ we have
\[
\nu_{kl} = (\mu_{kl} - \mu_{k^-l}) - (\mu_{kl^-} - \mu_{k^-l^-}), 
\]
where we use the convention that $\mu_{k^-l} = 0$ if $k^-=0$, 
$\mu_{kl^-} = 0$ if $l^-=0$, and $\mu_{k^-l^-} = 0$ if $k^-=0$ or $l^-=0$.
Since $k^-<k<l$, we have
\begin{equation}
\mu_{kl} - \mu_{k^-l} = (\vpi_{i_k}-\la_{k},\  \vpi_{i_l}+\la_{l}) 
- (\vpi_{i_k}-\la_{k^-},\  \vpi_{i_l}+\la_{l})
= (\b_k, \vpi_{i_l}+\la_{l}).  
\end{equation}

(a)\quad If $k<l^-$ we have similarly 
\[
\mu_{kl^-} - \mu_{k^-l^-} = ((\b_k, \vpi_{i_l}+\la_{l^-}) 
\]
and so
\[
\nu_{kl} = (\b_k, \vpi_{i_l}+\la_{l}) -  (\b_k, \vpi_{i_l}+\la_{l^-}) = -(\b_k,\b_l),
\]
as required.

(b)\quad If $k=l^-$ then $\mu_{kl^-}=0$ and 
$\mu_{k^-l}= (\vpi_{i_k}-\la_{k^-},\ \vpi_{i_k}+\la_{k})$.
Hence
\begin{eqnarray*}
\nu_{kl}& = &(\b_k, \vpi_{i_l}+\la_{l}) +  (\vpi_{i_k}-\la_{k^-},\ \vpi_{i_k}+\la_{k})\\
        & = &(\b_k, \vpi_{i_k}+\la_{k}) - (\b_k,\b_l) +  (\vpi_{i_k}-\la_{k^-},\ \vpi_{i_k}+\la_{k})\\
       & = & - (\b_k,\b_l) +  (\vpi_{i_k}-\la_{k^-}+\b_k,\ \vpi_{i_k}+\la_{k})\\
& = & - (\b_k,\b_l) +  (\vpi_{i_k}-\la_{k},\ \vpi_{i_k}+\la_{k})\\
& = & - (\b_k,\b_l)
\end{eqnarray*}
as required, because $(\vpi_{i_k}-\la_{k},\ \vpi_{i_k}+\la_{k})=(\vpi_{i_k},\vpi_{i_k})-(\la_k,\la_k) = 0$.

(c)\quad If $k>l^-$ then $\mu_{kl^-} = - \mu_{l^-k} = - (\vpi_{i_l}-\la_{l^-},\ \vpi_{i_k}+\la_{k})$.
Hence
\begin{eqnarray*}
-\mu_{kl^-}+\mu_{k^-l^-}& = &(\vpi_{i_l}-\la_{l^-},\ \vpi_{i_k}+\la_{k}) + 
(\vpi_{i_k}-\la_{k^-},\ \vpi_{i_l}+\la_{l^-})\\
&=& -(\b_k, \ \vpi_{i_l}+\la_{l^-}) + (\vpi_{i_l}-\la_{l^-},\ \vpi_{i_k}+\la_{k})
+(\vpi_{i_k}-\la_{k},\ \vpi_{i_l}+\la_{l^-})\\
&=& -(\b_k, \ \vpi_{i_l}+\la_{l^-}) + 2(\vpi_{i_k},\vpi_{i_l}) -2(\la_k,\la_{l^-})\\
&=& -(\b_k, \ \vpi_{i_l}+\la_{l^-}).
\end{eqnarray*}
Indeed, since $l^-<k<l$, we have
\[
(\la_k,\la_{l^-}) = (s_{i_{l^-+1}}\cdots s_{i_k}(\vpi_{i_k}),\ \vpi_{i_l})
= (\vpi_{i_k}, s_{i_{k}}\cdots s_{i_{l^-+1}}(\vpi_{i_l})) = (\vpi_{i_k},\vpi_{i_l}).
\]
It follows that again, 
$\nu_{kl} = (\b_k, \vpi_{i_l}+\la_{l}) -  (\b_k, \vpi_{i_l}+\la_{l^-}) = -(\b_k,\b_l)$,
as required.
\cqfd

\subsection{An isomorphism}\label{def_tilde_k}

It is well known that the roots $\b_k\ (1\le k\le r)$ give an
enumeration of $\De_+$. Therefore, for every $(i,p)\in\hI_Q$
there is a unique $k$ such that $\varphi(i,p) = (\b_k,0)$.

\begin{Prop}\label{prop2.2}
The assignment 
\[
\begin{array}{ll}
t \mapsto v,& \\[3mm]
Y_{i,p} \mapsto Z_k,  &\mbox{where } (i,p)\in \hI_Q,\ \mbox{and } \varphi(i,p) = (\b_k,0),
\end{array}
\]
extends to an isomorphism of quantum tori from $\YY_{t,\,Q}$ to $\TT_{v,\,Q}$.
\end{Prop}
\proof
This follows immediately from Proposition~\ref{expression_N} and
Proposition~\ref{lemcruc} if we note that when 
$\varphi(i,p) = (\b_k,0)$ and $\varphi(j,s) = (\b_l,0)$,
$p<s$ implies that $k>l$.
\cqfd

\subsection{The involution $\si$ and the rescaled generators $X_k$}\label{sigma}

Let $\T_{v,\,Q}:= \C(v^{1/2})\otimes_{\C(v)} \TT_{v,\,Q}$.
For $\gamma = \sum_{i}c_i\a_i$ in the root lattice of $\g$, we set
\begin{equation}\label{defNgamma} 
\deg\ga := \sum_i c_i,\qquad N(\ga) := \frac{(\ga,\ga)}{2} - \deg\ga. 
\end{equation}
Following \cite{GLS}, we introduce an involution $\si$ of $\T_{v,\,Q}$, defined as 
the $\C$-algebra anti-automor\-phism satisfying
\begin{equation}
\si(v^{1/2}) = v^{-1/2},\qquad 
\si\left(D_{\vpi_{i_k},\la_k}\right) = v^{N(\vpi_{i_k}-\la_k)} D_{\vpi_{i_k},\la_k}. 
\end{equation}

We rescale the generators $Z_k$ of $\TT_{v,\,Q}$ by defining
\begin{equation}
X_k 
:=
\left\{
\begin{array}{cl}
v^{N(\b_k)/2} Z_k & \mbox{if $1\le k\le n$,}\\[2mm]
v^{N(\b_k)/2+(\vpi_{i_k}-\la_{k-n},\ \b_k)} Z_k & \mbox{if $n+1\le k\le r$}.
\end{array}
\right.  
\end{equation}
Note that these elements live in $\T_{v,\,Q}$.

\begin{Lem} 
For $1\le k\le r$ we have:
\[ 
\si(X_k) = X_k.
\] 
\end{Lem}

\proof
For convenience, we set $\la_{k-n} = \vpi_{i_k}$ if $k-n\le 0$. Using (\ref{defEk}), (\ref{commutD}), and 
the definition of~$\si$, we have
\[
\si(Z_k) = v^{N(\vpi_{i_k}-\la_k) - N(\vpi_{i_k}-\la_{k-n}) - (\vpi_{i_k}-\la_{k-n},\ \vpi_{i_k}+\la_{k})}Z_k.
\]
A simple calculation using (\ref{relation-n}) shows that
\[
N(\vpi_{i_k}-\la_k) - N(\vpi_{i_k}-\la_{k-n}) - (\vpi_{i_k}-\la_{k-n},\ \vpi_{i_k}+\la_{k})
=
N(\b_k)+2(\vpi_{i_k}-\la_{k-n},\ \b_k),
\]
and the lemma follows.
\cqfd

Clearly, the rescaled generators $X_k$ satisfy the same commutation relations as the
$Z_k$. Therefore, if we define for $\aa:=(a_1,\ldots,a_r)\in\Z^r$,
\begin{equation}
X^{\aa}:= v^{\frac{1}{2}\sum_{i<j}a_ia_j(\b_i,\,\b_j)}\,X_1^{a_1}\cdots X_r^{a_r}, 
\end{equation}
we have by Proposition~\ref{lemcruc},
\begin{equation}
\si(X^{\aa}) = v^{-\frac{1}{2}\sum_{i<j}a_ia_j(\b_i,\,\b_j)}\,X_r^{a_r}\cdots X_1^{a_1}
= X^{\aa}. 
\end{equation}
Thus, $X^{\aa}$ is $\si$-invariant, and more generally an element of $\T_{v,\,Q}$ is $\si$-invariant 
if and only if all the coefficients of its
expansion with respect to the basis $\{X^{\aa}\mid \aa\in\Z^r\}$ are
invariant under the map $v^{1/2}\mapsto v^{-1/2}$. 
Moreover, one checks easily that
\begin{equation}\label{multXa}
X^{\aa}X^{\bb} =  v^{\frac{1}{2}\sum_{i<j}(a_jb_i-a_ib_j)(\b_i,\,\b_j)}\,X^{\aa+\bb}
= v^{\sum_{i<j}(a_jb_i-a_ib_j)(\b_i,\,\b_j)}\,X^{\bb}X^{\aa}.
\end{equation}

\subsection{The isomorphism $\Phi$}

We can now state the main result of this section, which follows immediately from
Proposition~\ref{prop2.2} and Equations (\ref{mult_mono}), (\ref{multXa}).

\begin{Prop}\label{main_prop}
There is a $\C$-algebra isomorphism $\Phi\colon \Y_{t,\,Q} \to \T_{v,\,Q}$ given
by
\[
\Phi(t^{1/2}) = v^{1/2},\qquad 
\Phi(Y_{i,p}) = X_k  \quad\mbox{for}\quad (i,p)\in \hI_Q\quad \mbox{and}\quad \varphi(i,p) = (\b_k,0).
\]
More generally, let 
\[
m = \prod_{(i,p)\in \hI_Q} Y_{i,p}^{u_{i,p}(m)}
\]
be a \emph{commutative} monomial in $\Y_{t,\,Q}$,
as in \S\ref{commut_prod}, and let
$\aa = (a_1,\ldots,a_r)$ where $a_k = u_{(i,p)}(m)$ if $\varphi(i,p)=(\b_k,0)$.
Then we have
\[
\Phi(m)
= X^{\aa}.
\]
\cqfd
\end{Prop}


\section{Quantum groups}\label{sect4}


\subsection{Background}\label{defAqn}

Let $\n$ denote a maximal nilpotent subalgebra of $\g$.
Let $U_v(\n)$ be the Drinfeld-Jimbo quantum enveloping algebra of $\n$
over $\C(v)$,
with Chevalley generators $e_i\ (i\in I)$ subject to the quantum Serre
relations:
\[
\begin{array}{lll}
&e_i\,e_j - e_j\,e_i = 0 &\mbox{if $C_{ij} = 0$,}\\[2mm]
&e_i^2\,e_j - (v+v^{-1})e_i\,e_j\,e_i + e_j\,e_i^2= 0\qquad &\mbox{if $C_{ij} = -1$.} 
\end{array}
\] 
It is endowed with a natural scalar product $(\cdot,\cdot)$ 
which we normalize by $(e_i,e_i) = 1$ (see \eg \cite[\S4.3]{GLS}).
We denote by $A_v(\n)$ the graded dual vector space of $U_v(\n)$.
The map $x \mapsto (x,\cdot)$ is a vector space isomorphism from
$U_v(\n)$ to $A_v(\n)$, which allows to define a multiplication
on $A_v(\n)$ by transporting the multiplication of $U_v(\n)$.

Thus, $U_v(\n)$ and $A_v(\n)$ are isomorphic algebras, but they have
dual integral forms and therefore they specialize differently
at $v=1$.
One should regard $A_v(\n)$ as a quantum coordinate ring of the
unipotent group $N$ with Lie algebra $\n$.
For example, the elements $D_{\vpi_{i_k},\,\la_k}$ of \S\ref{qtorusA} are 
quantum analogues of certain generalized minors on $N$.
We set
\[
\U_v(\n) := \C(v^{1/2})\otimes_{\C(v)}U_v(\n),
\qquad
\A_v(\n) := \C(v^{1/2})\otimes_{\C(v)}A_v(\n). 
\]
Since the basis involved in Theorem~\ref{mainTh}~(b) is the dual canonical
basis $\bB^*$, it is more natural to think of the quantum algebra of
Theorem~\ref{mainTh}~(a) as being $\A_{v}(\n)$ rather than $\U_{v}(\n)$.

The algebra $U_v(\n)$ has a natural grading by the root lattice of $\g$,
given by $\deg(e_i)=\a_i$. The above isomorphism allows to transfer this
grading to $A_v(\n)$.

\subsection{Determinantal identities}
In \cite{GLS}, it is shown that $A_v(\n)$ has a quantum cluster algebra
structure. In particular, an explicit realization of $A_v(\n)$ as a subalgebra
of the quantum torus $\TT_{v,\,Q}$ is given.
This goes as follows.

For $u, w \in W$ and $\la\in P_+$, one has unipotent quantum minors
$D_{u(\la),w(\la)}\in A_v(\n)$ (see \cite[\S5.2]{GLS}).
They satisfy
\[
D_{u(\la),w(\la)} = 
\left\{
\begin{array}{cl}
1 & \mbox{if $u(\la)=w(\la)$,}\\[2mm]
0 & \mbox{if $u(\la)\not \le w(\la)$}.
\end{array}
\right.  
\]
Let $\ii = (i_1,\ldots,i_r)$ be as in \S\ref{qtorusA}.
In \cite[\S5.4]{GLS}, a system of identities
relating the quantum minors
\begin{equation}
D(b,d) := D_{s_{i_1}\cdots s_{i_b}(\vpi_{i_b}),\ s_{i_1}\cdots s_{i_d}(\vpi_{i_b})},
\qquad
(0\le b \le d \le r,\ \ i_b=i_d\in I), 
\end{equation}
is derived, which we now recall.
By convention, we write $D(0,b) = D_{\vpi_{i_b},\ s_{i_1}\cdots s_{i_b}(\vpi_{i_b})}$. 
Note that the minors $D(0,b)\ (1\le b\le r)$ form by definition
a system of generators of $\TT_{v,\,Q}$.
We will also use the following shorthand notation:
\begin{eqnarray}
b^-(j) &:=& \max\left(\{s < b \mid i_s = j\}\cup\{0\}\right),\\
b^- &:=& \max\left(\{s < b \mid i_s = i_b\}\cup\{0\}\right),\label{defb-}\\
\mu(b,j) &:=&  s_{i_1}\cdots s_{i_b}(\vpi_j).\label{equamu}
\end{eqnarray}
In (\ref{equamu}) we understand that $\mu(0,j) = \vpi_j$. 
\begin{Prop}[\cite{GLS}]\label{Tsystem}
Let $1\le b<d\le r$ be such that $i_b = i_d = i$. 
There holds
\begin{equation}\label{eqTsystem}
v^A\, D(b,d) D(b^-,d^-) = 
v^{-1+B}\, D(b,d^-)D(b^-,d)
\ +\ v^C
\prod_{j\sim i}^{\longrightarrow}D(b^-(j),d^-(j))
\end{equation}
where
\[
A = (\mu(d,i),\,\mu(b^-,i)-\mu(d^-,i)),\qquad
B = (\mu(d^-,i),\,\mu(b^-,i)-\mu(d,i)), 
\]
and 
\[
C 
= 
\sum\limits_{\substack{j<k\\ j\sim i,\ k\sim i}}
\left(\mu(d,j),\,\mu(b,k)-\mu(d,k)\right).
\]
\end{Prop}

This system of identities allows to express inductively every minor $D(b,d)$
as a rational function of the flag minors $D(0,c)$.
Moreover, it follows from \cite[Theorem 12.3]{GLS} that all these rational
functions belong in fact to $\TT_{v,\,Q}$, and that 
$A_v(\n)$ is the subalgebra of $\TT_{v,\,Q}$
generated by the minors $D(b^-,b)\ (1\le b \le r)$.

\subsection{The dual canonical basis $\bB^*$}\label{bases}
Let us write 
\begin{equation}
E^*(\b_k) := D(k^-,k),\qquad (1\le k \le r),
\end{equation}
and for $\aa = (a_1,\ldots,a_r)\in\Z^r$,
\begin{equation}
E^*(\aa) = v^{-\sum_{k=1}^r a_k(a_k-1)/2}E^*(\b_1)^{a_1}\cdots E^*(\b_r)^{a_r}. 
\end{equation}
Then $\bE^*=\{E^*(\aa) \mid \aa \in\Z^r\}$ is a $\C(v)$-basis of $A_v(\n)$, dual to
a basis of $U_v(\n)$ of PBW-type, as defined by Lusztig.
The basis $\bE^*$ is called the \emph{dual PBW-basis} of $A_v(\n)$.

The involution $\si$ of $\T_{v,\,Q}$ (see \S\ref{sigma}) can be 
restricted to $A_v(\n)$.
Lusztig \cite{Lu} has constructed a canonical basis $\bB$ of $U_v(\n)$. The dual basis 
$\bB^* = \{B^*(\aa) \mid \aa\in\N^r\}$ of $A_v(\n)$ can be characterized as follows
(see \eg \cite{GLS}).

\begin{Prop}\label{caracter}
For $\aa=(a_1,\ldots,a_r)\in\N^r$, the vector $B^*(\aa)$ is uniquely determined by the following
conditions:
\begin{itemize}
\item[(a)] $B^*(\aa) \in E^*(\aa) + \sum_{\cc\not = \aa} v^{-1}\Z[v^{-1}]  E^*(\cc)$;
\item[(b)] let $\b(\aa) := \sum_{1\le k\le r}a_k \b_k$. Then $\si(B^*(\aa)) = v^{N(\b(\aa))}B^*(\aa)$.
\end{itemize}
\end{Prop}
The integer $N(\ga)$ of (b) is defined in (\ref{defNgamma}).
Note that $\b(\aa)$ is just the weight of $B^*(\aa)$ or $E^*(\aa)$ in the natural
grading of $A_v(\n)$ by the root lattice of $\g$.
The basis $\bB^*$ is called the \emph{dual canonical basis} of $A_v(\n)$.


\section{Quantum Grothendieck rings}\label{sect5}


\subsection{Background}

For recent surveys on the representation theory of quantum loop algebras,
we invite the reader to consult \cite{CH} or \cite{L}.
 
Let $L\g$ be the loop algebra attached to $\g$, and let $U_q(L\g)$ be the
associated quantum enveloping algebra. 
We assume that the deformation parameter $q\in\C^*$ is not a root of unity.

By \cite{FR}, every finite-dimensional $U_q(L\g)$-module $M$
(of type 1) has a $q$-character $\chi_q(M)$. 
These $q$-characters generate a commutative $\C$-algebra
isomorphic to the complexified Grothen\-dieck ring of the category of  
finite-dimensional irreducible $U_q(L\g)$-modules.
Nakajima \cite{N}, Varagnolo and Vasserot \cite{VV}, and Hernandez \cite{H1}, have
studied $t$-deformations 
of the $q$-charac\-ters of the standard modules
and of the simple modules, as well as
corresponding $t$-deforma\-tions of the Grothendieck ring.
Although slightly different, these $t$-deformed Grothendieck rings are 
essentially equivalent, and in particular they give rise to the same
$(q,t)$-characters for the simple modules. In what follows, we will use
the version of \cite{H1}. Its definition will be recalled in the next
sections.

\subsection{The subcategory $\CC_\Z$}\label{subsectCZ}
The simple finite-dimensional irreducible $U_q(L\g)$-modules (of type 1)
are usually labelled by Drinfeld polynomials. Here we shall use an 
alternative labelling by dominant monomials (see  \cite{FR}).
Moreover, as in \cite{HL}, we shall restrict our attention to a certain tensor subcategory 
$\CC_\Z$ of the category of finite-dimensional $U_q(L\g)$-modules. 
The simple modules in 
$\CC_\Z$ are labelled by the dominant monomials
in $\Y$, or equivalently, by the dominant commutative monomials in
$(\Y_{t^{1/2}},*)$ (see \S\ref{commut_prod}), and their $q$-characters
belong to $\Y$.
We shall denote by $L(m)$ the simple module labelled by the dominant monomial $m$. 
When $m=Y_{i,p}$ is reduced to a single variable, $L(m)$
is called a \emph{fundamental module}.
When $m$ is the only dominant monomial occuring in $\chi_q(L(m))$,
$L(m)$ is said to be \emph{minuscule}.
Fundamental modules are examples of minuscule modules \cite{FM}.

\subsection{Standard modules}

To a dominant commutative monomial
$m$ is also attached a tensor product
of fundamental modules called
a \emph{standard module} $M(m)$ defined by
\begin{equation}\label{standard}
M(m) := \ \stackrel{\longrightarrow}{\bigotimes_{(i,p)\in\hI}} L( Y_{i,p})^{\otimes\, u_{i,p}(m)},
\end{equation}
where the product is ordered according to the following partial order on $\hI$:
\[
(i,p) < (j,s) \quad \Longleftrightarrow \quad p < s. 
\]
Note that for any fixed $p\in \Z$ and any total order on $I$, the tensor product 
\[
\stackrel{\longrightarrow}{\bigotimes_{i\in I}} L( Y_{i,p})^{\otimes\, u_{i,p}(m)} 
\]
is irreducible, and its isomorphism class
$L\left(\prod_{i\in I}Y_{i,p}^k\right)$
does not depend on the order of the factors, hence (\ref{standard}) is well defined
up to isomorphism (see \cite[Proposition 6.15]{FM}).
The classes $[M(m)]$ of the standard modules $M(m)$ form a second basis of the Grothendieck
group of $\CC_\Z$.

\subsection{The ring $\K_t$}
We introduce the commutative monomials \cite{FR}
\begin{equation}\label{defA}
A_{i,p} = Y_{i,p+1}Y_{i,p-1}\prod_{j \sim i} Y_{j,p}^{-1},\qquad ((i,p-1)\in \hI). 
\end{equation}
Recall from \S\ref{commut_prod} that commutative monomials in $\Y$ can be 
regarded as elements of $(\Y_{t},*)$. 
More generally, the commutative polynomials
\[
Y_{i,p}\left(1+A_{i,\,p+1}^{-1}\right) 
=  Y_{i,p}\ + \ Y_{i,p+2}^{-1}\prod_{j \sim i} Y_{j,p+1},\qquad ((i,p)\in \hI)
\]
can be regarded as elements of $(\Y_{t},*)$.
For $i\in I$, let $\mathcal{K}_{i,t}$ be the $\C(t^{1/2})$-subalgebra of $\Y_{t}$ (for the noncommutative
product $*$) generated by 
\[
Y_{i,p}\left(1+A_{i,\,p+1}^{-1}\right),\quad
Y_{j,s}^{\pm 1},
\qquad \left((i,p), (j,s) \in \hI,\ j\not = i\right).
\]
(In \cite{H0}, $\mathcal{K}_{i,t}$ is identified with the kernel of a $t$-deformed screening
operator.)
Define 
\[
\mathcal{K}_t := \bigcap_{i\in I}\mathcal{K}_{i,t}.
\] 
It is shown in \cite{H1} that
an element of $\mathcal{K}_t$ is uniquely determined by the coefficients
of its dominant monomials. Moreover, for any dominant monomial $m$, there is a unique
$F(m)\in\mathcal{K}_t$ such that $m$ occurs in $F(m)$ with multiplicity $1$ and
no other dominant monomial occurs in $F(m)$. 
These $F(m)$ form a $\C(t^{1/2})$-basis of $\mathcal{K}_t$.

\subsection{Comparison with other $t$-deformations}\label{differ}
The product $*$ used in this paper is the same as that of \cite{H1}, except that we have replaced
$t$ by~$t^{-1}$. The product of \cite{H1} is slightly different from the products of \cite{N} and \cite{VV}
(see Remark~\ref{Remark_differ}).
However, as shown in \cite[Proposition 3.16]{H1}, for every $(i,p), (j,s) \in \hI$ the pairs $(Y_{i,p}, A_{j,s})$ and $(A_{i,p}, A_{j,s})$ 
are $t$-commutative with the \emph{same exponents of $t$} for the three products of \cite{N,VV,H1}. 
This implies that the $t$-deformations of the Grothendieck ring $\R$ of $\CC_\Z$
associated with the three products are 
essentially equivalent, as will be explained below.

\subsection{$(q,t)$-characters of standard modules} \label{sectqtstandard}
For a dominant commutative monomial $m\in \Y_{t^{1/2}}$, define
\begin{equation}\label{defMt}
[M(m)]_t := t^{\alpha(m)}\, 
\underset{p\in\Z}{\overset{\longleftarrow}{\mbox{\Large *}}} 
\,F\left(\prod_{i\in I} Y_{i,p}^{u_{i,p}(m)}\right)
\in \mathcal{K}_t. 
\end{equation}
Here, $\alpha(m)\in\frac{1}{2}\Z$ is chosen so that $m$ occurs with coefficient $1$
in the expansion of $[M(m)]_t$ on the basis of commutative monomials of $\Y_{t^{1/2}}$.
The coefficients of $[M(m)]_t$ on this basis belong to $\Z[t^{\pm1}]$ and may
therefore be specialized at $t=1$.  
The obtained specialization of $[M(m)]_t$ at $t=1$ is equal to $\chi_q(M(m))$, the $q$-character
of the standard module $M(m)$.
Therefore we may use the alternative notation
\[
\chi_{q,t}(M(m)) := [M(m)]_t, 
\]
and call this element of $\K_t$ the \emph{$(q,t)$-character of $M(m)$}.

\subsection{The bar involution}\label{susbsectinvol}
One shows that there is a unique $\C$-algebra anti-automorphism of $(\Y_{t},*)$ such that
\[
\overline{t^{1/2}} = t^{-1/2}, \quad
\overline{Y_{i,p}}=Y_{i,p},\qquad ((i,p)\in \hI).
\]
Clearly, the \emph{commutative} monomials are bar-invariant, as in \cite{N,VV}.
The subring $\mathcal{K}_t$ is stable under the bar involution,
since each $\K_{i,t}$ is obviously stable. 
It follows that the elements $F(m)$ are bar-invariant (since $m$ is the unique 
dominant monomial of $\overline{F(m)}$). 
Hence the coefficients of the expansion of $F(m)$ on the basis of commutative monomials
are unchanged under the replacement of $t$ by $t^{-1}$.
Therefore, $F(m)$ is the same as in \cite{H1}.
Since we have used in (\ref{defMt}) the reverse product~$\overset{\leftarrow}{*}$,
the elements $\chi_{q,t}(M(m))$ also coincide with the corresponding elements
of \cite{N,VV,H1}, \ie the coefficients of their expansion on the basis of 
commutative monomials are the same.

\subsection{$(q,t)$-characters of simple modules} \label{sectqtsimple}

\begin{Prop}[\cite{N}]\label{defLt}
For every dominant monomial $m$, there is a unique element $[L(m)]_t$ of $\K_t$ 
satisfying 
\begin{itemize}
\item[(a)] $\overline{[L(m)]_t} = [L(m)]_t$, 
\item[(b)] $\ds [L(m)]_t \in [M(m)]_t + \sum_{m' < m} t^{-1}\Z[t^{-1}]\, [M(m')]_t.$
\end{itemize}
Here $m'\leq m$ means that $m(m')^{-1}$ is a product of elements $A_{i,p}$ in $\Y$. 
\end{Prop}
By \S\ref{susbsectinvol}, the elements $[L(m)]_t$ coincide with the corresponding elements
of \cite{N,VV,H1}.
Using the geometry of quiver varieties, Nakajima has shown:
\begin{Thm}[\cite{N}]\label{NakaTh}
The specialization 
of $[L(m)]_t$ at $t=1$ is equal to $\chi_q(L(m))$, and 
the coefficients of the expansion of $[L(m)]_t$ as a linear combination of 
monomials in the $Y_{i,p}$'s belong to $\N[t^{\pm 1}]$. 
\end{Thm}
Therefore we may use the alternative notation
\[
\chi_{q,t}(L(m)) := [L(m)]_t,
\]
and call this element of $\K_t$ the \emph{$(q,t)$-character of $L(m)$}.

\begin{Cor}
\begin{itemize}
\item[(a)] If $L(m)$ is minuscule, $\chi_{q,t}(L(m)) = F(m)$.
\item[(b)] If $\chi_q(L(m))$ is multiplicity-free, then 
$\chi_{q,t}(L(m)) = \chi_q(L(m))$ does not depend on $t$ when expressed on the basis of
commutative monomials.
\end{itemize}
\end{Cor}
\proof By the positivity statement of Theorem~\ref{NakaTh}, every monomial occuring in $\chi_{q,t}(L(m))$ 
already occurs in $\chi_{q}(L(m))$. Thus, if $L(m)$ is minuscule then $\chi_{q,t}(L(m))$ is an element of $\K_t$
containing the unique dominant monomial $m$, which proves (a).
If $\chi_q(L(m))$ is multiplicity-free, then the coefficient of every commutative monomial in 
$\chi_{q,t}(L(m))$ is of the form $t^k$ for some $k\in\Z$. But since $\chi_{q,t}(L(m))$ is
bar-invariant, we must have $k=0$, which proves (b). \cqfd

\subsection{Multiplicative structure}
We shall regard the noncommutative ring $(\K_t,*)$ as a 
\emph{$t$-deformed version of the Grothen\-dieck ring $\R$}. 
But one should be aware that only the simple modules $L(m)$ and the standard modules $M(m)$
have well-defined ``classes'' $\chi_{q,t}(L(m))$ and $\chi_{q,t}(M(m))$ in $\K_t$.

For any dominant monomials $m_1$ and $m_2$, write
\[
\chi_{q,t}(L(m_1))*\chi_{q,t}(L(m_2)) = \sum_m c_{m_1,\,m_2}^{m}(t^{1/2})\, \chi_{q,t}(L(m)). 
\]
Note that every irreducible $(q,t)$-character is of the form
$\chi_{q,t}(L(m)) = m(1 + \sum_{k} M_k)$,
where the $M_k$ are monomials in the $A_{i,p}^{-1}$ with coefficients in $\N[t,t^{-1}]$
(see \cite{H1}).
So, by \S\ref{differ}, the above coefficients $c_{m_1,\,m_2}^{m}(t^{1/2})$ are obtained 
from the corresponding ones in \cite{N,VV} by multiplying by some $t^k$ with $k\in\Z/2$.
Varagnolo and Vasserot have shown the following
positivity result:
\begin{Thm}[\cite{VV}]\label{posmult}
The structure constants $c_{m_1,\,m_2}^{m}(t^{1/2})$ 
belong to $\N[t^{1/2},t^{-1/2}]$. 
\end{Thm}

\begin{Cor}\label{tcomm}
$L(m_1)\otimes L(m_2) \simeq L(m)$ is a simple module if and only if 
\[
\chi_{q,t}(L(m_1))*\chi_{q,t}(L(m_2))=t^{2k}\chi_{q,t}(L(m_2))*\chi_{q,t}(L(m_1)) 
= t^k\chi_{q,t}(L(m))
\]
for some $k\in\Z/2$.
\end{Cor}

\proof If $L(m_1)\otimes L(m_2) \simeq L(m)$
then $\chi_{q}(L(m_1))*\chi_{q}(L(m_2)) = \chi_{q}(L(m))$. 
Hence $c_{m_1,\,m_2}^{m}(1) = 1$, and
it follows from Theorem~\ref{posmult} that $\chi_{q,t}(L(m_1))*\chi_{q,t}(L(m_2)) = t^k \chi_{q,t}(L(m))$
for some $k\in\Z/2$.
Applying the bar involution, we get 
$\chi_{q,t}(L(m_2))*\chi_{q,t}(L(m_1)) = t^{-k} \chi_{q,t}(L(m))$.
If conversely $\chi_{q,t}(L(m_1))*\chi_{q,t}(L(m_2))= t^k\chi_{q,t}(L(m))$, then
specializing $t$ to $1$ we get 
\[
\chi_{q}(L(m_1))\chi_{q}(L(m_2))= \chi_{q}(L(m_1)\otimes L(m_2)) =\chi_{q}(L(m)),
\]
hence $L(m_1)\otimes L(m_2) \simeq L(m)$.
\cqfd

\subsection{Quantum $T$-system}
For $(i,p)\in \hI$ and $k\in\N$, let $m^{(i)}_{k,p}:= Y_{i,p}Y_{i,p+2}\cdots Y_{i,p+2k-2}$.
The simple $U_q(L\g)$-module 
\[
W^{(i)}_{k,p} := L\left(m^{(i)}_{k,p}\right)
\]
is called a \emph{Kirillov-Reshetikhin module}. (By convention, if $k=0$ then 
$W^{(i)}_{k,p}$ is the trivial one-dimensional module.)
The $q$-characters of the Kirillov-Reshetikhin
modules satisfy the following system of algebraic identities called \emph{$T$-system} \cite{KNS,N-KR,H2}.
For every $(i,p)\in\hI$ and $k>0$, there holds
\[
\chi_q\left(W^{(i)}_{k,p}\right)\chi_q\left(W^{(i)}_{k,p+2}\right) = 
\chi_q\left(W^{(i)}_{k-1,p+2}\right)\chi_q\left(W^{(i)}_{k+1,p}\right)
+
\prod_{j\not= i} \chi_q\left(W^{(j)}_{k,p+1}\right)^{-C_{ij}}. 
\]
This can be lifted to a $t$-deformed $T$-system in $\K_t$, as shown by the next proposition
(see also \cite[\S4]{N-KR}, where a different $t$-deformed product is used, as explained in 
Remark~\ref{Remark_differ} and \S\ref{differ}).
Before stating it, we note that 
$\bigotimes_{j\sim i} W^{(j)}_{k,p+1}$ is a simple module, hence by Corollary~\ref{tcomm}
the $(q,t)$-characters $\chi_{q,t}(W^{(j)}_{k,p+1})$ pairwise $t$-commute in $\K_t$.
Moreover, it is easy to check that, since $\widetilde{C}(z)$ is symmetric,
 for any $j\sim i$ and $j'\sim i$, one has
$m^{(j)}_{k,p+1}*m^{(j')}_{k,p+1}=m^{(j')}_{k,p+1}*m^{(j)}_{k,p+1}$,
hence the $(q,t)$-characters $\chi_{q,t}(W^{(j)}_{k,p+1})$ do in fact pairwise commute in $\K_t$.
So we may write 
$\underset{j\sim i}{*} \chi_{q,t}\left(W^{(j)}_{k,p+1}\right)$ without specifying an ordering of the factors.

\begin{Prop}\label{Tsystem2} 
In $\K_t$ there holds:
\[
\chi_{q,t}\left(W^{(i)}_{k,p}\right)*\chi_{q,t}\left(W^{(i)}_{k,p+2}\right) = 
t^{\a(i,k)}\chi_{q,t}\left(W^{(i)}_{k-1,p+2}\right)*\chi_{q,t}\left(W^{(i)}_{k+1,p}\right)
+
t^{\ga(i,k)}
\underset{j\sim i}{*} 
\chi_{q,t}\left(W^{(j)}_{k,p+1}\right), 
\]
where
\begin{equation}\label{qTsyst}
\a(i,k) = -1 + \frac{1}{2}\left(\widetilde{C}_{ii}(2k-1)+\widetilde{C}_{ii}(2k+1)\right),
\qquad
\ga(i,k) = \a(i,k)+1.
\end{equation}
\end{Prop}
\proof Using Theorem~\ref{posmult}, we see that the claimed identity holds for some integers
$\a(i,k)$ and $\ga(i,k)$, and we only have to check (\ref{qTsyst}).
To do so it is enough to compare the coefficients of some particular monomials on both sides.
We have
$m_{k,p}^{(i)}*m_{k,p+2}^{(i)}= t^{\a}m_{k-1,p+2}^{(i)}*m_{k+1,p}^{(i)}$, where
\begin{eqnarray*}
\a &=& \sum_{a=1}^{k-1}\cN(i,p;\,i,p+2a) + \frac{1}{2}\cN(i,p;\,i,p+2k) \\
&=& \sum_{a=1}^{k-1} \left(\widetilde{C}_{ii}(2a+1) - \widetilde{C}_{ii}(2a-1)\right)
+ \frac{1}{2} \left(\widetilde{C}_{ii}(2k+1) - \widetilde{C}_{ii}(2k-1)\right)\\
&=& -\widetilde{C}_{ii}(1)+\frac{1}{2}\left(\widetilde{C}_{ii}(2k-1) + \widetilde{C}_{ii}(2k+1)\right).
\end{eqnarray*}
Thus $\a(i,k) = \a = -1 + (\widetilde{C}_{ii}(2k-1)+\widetilde{C}_{ii}(2k+1))/2$,
as claimed. 

Similarly, $\chi_q(W^{(i)}_{k,p})$ contains the monomial 
$m:=m_{k,p}^{(i)}A_{i,p+2k-1}^{-1}\cdots A_{i,p+3}^{-1} A_{i,p+1}^{-1}$
with coefficient~1,
and we have $m\,m_{k,p+2}^{(i)} = \prod_{j\sim i} m^{(j)}_{k,p+1}$.
Now 
\[
m * m_{k,p+2}^{(i)} = 
\left(\left(m_{k,p+2}^{(i)}\right)^{-1}\prod_{j\sim i} m^{(j)}_{k,p+1}\right) * m_{k,p+2}^{(i)}=t^{\ga} 
\prod_{j\sim i}m^{(j)}_{k,p+1},
\] 
where
\begin{eqnarray*}
\ga &=& 
\frac{1}{2}\sum_{j\sim i}\sum_{a=1}^k\sum_{b=1}^{k}\cN(j,p+2a-1;\, i,p+2b)\\
&=&
\frac{1}{2}\sum_{j\sim i}\sum_{a=1}^k\sum_{b=1}^{k}
\left( 
\widetilde{C}_{ji}(2(a-b)-2) - \widetilde{C}_{ji}(2(a-b)) - \widetilde{C}_{ji}(2(b-a)) + \widetilde{C}_{ji}(2(b-a)+2)
\right)\\
&=&
\frac{1}{2}\sum_{j\sim i}\sum_{a=1}^{k}
\left(
\widetilde{C}_{ji}(2(a-k)-2) - \widetilde{C}_{ji}(2(a-1)) - \widetilde{C}_{ji}(2(1-a)) + \widetilde{C}_{ji}(2(k-a)+2)
\right)\\
&=& 
\frac{1}{2}\sum_{j\sim i}\sum_{a=1}^{k}
\left(
- \widetilde{C}_{ji}(2(a-1)) + \widetilde{C}_{ji}(2(k-a)+2)
\right)\\
&=& 
\frac{1}{2}\sum_{j\sim i}
\widetilde{C}_{ji}(2k).
\end{eqnarray*}
Thus $\ga(i,k) = \ga = 
\left(\sum_{j\sim i}\widetilde{C}_{ji}(2k)\right)/2 = 
\left(\sum_{j\sim i}\widetilde{C}_{ij}(2k)\right)/2 = 
(\widetilde{C}_{ii}(2k-1)+\widetilde{C}_{ii}(2k+1))/2$,
as clai\-med. 
Here, the last equality comes from the definition of $\widetilde{C}(z)$
(see the proof of Proposition~\ref{formula_inverse}).
\cqfd

\begin{example}
{\rm 
(a)\ \ Take $\g$ of type $A_1$.
We have 
\[
\widetilde{C}(z) = z-z^3+z^5-z^{7}+z^9 -\cdots,\qquad
\] 
hence $\a(k) = -1$ for all $k>0$. 
Thus we get
\[
\chi_{q,t}\left(W_{k,p}\right)*\chi_{q,t}\left(W_{k,p+2}\right) = 
t^{-1}\chi_{q,t}\left(W_{k-1,p+2}\right)*\chi_{q,t}\left(W_{k+1,p}\right)
+ 1.
\]
\bigskip\noindent
(b)\ \ Take $\g$ of type $A_3$. Choose $i=1$, $k=1$, and $p=0$.
Using for example Proposition~\ref{formula_inverse}, we can
calculate 
\[
\widetilde{C}_{11}(z) = z-z^7+z^9-z^{15}+\cdots,\qquad
\widetilde{C}_{12}(z) = z^2-z^6+z^{10}-z^{14}+\cdots, 
\]
hence
\[
\a(1,1) = -1 + \frac{1}{2}\left(\widetilde{C}_{11}(1)+\widetilde{C}_{11}(3)\right) = -\frac{1}{2},
\qquad
\ga(1,1) = \frac{1}{2}\widetilde{C}_{12}(2) = \frac{1}{2}. 
\]
Thus Proposition~\ref{Tsystem2} gives
\[
\chi_{q,t}\left(W^{(1)}_{1,0}\right)*\chi_{q,t}\left(W^{(1)}_{1,2}\right) = 
t^{-1/2} \chi_{q,t}\left(W^{(1)}_{2,0}\right) +
t^{1/2} \chi_{q,t}\left(W^{(2)}_{1,1}\right).
\]
} 
\end{example}

\subsection{The subcategory $\CC_Q$}\label{sectCh'}
Recall the quantum torus $\Y_{t,\,Q}$ of \S\ref{torusT}.
The dominant commutative monomials in $\Y_{t,\,Q}$ para\-me\-tri\-ze the 
simple objects
of an abelian subcategory $\CC_Q$ of $\CC_\Z$.
More precisely, we define $\CC_Q$ as the full subcategory of $\CC_\Z$ whose
objects have all their composition factors of the form $L(m)$
where $m$ is a dominant commutative monomial in $\Y_{t,\,Q}$.
When $Q$ is a sink-source orientation of the Dynkin diagram and
the Coxeter number
$h$ is even, $\CC_Q$ is one of the subcategories $\CC_\ell$ introduced in \cite{HL};
namely, $\CC_Q = \CC_{h'}$ where $h' = h/2 - 1$.
\begin{Lem}\label{tensorCQ}
$\CC_Q$ is closed under tensor products, hence is a tensor subcategory
of $\CC_\Z$. 
\end{Lem}

\proof
This is a slight
modification of the proof of \cite[Proposition 3.2]{HL}.
Let $L(m)$ and $L(m')$ be in $\mathcal{C}_Q$.
This means that $m$ and $m'$ are monomials in 
the variables $Y_{i,p}$, $(i,p)\in\hat{I}_Q$.
If $L(m'')$ is a composition factor of $L(m)\otimes L(m')$
then $m''$ is a product of monomials of $\chi_q(L(m))$ and
$\chi_q(L(m'))$.
So we have $m''=m m' M$ where $M$ is a monomial in the $A_{j,r}^{-1}$. 
Then it is checked as in \cite[Section 5.2.4]{HL} that, for $m''$ to be dominant,
these $(j,r)$ have to satisfy $(j,r-1)\in\hat{I}_Q$ and $(j,r+1)\in\hat{I}_Q$.
It follows that $m''$ depends only on the variables 
$Y_{i,p}$, $(i,p)\in\hat{I}_Q$, because $\hI_Q$ 
is a ``convex slice'' of~$\hI$, that is, it satisfies:
\begin{itemize}
 \item[(i)] if $(i,p),\, (i,p+2k) \in \hI_Q$ for $i\in I,\, p\in\Z,\, k>0$, 
then $(i,p+2j)\in \hI_Q$ for $1\le j\le k-1$;
 \item[(ii)] if $(i,p),\, (i,p+2) \in \hI_Q$ for $i\in I,\, p\in\Z$,
then for every $j\sim i$ we have $(j,p+1)\in \hI_Q$.
\end{itemize}
Hence the result.
\cqfd

\begin{example}
{\rm We continue Example~\ref{exampleYQ}.
We take $\g$ of type $D_4$ and choose $Q$ as in Figure~\ref{FigD4_1}. 
The simple objects of $\CC_Q$ are of the form $L(m)$,
where  
\[
m =Y_{1,0}^{u_{1,0}} Y_{1,-2}^{u_{1,-2}} Y_{1,-4}^{u_{1,-4}}
Y_{2,0}^{u_{2,0}} Y_{2,-2}^{u_{2,-2}} Y_{2,-4}^{u_{2,-4}}
Y_{3,1}^{u_{3,1}} Y_{3,-1}^{u_{3,-1}} Y_{3,-3}^{u_{3,-3}}
Y_{4,2}^{u_{4,2}} Y_{4,0}^{u_{4,0}} Y_{4,-2}^{u_{4,-2}}. 
\]
and $u_{i,p}\in\N$.
}
\end{example}

\subsection{The ring $\K_{t,\,Q}$ and the truncated $(q,t)$-characters}\label{sectKtCh'}

We denote by $\K_{t,\,Q}$ the $\C(t^{1/2})$-subalgebra of $\K_t$ spanned by 
the $(q,t)$-characters $\chi_{q,t}(L(m))$ of the simple objects $L(m)$ in $\CC_{Q}$.
We call $\K_{t,\,Q}$ the \emph{$t$-deformed Grothendieck ring of $\CC_{Q}$}.

The $(q,t)$-character of a simple object $L(m)$ of $\CC_{Q}$ contains in general many
monomials $m'$ which do not belong to $\Y_{t,\,Q}$. By discarding these monomials we
obtain a \emph{truncated $(q,t)$-charac\-ter}. 
We shall denote by $\bchi_{q,t}(L(m))$ the truncated 
$(q,t)$-character of $L(m)$.
One checks that for a simple object $L(m)$ of $\CC_{Q}$,
all the dominant monomials occuring in $\chi_{q,t}(L(m))$
belong to the truncated $(q,t)$-character $\bchi_{q,t}(L(m))$
(the proof is similar to that of \cite{HL} for the category $\CC_1$, as for the proof 
of Lemma~\ref{tensorCQ} above).
Therefore the truncation map 
\[
\chi_{q,t}(L(m))\mapsto \bchi_{q,t}(L(m))
\] 
extends to an injective algebra homomorphism from $\K_{t,\,Q}$ to $\Y_{t,\,Q}$.
In the sequel we shall identify $\K_{t,\,Q}$ with the subalgebra of $\Y_{t,\,Q}$ 
given by the image of this homomorphism.


\section{An isomorphism between quantum Grothendieck rings and quantum groups}\label{sect6}


\subsection{The isomorphism between $\K_{t,\,Q}$ and $\A_{v}(\n)$}

Recall the isomorphism
$\Phi \colon  \Y_{t,\,Q} \to \T_{v,\,Q}$ of Proposition~\ref{main_prop},
and the notation
$$
\A_{v}(\n):=\C(v^{1/2})\otimes A_v(\n).
$$
Define the \emph{rescaled dual canonical basis} of $\A_{v}(\n)$:
\[
\widetilde{\bB}^*:=\left\{\widetilde{B}^*(\aa):=v^{N(\b(\aa))/2}B^*(\aa)\mid B^*(\aa) \in \bB^*\right\}. 
\]
Clearly, the elements of $\widetilde{\bB}^*$ are invariant under the involution $\si$.
The next theorem is Theorem~\ref{mainTh} in a slightly more precise formulation.

\begin{Thm}\label{mainth}
\begin{itemize}
 \item[(a)] $\Phi$ restricts to an 
isomorphism 
\[
\K_{t,\,Q}\overset{\sim}{\longrightarrow} \A_{v}(\n).
\]
 \item[(b)] 
The basis of $\K_{t,\,Q}$ consisting of
the irreducible truncated $(q,t)$-characters $\bchi_{q,\,t}(L(m))$
is mapped by $\Phi$ onto $\widetilde{\bB}^*$.
\end{itemize}
\end{Thm}

\proof
We introduce some necessary notation. For $1\le k \le r$, let
$k_{\rm min} := \min\{1\le s \le r\mid i_s = i_k\}$.
Set $k^{(0)}:= k$ and, for a negative integer $j$, define $k^{(j)}= (k^{(j+1)})^-$, 
where the notation $b^-$ is as in Equation~(\ref{defb-}).
We also note that, by definition of $\la_k$ and $\b_k$, if $k^-\not = 0$ then 
$\tau^{-1}(\la_k)=\la_{k^-}$ and $\tau^{-1}(\b_k)=\b_{k^-}$.

Let us fix some $(i,p)\in \hI_Q$.
By definition of $\Phi$, we have: 
\[
\Phi(Y_{i,p}) = X_k  \quad\mbox{for}\quad (i,p)\in \hI_Q\quad \mbox{and}\quad \varphi(i,p) = (\b_k,0). 
\]
Note that this relation between $(i,p)$ and $\b_k$ implies in particular that $i_k=i$.
Since if $k^-\not = 0$,
\[
\varphi(i,p+2) = (\tau^{-1}(\b_k),0) = (\b_{k^-},0),
\]
we deduce that
$\Phi\left(Y_{i,p}Y_{i,p+2}Y_{i,p+4}\cdots Y_{i,\xi_i}\right)$ is
equal up to a power of $v$ to $X_kX_{k^-}X_{k^{(-2)}}\cdots X_{k_{\rm min}}$,  
that is, up to a power of $v$, to $D(0,k)$.
Since the commutative monomial $Y_{i,p}Y_{i,p+2}\cdots Y_{i,\xi_i}$ is bar-invariant,
its image is $\si$-invariant, so it has to be equal to $v^{N(\vpi_{i_k}-\la_k)/2}D(0,k)$.
Now $Y_{i,p}Y_{i,p+2}\cdots Y_{i,\xi_i}$ is equal to the truncated $(q,t)$-character
of the Kirillov-Reshetikhin module $W^{(i)}_{1+(\xi_i-p)/2,\,p}$. Hence we have shown that
\[
\Phi\left(\bchi_{(q,t)}\left(W^{(i)}_{1+(\xi_i-p)/2,\,p}\right)\right)
=
v^{N(\vpi_{i_k}-\la_k)/2}D(0,k). 
\]

We now want to show that, more generally, for $1\le s\le (\xi_i-p)/2+1$ we have 
\begin{equation}\label{toshow}
\Phi\left(\bchi_{(q,t)}\left(W^{(i)}_{s,\,p}\right)\right)
=
v^{N(\la_{k^{(-s)}}-\la_k)/2}D(k^{(-s)},k). 
\end{equation}
This will be proved by comparing Proposition~\ref{Tsystem} and Proposition~\ref{Tsystem2}.
Let us denote by 
\[
\bD(b,d) := v^{N(\la_b-\la_d)/2}D(b,d) 
\]
the rescaled quantum minors. Note that
\[
N(\la_b-\la_d) = \frac{1}{2}(\la_b -\la_d,\, \la_b - \la_d) -  \deg(\la_b-\la_d) = (\la_b,\,\la_b-\la_d) - \deg(\la_b-\la_d).
\]
We can rewrite Proposition~\ref{Tsystem} as
\begin{equation}\label{TD}
\bD(b,d) \bD(b^-,d^-) = 
v^{X}\, \bD(b,d^-)\bD(b^-,d)
\ +\ v^Y
\prod_{j\sim i}^{\longrightarrow}\bD(b^-(j),d^-(j))
\end{equation}
where
\[
X:= -1+B-A+\frac{1}{2}
\left((\la_b,\,\la_b-\la_d) + (\la_{b^-},\,\la_{b^-}-\la_{d^-}) 
-(\la_{b^-},\,\la_{b^-}-\la_d) - (\la_{b},\,\la_{b}-\la_{d^-})
\right), 
\]
and
\[
Y:=C-A+\frac{1}{2}\left((\la_b,\,\la_b-\la_d) + (\la_{b^-},\,\la_{b^-}-\la_{d^-})
-\sum_{j\sim i} \left(\la_{b^-(j)},\,\la_{b^-(j)}-\la_{d^-(j)}\right)
\right).
\]
Replacing $A$ and $B$ by their values from Proposition~\ref{Tsystem}, and simplifying
the resulting expression, we easily get
\[
X=-1+\frac{1}{2}(\la_b+\la_{b^-},\,\la_{d^-}-\la_d).
\]
Now, writing $i_b=i_d=i$ and $b = d^{(-s)}$,
\[
(\la_b+\la_{b^-},\,\la_{d^-}-\la_d)
=(\la_b,\la_{d^-})-(\la_{b^-},\la_d) 
=(\vpi_i,\tau^{s-1}(\vpi_i))-(\vpi_i,\tau^{s+1}(\vpi_i)).
\]
Hence, using that $\tau^{s-1}(\vpi_i)-\tau^{s+1}(\vpi_i)= \tau^s(\ga_i)+\tau^{s-1}(\ga_i)$,
by Proposition~\ref{prop1.1} we get
\[
X = -1 + \frac{1}{2}\left(\widetilde{C}_{ii}(2s-1)+\widetilde{C}_{ii}(2s+1) \right). 
\]
Similarly, replacing $A$ and $C$ by their values from Proposition~\ref{Tsystem}, and simplifying
the resulting expression, we get
\[
 Y=(\vpi_i,\,\tau(\vpi_i)-\tau^{s+1}(\vpi_i))
+(\vpi_i,\,\vpi_i-\tau^s(\vpi_i))-\frac{1}{2}\sum_{j\sim i}\sum_{k\sim i} (\tau^{(\xi_j-\xi_k)/2}\vpi_j,\,\vpi_k-\tau^s(\vpi_k)).
\]
Using the identities
\[
\vpi_i - \tau^s(\vpi_i) = \sum_{l=0}^{s-1}\tau^l(\ga_i),
\qquad
\tau^l(\ga_i)+\tau^{l+1}(\ga_i)=\sum_{k\sim i}\tau^{l+(1+\xi_k-\xi_i)/2}(\ga_k), 
\]
we get
\[
(\vpi_i,\,\tau(\vpi_i)-\tau^{s+1}(\vpi_i))
+(\vpi_i,\,\vpi_i-\tau^s(\vpi_i))
=
\sum_{k\sim i} \left(\tau^{(\xi_i-\xi_k-1)/2}(\vpi_i),\, \vpi_k-\tau^s(\vpi_k)\right), 
\]
hence, 
\[
Y=
\frac{1}{2}\,\sum_{k\sim i}\left(
\tau^{(\xi_i-\xi_k-1)/2}\left(2\vpi_i-\sum_{j\sim i}\tau^{(\xi_j-\xi_i+1)/2}(\vpi_j)\right),\,\vpi_k-\tau^s(\vpi_k)
\right).
\]
Now,
\[
2\vpi_i-\sum_{j\sim i}\tau^{(\xi_j-\xi_i+1)/2}(\vpi_j)
=
2\vpi_i-\sum_{j\sim i;\ \xi_j-\xi_i=1} \tau(\varpi_j) -\sum_{j\sim i;\ \xi_j-\xi_i=-1} \varpi_j
=
\a_i+ \sum_{j\sim i;\ \xi_j-\xi_i=1} \ga_j
= \ga_i.
\]
Hence
\begin{eqnarray*}
Y&=&
\frac{1}{2}\,\sum_{k\sim i}\left(
\tau^{(\xi_i-\xi_k-1)/2}(\ga_i),\,\vpi_k-\tau^s(\vpi_k)
\right)\\
&=& - \frac{1}{2}\,\sum_{k\sim i}\left(\tau^{-s+(\xi_i-\xi_k-1)/2}(\ga_i),\,\vpi_k\right)\\
&=& \frac{1}{2}\,\sum_{k\sim i}\left(\tau^{s+(\xi_i-\xi_k-1)/2}(\ga_i),\,\vpi_k\right)\\
&=&
\frac{1}{2}\,\sum_{k\sim i}
\widetilde{C}_{ik}(2s). 
\end{eqnarray*}
Therefore $X=\a(i,s)$, $Y=\ga(i,s)$, and by Proposition~\ref{Tsystem2} we see that, 
for any $(i,p')\in \hI_Q$, there holds
in $\K_t$:
\begin{equation}\label{TW}
\chi_{q,t}\left(W^{(i)}_{s,\,p'}\right)*\chi_{q,t}\left(W^{(i)}_{s,\,p'+2}\right) = 
t^{X}\chi_{q,t}\left(W^{(i)}_{s-1,\,p'+2}\right)*\chi_{q,t}\left(W^{(i)}_{s+1,\,p'}\right)
+
t^{Y}
\underset{j\sim i}{*} 
\chi_{q,t}\left(W^{(j)}_{s,\,p'+1}\right).
\end{equation}
It was shown in \cite{GLS} that one can express every quantum minor $D(b,d)$ as a 
(noncommutative) Laurent polynomial in the quantum flag minors $D(0,k)$, by means
of an explicit sequence of applications of Proposition~\ref{Tsystem}. 
Equivalently, every rescaled quantum minor $\bD(b,d)$ can be expressed as a 
Laurent polynomial in the $\bD(0,k)$'s with coefficients in $\Z[v^{\pm1/2}]$, by means
of an explicit sequence of applications of (\ref{TD}). 
By comparing (\ref{TD}) and (\ref{TW}), we see that the $(q,t)$-character
of $W_{s,\ p}^{(i)}$ (where $\varphi(i,p) = (\b_d,0)$ and $b = d^{(-s)}$)
can be expressed by the \emph{same} Laurent polynomial in the $(q,t)$-characters
\[
\chi_{q,t}\left(
W_{1+(\xi_j-p')/2,\ p'}^{(j)}
\right),
\qquad ((j,p')\in\hI_Q), 
\]
where $v^{1/2}$ is replaced by $t^{1/2}$. This proves (\ref{toshow}).
In particular, we have 
\[
\Phi\left(\bchi_{q,t}\left(
L(Y_{i,p})
\right)\right)
=
\Phi\left(\bchi_{q,t}\left(
W_{1,p}^{(i)}
\right)\right)
=
\bD(d^-,d)
=v^{N(\b_d)/2}E^*(\b_d)
,
\qquad
((i,p)\in\hI_Q).
\]
Thus, $\Phi$ maps the truncated $(q,t)$-characters of the fundamental
modules of $\CC_Q$, that is, a set of algebra generators of $\K_{t,\,Q}$,
to the rescaled dual PBW generators of $\A_{v}(\n)$.
This proves (a).

It follows that $\Phi$ maps the truncated $(q,t)$-characters of the standard
modules of $\CC_Q$ to the elements of the dual PBW-basis of $\A_v(\n)$
up to some power of $v$. Let us calculate this power of~$v$.
By Proposition~\ref{defLt}, we have that $[M(m)]_t - [L(m)]_t$ is a linear
combination of $[L(m')]_t$ with coefficients in $t^{-1}\Z[t^{-1}]$, where
$[L(m)]_t$ and the $[L(m')]_t$ are bar-invariant.
On the other hand, note that the rescaling factor $v^{N(\b(\aa))/2}$ of
the dual canonical basis depends only on the weight of the vector $B^*(\aa)$. 
Hence if we write $\widetilde{E}^*(\aa) = v^{N(\b(\aa))/2} E^*(\aa)$,
the transition matrix between the rescaled dual 
PBW-basis $\{\widetilde{E}^*(\aa)\}$ and the rescaled dual canonical basis 
$\{\widetilde{B}^*(\aa)\}$ is identical to the 
transition matrix between the original bases.
Thus, by Proposition~\ref{caracter}, $\widetilde{E}^*(\aa) - \widetilde{B}^*(\aa)$ 
is a linear
combination of $\widetilde{B}^*(\aa')$ with coefficients in $v^{-1}\Z[v^{-1}]$, where
$\widetilde{B}^*(\aa)$ and the $\widetilde{B}^*(\aa')$ are $\si$-invariant.
By Proposition~\ref{main_prop}, $\Phi$ maps the set of bar-invariant elements of $\K_{t,\,Q}$
to the set of $\si$-invariant elements of $\A_{v}(\n)$. This implies that
$\Phi$ maps the 
basis of $\K_{t,\,Q}$ given by the truncated $(q,t)$-characters of the standard 
modules of $\CC_{Q}$, to the rescaled dual PBW-basis
$\{\widetilde{E}^*(\aa)\}$ of $\A_{v}(\n)$.
Finally, using again Proposition~\ref{caracter} and Proposition~\ref{defLt},
this yields (b). \cqfd 

\begin{example}
{\rm 
Let~$\g$ be of type $A_2$. 
Let $Q$ be the quiver of type $A_2$ with height function
$\xi_1=2$ and $\xi_2=1$. 
We have $\bi = (1,2,1)$, and 
\[
D(0,1) = D_{\vpi_1,\,s_1(\vpi_1)},
\quad
D(0,2) = D_{\vpi_2,\,s_1s_2(\vpi_2)},
\quad
D(1,3) = D_{s_1(\vpi_1),\,s_1s_2s_1(\vpi_1)} = D_{\vpi_2,\,s_2(\vpi_2)}. 
\]
Let $e_1$ and $e_2$ be the Chevalley generators of $U_v(\n)$.
In the identification $A_v(\n) \equiv U_v(\n)$ we have
$D(0,1)\equiv e_1$ and $D(1,3)\equiv e_2$.

In this case the quantum torus
$\Y_{t,\,Q}$ is generated by $Y_{1,0}, Y_{1,2}, Y_{2,1}$, so 
$\K_{t,\,Q}$ is generated by $\bchi_{q,t}(L(Y_{1,0})), \bchi_{q,t}(L(Y_{1,2})), \bchi_{q,t}(L(Y_{2,1}))$. 
The isomorphism $\Phi$ of Theorem~\ref{mainth} satisfies
\[
\Phi\left(\bchi_{q,t}(L(Y_{1,2}))\right) = D(0,1),
\quad
\Phi\left(\bchi_{q,t}(L(Y_{1,0}))\right)= D(1,3). 
\]
Thus Theorem~\ref{mainth} implies that $\bchi_{q,t}(L(Y_{1,2}))$ and
$\bchi_{q,t}(L(Y_{1,0}))$ generate $\K_{t,\,Q}$ and satisfy the 
quantum Serre relations.

This can easily be checked by means of the quantum $T$-system.
Indeed we have by Proposition~\ref{Tsystem2}:
\[
\chi_{q,t}(L(Y_{1,0}))*\chi_{q,t}(L(Y_{1,2})) 
= t^{-1/2} \chi_{q,t}(L(Y_{1,0}Y_{1,2})) + t^{1/2} \chi_{q,t}(L(Y_{2,1})),
\]
and by applying the bar-involution
\[
\chi_{q,t}(L(Y_{1,2}))*\chi_{q,t}(L(Y_{1,0})) 
= t^{1/2} \chi_{q,t}(L(Y_{1,0}Y_{1,2})) + t^{-1/2} \chi_{q,t}(L(Y_{2,1})).
\]
Eliminating $\chi_{q,t}(L(Y_{1,0}Y_{1,2}))$ we get 
\[
(t^{-1/2}-t^{3/2})\chi_{q,t}(L(Y_{2,1})) =  
\chi_{q,t}(L(Y_{1,2}))*\chi_{q,t}(L(Y_{1,0})) 
- t \chi_{q,t}(L(Y_{1,0}))*\chi_{q,t}(L(Y_{1,2})),
\]
which shows that $\K_{t,\,Q}$ is generated by $\bchi_{q,t}(L(Y_{1,2}))$ and
$\bchi_{q,t}(L(Y_{1,0}))$. Finally, using that 
\[
\chi_{q,t}(L(Y_{2,1}))*\chi_{q,t}(L(Y_{1,0})) = t^{-1} \chi_{q,t}(L(Y_{1,0}))*\chi_{q,t}(L(Y_{2,1}))  
\]
we obtain that
\begin{eqnarray*}
&&\chi_{q,t}(L(Y_{1,2}))*\chi_{q,t}(L(Y_{1,0}))^2 
- t \chi_{q,t}(L(Y_{1,0}))*\chi_{q,t}(L(Y_{1,2}))*\chi_{q,t}(L(Y_{1,0}))\\[2mm]
&&\qquad=\ 
t^{-1}\chi_{q,t}(L(Y_{1,0}))*\chi_{q,t}(L(Y_{1,2}))*\chi_{q,t}(L(Y_{1,0})) 
-\chi_{q,t}(L(Y_{1,0}))^2*\chi_{q,t}(L(Y_{1,2})), 
\end{eqnarray*}
which is the first quantum Serre relation. The second one is obtained similarly.
}
\end{example}

\begin{example}
{\rm
In this example, we illustrate the calculations behind the proof of
Theorem~\ref{mainth}.
Let~$\g$ be of type $A_3$. 
Let $Q$ be the quiver of type $A_3$ with height function
$\xi_1=\xi_3=2$ and $\xi_2=3$. Thus $Q$ has source $2$ and sinks $1$, $3$.
We take $\bi = (2,1,3,2,1,3)$, hence
\[
\b_1=\a_2,\quad
\b_2=\a_1+\a_2,\quad
\b_3=\a_2+\a_3,\quad 
\b_4=\a_1+\a_2+\a_3,\quad
\b_5=\a_3,\quad
\b_6=\a_1,
\]
and
\[
\begin{array}{lll}
\la_1=\vpi_2-\a_2,&
\la_2=\vpi_1-\a_1-\a_2,&
\la_3=\vpi_3-\a_2-\a_3,\\[2mm] 
\la_4=\vpi_2-\a_1-2\a_2-\a_3,&
\la_5=\vpi_1-\a_1-\a_2-\a_3,&
\la_6=\vpi_3-\a_1-\a_2-\a_3.
\end{array}
\]
Note that in this case $w_0 = c^2$ where $c=s_2s_1s_3$ is a Coxeter element.
Thus, this example illustrates also Corollary~\ref{Cor_cluster}.
The quantum unipotent minors generating $\T_{v,\,Q}$ are
\[
\begin{array}{lll}
D(0,1) = D_{\vpi_2,\,s_2(\vpi_2)},&
D(0,2) = D_{\vpi_1,\,s_2s_1(\vpi_1)},&
D(0,3) = D_{\vpi_3,\,s_2s_1s_3(\vpi_3)},\\[2mm]
D(0,4) = D_{\vpi_2,\,s_2s_1s_3s_2(\vpi_2)},&
D(0,5) = D_{\vpi_1,\,s_2s_1s_3s_2s_1(\vpi_1)},&
D(0,6) = D_{\vpi_3,\,w_0(\vpi_3)}. 
\end{array}
\]
The generators of the dual PBW-basis are
\[
\begin{array}{lll}
E^*(\b_1) = D(0,1),&
E^*(\b_2) = D(0,2),&
E^*(\b_3) = D(0,3),\\[2mm]
E^*(\b_4) = D_{s_2(\vpi_2),\,s_2s_1s_3s_2(\vpi_2)},&
E^*(\b_5)  = D_{s_2s_1(\vpi_1),\,s_2s_1s_3s_2s_1(\vpi_1)},&
E^*(\b_6)  = D_{s_2s_1s_3(\vpi_3),\,w_0(\vpi_3)}. 
\end{array}
\]
The new generators $X_i$ of $\T_{v,\,Q}$ are
\[
\begin{array}{lll}
X_1 = D(0,1),&
X_2 = v^{-1/2}D(0,2),&
X_3 = v^{-1/2}D(0,3),\\[2mm]
X_4 = v^{-1}D(0,4)D(0,1)^{-1},&
X_5 = v^{-1}D(0,5)D(0,2)^{-1},&
X_6 = v^{-1}D(0,6)D(0,3)^{-1}. 
\end{array}
\]
Let us define integers $\la_{ij}$ and $\mu_{ij}$ by 
\[
D(0,i)D(0,j) = v^{\la_{ij}} D(0,j)D(0,i),
\qquad 
X_i X_j = v^{\mu_{ij}} X_jX_i,
\qquad 
(1\le i,j\le 6). 
\]
The matrices $L=[\la_{ij}]$ and $M=[\mu_{ij}]$ are given by
\[
L = \left[ 
\begin{matrix}
0&-1&-1&0&0&0\cr
1&0&0&0&1&-1\cr
1&0&0&0&-1&1\cr
0&0&0&0&0&0\cr
0&-1&1&0&0&0\cr
0&1&-1&0&0&0 
\end{matrix}
\right],
\qquad
M = \left[ 
\begin{matrix}
0&-1&-1&0&1&1\cr
1&0&0&-1&1&-1\cr
1&0&0&-1&-1&1\cr
0&1&1&0&-1&-1\cr
-1&-1&1&1&0&0\cr
-1&1&-1&1&0&0 
\end{matrix}
\right].
\]
The generators of $\Y_{t,\,Q}$ are
$Y_{1,0}, Y_{3,0}, Y_{2,1}, Y_{1,2}, Y_{3,2}, Y_{2,3}$.
The isomorphism $\Phi$ is defined by
\[
\Phi(Y_{1,0}) = X_5, \quad
\Phi(Y_{3,0}) = X_6, \quad
\Phi(Y_{2,1}) = X_4, \quad
\Phi(Y_{1,2}) = X_2, \quad
\Phi(Y_{3,2}) = X_3, \quad
\Phi(Y_{2,3}) = X_1.
\]
The truncated $(q,t)$-characters of the fundamental modules of $\CC_{Q}$
are expressed in terms of commutative monomials by
\[
\begin{array}{rclrcl}
\bchi_{q,t}(Y_{1,2})&=&Y_{1,2},
&\bchi_{q,t}(Y_{1,0})&=&Y_{1,0}+Y_{1,2}^{-1}Y_{2,1}+Y_{2,3}^{-1}Y_{3,2},
\\[2mm]
\bchi_{q,t}(Y_{2,1})&=&Y_{2,1}+Y_{1,2}Y_{2,3}^{-1}Y_{3,2},\ \
&\bchi_{q,t}(Y_{2,3})&=&Y_{2,3},
\\[2mm]
\bchi_{q,t}(Y_{3,2})&=&Y_{3,2},
&\bchi_{q,t}(Y_{3,0})&=&Y_{3,0}+Y_{3,2}^{-1}Y_{2,1}+Y_{2,3}^{-1}Y_{1,2},
\end{array}
\]
Here, we have used the shorthand notation $\bchi_{q,t}(m)$ instead
of $\bchi_{q,t}(L(m))$.
We also have
\[
\bchi_{q,t}(Y_{1,0}Y_{1,2})=Y_{1,0}Y_{1,2},\qquad
\bchi_{q,t}(Y_{2,1}Y_{2,3})=Y_{2,1}Y_{2,3},\qquad
\bchi_{q,t}(Y_{3,0}Y_{3,2})=Y_{3,0}Y_{3,2}.
\] 
Using the expression of $D(0,k)$ in terms of $X_j$'s and the definition of $\Phi$,
one checks that
\[
\begin{array}{lll}
\Phi(Y_{2,3}) = D(0,1),&
\Phi(Y_{1,2}) = v^{-1/2}D(0,2),&
\Phi(Y_{3,2}) = v^{-1/2}D(0,3),\\[2mm]
\Phi(Y_{2,1}Y_{2,3})=v^{-1}D(0,4),&
\Phi(Y_{1,0}Y_{1,2})=v^{-1}D(0,5),&
\Phi(Y_{3,0}Y_{3,2}) = v^{-1}D(0,6),
\end{array} 
\]
in agreement with Theorem~\ref{mainth}.
By Proposition~\ref{Tsystem}, we have
\[
v^{-1}D(1,4)D(0,1) = v^{-1}D(1,1)D(0,4) + D(0,2)D(0,3), 
\]
hence
\[
v^{-1}D(1,4) = (v^{-1}D(0,4) + D(0,2)D(0,3))D(0,1)^{-1}. 
\]
Therefore
\[
\Phi^{-1}\left(v^{-1}D(1,4)\right) = (Y_{2,1}Y_{2,3} + tY_{1,2}*Y_{3,2})*Y_{2,3}^{-1}
= Y_{2,1} + t Y_{1,2}*Y_{3,2}*Y_{2,3}^{-1} 
= Y_{2,1} + Y_{1,2}Y_{3,2}Y_{2,3}^{-1}, 
\]
where the last equality follows from (\ref{mult_mono}).
Thus we have 
\[
\Phi^{-1}\left(v^{-1}D(1,4)\right) =\Phi^{-1}\left(v^{-1}E^*(\b_4)\right) = \bchi_{q,t}(Y_{2,1}), 
\]
in agreement with Theorem~\ref{mainth}.
Next, we have again by Proposition~\ref{Tsystem},
\[
 D(2,5)D(0,2) = v^{-1}D(0,5) + D(1,4).
\]
Hence 
\[
\Phi^{-1}(D(2,5))=\left(Y_{1,0}Y_{1,2} + t(Y_{2,1}+ Y_{1,2}Y_{2,3}^{-1}Y_{3,2})\right)*(t^{1/2}Y_{1,2})^{-1}. 
\]
Now, 
\[
(Y_{1,0}Y_{1,2})*Y_{1,2}^{-1}=t^{1/2}Y_{1,0}*Y_{1,2}*Y_{1,2}^{-1} = t^{1/2}Y_{1,0},
\]
and similarly
\[
Y_{2,1}*Y_{1,2}^{-1} = t^{-1/2}Y_{2,1}Y_{1,2}^{-1},\quad
(Y_{1,2}Y_{2,3}^{-1}Y_{3,2})*Y_{1,2})^{-1} = t^{-1/2} Y_{2,3}^{-1}Y_{3,2}.
\]
Therefore,
\[
\Phi^{-1}(D(2,5))= \Phi^{-1}\left(E^*(\b_5)\right) =
 Y_{1,0} + Y_{1,2}^{-1}Y_{2,1} + Y_{2,3}^{-1}Y_{3,2} =  \bchi_{q,t}(Y_{1,0}).
\]
Similarly, starting from the minor identity
\[
 D(3,6)D(3,0) = v^{-1}D(0,6) + D(1,4).
\]
we deduce that
\[
\Phi^{-1}(D(3,6))=\Phi^{-1}\left(E^*(\b_6)\right)= Y_{3,0} + Y_{3,2}^{-1}Y_{2,1} + Y_{2,3}^{-1}Y_{1,2} =  \bchi_{q,t}(Y_{3,0}).
\]
Thus we have checked that $\Phi$ maps the fundamental characters
\[
\bchi_{q,t}(Y_{1,0}),\quad
\bchi_{q,t}(Y_{3,0}),\quad
\bchi_{q,t}(Y_{2,1}),\quad
\bchi_{q,t}(Y_{1,2}),\quad
\bchi_{q,t}(Y_{3,2}),\quad
\bchi_{q,t}(Y_{2,3}),
\]
to the rescaled dual PBW generators $v^{N(\b_k)/2}E^*(\b_k)$,
in agreement with Theorem~\ref{mainth}.
} 
\end{example}

\subsection{Proof of Corollary~\ref{corollary_binary}}

Let $b_1,\ldots,b_k\in\bB^*$, and let $L_1,\ldots, L_k$ be
the simple objects of $\CC_Q$ such that 
\[
\Phi(\bchi_{q,t}(L_i)) \in v^{\Z/2} b_i,\qquad (1\le i\le k). 
\]
We have 
$\Phi(\bchi_{q,t}(L_1)*\cdots *\bchi_{q,t}(L_k)) \in v^{\Z/2} b_1\cdots b_k$,
thus, by Theorem~\ref{mainth}, $b_1\cdots b_k \in v^\Z\bB^*$ if and only if
$\bchi_{q,t}(L_1)*\cdots *\bchi_{q,t}(L_k)$ is the $(q,t)$-character of
a simple module up to a power of $v$, that is by Corollary~\ref{tcomm},
if and only if $L_1\otimes\cdots\otimes L_k$ is simple. Hence Corollary~\ref{corollary_binary}
follows from \cite{H3}.


\section{A presentation of quantum Grothendieck rings}\label{sect_pres}


In the remaining sections we drop the symbol $*$ for the $t$-deformed product of $\K_t$,
and simply write $xy$ instead of $x*y$.

\subsection{The generators}\label{generators}

Fix an orientation $Q$ of the Dynkin diagram of $\g$.
Define an involution $\nu$ of $I$ by $w_0(\a_i) = -\a_{\nu(i)}$.
For $i\in I$ write $\varphi^{-1}(\a_i,0) = (k_i,p_i)\in\hI_Q$.
Define the following elements of $\K_t$:
\begin{equation}
x_{i,m}^Q := \chi_{q,t}\left(L\left(Y_{\nu^m(k_i),p_i+mh}\right)\right),
\qquad
(i\in I,\ m\in\Z).
\end{equation}
The elements $x_{i,0}^Q$ belong to $\K_{t,\,Q}$ and map to the Chevalley generators
$D_{\vpi_i,\,s_i(\vpi_i)}\equiv e_i$ of $\A_{v}(\n)\equiv \U_{v}(\n)$ under 
the isomorphism $\Phi$ of Theorem~\ref{mainth}.
Hence $\K_{t,\,Q}$ has a presentation given by the generators $x_{i,0}^Q\ (i\in I)$
subject to the relations (see \S\ref{defAqn})
\[
\begin{array}{lll}
&x_{i,0}^Q\,x_{j,0}^Q - x_{j,0}^Q\,x_{i,0}^Q = 0 &\mbox{if $C_{ij} = 0$,}\\[2mm]
&(x_{i,0}^Q)^2\,x_{j,0}^Q - (t+t^{-1})x_{i,0}^Q\,x_{j,0}^Q\,x_{i,0}^Q + x_{j,0}^Q(\,x_{i,0}^Q)^2= 0\qquad &\mbox{if $C_{ij} = -1$.} 
\end{array}
\]
In particular, every $\chi_{q,t}\left(L(Y_{i,p})\right)$ with $(i,p)\in \hI_Q$ 
can be written as a noncommutative polynomial in the $x_{i,0}^Q$'s.

For $m\in\Z$, let $\K^{(m)}$ be the subalgebra of $\K_t$ generated by the 
$x_{i,m}^Q\ (i\in I)$.
Thus $\K^{(0)} = \K_{t,\,Q}$, 
and $\K^{(m)}$ is isomorphic to $\K^{(0)}$ for every $m\in\Z$.
This comes from the fact that $\K_t$ is generated by the fundamental $(q,t)$-characters
$\chi_{q,t}\left(L(Y_{i,p})\right)\ ((i,p)\in\hI)$, and that the assignment 
\[
\chi_{q,t}\left(L(Y_{i,p})\right) \mapsto \chi_{q,t}\left(L(Y_{\nu(i),\,p+h})\right)
\] 
extends to an algebra automorphism $\Sigma$ of $\K_t$. (In fact, $L(Y_{i,p})$ is the
$U_q(L\g)$-module dual to $L(Y_{\nu(i),\,p+h})$ \cite[\S 5]{CP}, see also \cite[Cor. 6.10]{FM}.) 
Let $\hI_{Q,m} := \varphi^{-1}(\De_+\times \{m\})$.
Thus $\chi_{q,t}\left(L(Y_{i,p})\right)\in \K^{(m)}$ for 
$(i,p)\in \hI_{Q,m}$.
Therefore, we have proved:
\begin{Lem}
The elements $x_{i,m}^Q\ (i\in I, m\in\Z)$ generate $\K_t$.
\cqfd
\end{Lem}

\subsection{The presentation}

We start with the following:
\begin{Lem}\label{lemme_commut}
Let $(i,p)\in\hI$ and $(j,p+h)\in\hI$. Write $V := L(Y_{i,p})\otimes L(Y_{j,\,p+h})$.
\begin{itemize}
 \item [(a)] If $j\not = \nu(i)$ then $V$ is simple.
\item [(b)] If $j=\nu(i)$ then $\chi_q(V) = \chi_q\left(L(Y_{i,p}Y_{\nu(i),\,p+h})\right) + 1$.
\item [(c)] In general we have
\[
\chi_{q,t}\left(L(Y_{i,p})\right) \chi_{q,t}\left(L(Y_{j,p+h})\right) 
= t^{-(\a_i,\a_{\nu(j)})} \chi_{q,t}\left(L(Y_{j,p+h})\right) \chi_{q,t}\left(L(Y_{i,p})\right) + \de_{i\,\nu(j)}(1-t^{-2}),
\]
where $\de_{ik}$ is the Kronecker symbol $\de$.
\end{itemize}
\end{Lem}

\proof
Consider the product  $\pi:=\chi_q\left(L(Y_{i,p})\right)\chi_q\left(L(Y_{\nu(i),\,p+h})\right)$.
By \cite[\S6]{FM}, $\chi_q\left(L(Y_{i,p})\right)$ contains only one dominant monomial, namely $Y_{i,p}$,
one anti-dominant monomial, namely $Y_{\nu(i),p+h}^{-1}$, and
all its other monomials involve only variables of the form $Y_{j,m}^{\pm 1}$ with $p<m< p+h$.
It follows that, if $j\not = \nu(i)$, then $\pi$ contains no other dominant monomial than $Y_{i,p}Y_{j,\,p+h}$,
hence $V$ is irreducible and isomorphic to $L(Y_{i,p}Y_{\nu(i),\,p+h})$. This proves (a).

If $j=\nu(i)$ then $\pi$ contains only two dominant monomials, that is, $Y_{i,p}Y_{\nu(i),\,p+h}$
and $1$. Therefore $V$ has at most two composition factors, $L(Y_{i,p}Y_{\nu(i),\,p+h})$ and
the trivial one-dimensional representation. Since $L(Y_{i,p}) = L(Y_{\nu(i),\,p+h})^*$, the trivial
representation is indeed a composition factor of $V$ because 
$U_q(L\g)$ is a Hopf algebra. This proves (b).

It follows that
\[
\chi_{q}\left(L(Y_{i,p})\right) \chi_{q}\left(L(Y_{j,p+h})\right) 
= \chi_{q}\left(L(Y_{i,p}Y_{j,p+h})\right) + \de_{i\,\nu(j)}. 
\]
In $\K_t$, this identity gets $t$-deformed as 
\[
\chi_{q,t}\left(L(Y_{i,p})\right) \chi_{q,t}\left(L(Y_{j,p+h})\right) 
= t^{\frac{1}{2}\cN(i,p;j,p+h)}\chi_{q,t}\left(L(Y_{i,p}Y_{j,p+h})\right) + \de_{i\,\nu(j)}. 
\]
Now using Proposition~\ref{expression_N} and a sink-source orientation $Q$
where $i$ is a source, we see that $\cN(i,p;j,p+h) = \cN(i,0;j,h) =-(\a_i,\a_{\nu(j)})$.
Using the bar involution, we also have
\[
\chi_{q,t}\left(L(Y_{j,p+h})\right) \chi_{q,t}\left(L(Y_{i,p})\right) 
= t^{-\frac{1}{2}\cN(i,p;j,p+h)}\chi_{q,t}\left(L(Y_{i,p}Y_{j,p+h})\right) + \de_{i\,\nu(j)}. 
\]
Then (c) follows by eliminating $\chi_{q,t}\left(L(Y_{i,p}Y_{j,p+h})\right)$ between these 
two equations.
\cqfd

We can now give a presentation of $\K_t$.
\begin{Thm}\label{presentation}
The algebra $\K_t$ is isomorphic to the $\C(t^{1/2})$-algebra $\AA$ presented by
generators $y_{i,m}\ (i\in I, m\in\Z)$ subject only
to the following relations:
\begin{itemize}
 \item[(R1)] for every $m\in\Z$,
\[
\begin{array}{lll}
&y_{i,m}\,y_{j,m} - y_{j,m}\,y_{i,m} = 0 &\mbox{if $(\a_i,\a_j) = 0$,}\\[2mm]
&y_{i,m}^2\,y_{j,m} - (t+t^{-1})y_{i,m}\,y_{j,m}\,y_{i,m} + y_{j,m}\,y_{i,m}^2= 0\qquad &\mbox{if $(\a_i,\a_j) = -1$;} 
\end{array}
\] 
\item[(R2)]
for every $m\in\Z$ and every $i,j\in I$,
\[
y_{i,m}\, y_{j,\,m+1} = t^{-(\a_i,\a_j)} y_{j,\,m+1}\, y_{i,m} + \de_{ij}(1-t^{-2});
\]
\item[(R3)]
for every $p>m+1$ and every $i,j\in I$,
\[
y_{i,m}\,y_{j,p} = t^{(-1)^{p-m}(\a_i,\a_j)} \,y_{j,p}\,y_{i,m}. 
\]
\end{itemize}
\end{Thm}
\proof 
We fix a sink-source orientation $Q$.
We first check that the $x_{i,m}^Q$ satisfy the above relations.
The relations (R1) are the Drinfeld-Jimbo relations for the subalgebra $\K^{(m)}$,
as explained in \S\ref{generators}.
The relations (R2) follow from 
Lemma~\ref{lemme_commut}~(c) when $\xi_i=\xi_j$.
If $\xi_i\not = \xi_j$, then $x_{i,m}^Q\, x_{j,\,m+1}^Q$
corresponds to a tensor product of the form
$L(Y_{i,p})\otimes L(Y_{j,p+1})$ or $L(Y_{i,p})\otimes L(Y_{j,p+2h-1})$.
These two types of tensor products are always irreducible \cite[Proposition 6.15]{FM}.
Using Corollary~\ref{tcomm}, it follows that $x_{i,m}$ and $x_{j,m+1}$ $t$-commute,
and the exponent of $t$ is calculated by means of Proposition~\ref{expression_N}. 
For the relations (R3) we note that $L(Y_{i,p})\otimes L(Y_{j,s})$ is 
irreducible if $s-p>h$ \cite[Proposition 6.15]{FM}, and we conclude similarly.

It follows that we have a surjective homomorphism $F$ from $\AA$ to $\K_t$ given by $y_{i,m}\mapsto x_{i,m}^Q$,
and we have to show that this is an isomorphism. Define $\AA^{(m)}$ as we have defined $\K^{(m)}$ before.
Then $\AA^{(m)}$ is presented by the relations (1) (with $x_{i,m}^Q$ replaced by $y_{i,m}$),
so $F$ restricts to an isomorphism from $\AA^{(m)}$ to $\K^{(m)}$. 
It follows from the relations (R2) and (R3) that every monomial $M$ in the $y_{i,m}$'s can be
rewritten as a linear combination of monomials of the form
$M_{k_1}M_{k_2}\cdots M_{k_s}$ with $M_{k_j}\in \AA^{(k_j)}$ and $k_1>k_2>\cdots>k_s$.
So we have $\AA = \ds\prod_{m\in\Z}^{\leftarrow} \AA^{(m)}$. 
Now each $\K^{(m)}$ has a basis $\B^{(m)}$ consisting of the
$(q,t)$-characters of standard modules that it contains. Taking
\[
\B':= \{ b_{k_1}b_{k_2}\cdots b_{k_s}\mid b_{k_j}\in\AA^{(k_j)},\ F(b_{k_j})\in \B^{(k_j)},\ 
k_1>\cdots >k_s,\ s\in \N\}, 
\]
we get a spanning set of $\AA$ such that $F(\B')$ is a basis of $\K_t$, consisting of
the $(q,t)$-characters of all the standard modules of $\CC_\Z$.  
Hence $\B'$ is a basis of $\AA$ and $F$ is an isomorphism.
\cqfd

\begin{example}
{\rm
Let $\g = \Sl_2$. By Theorem~\ref{presentation}, $\K_t$
is presented by generators $y_m := \chi_{q,t}(L(Y_{2m}))$ indexed
by $m\in\Z$, subject to
\[
\begin{array}{llll}
y_my_{m+1} &=& t^{-2}y_{m+1}y_m + 1-t^{-2},\\[2mm]
y_my_{p} &=& t^{2(-1)^{p-m}}\,y_{p}y_m, &\mbox{if $p>m+1$.} 
\end{array}
\]
} 
\end{example}

\begin{remark}
{\rm
(a)\ \ It was shown by Frenkel and Reshetikhin \cite[Corollary 2]{FR} that the (classical) Grothen\-dieck ring $\R$ of $\CC_\Z$ 
is the polynomial ring in the classes of all fundamental modules $L(Y_{i,p})\ ((i,p)\in \hI)$.
More recently, a presentation of $\R$ in terms of Kirillov-Reshetikhin modules and $T$-systems
was given in \cite[Corollary 2.9]{IIKNS}. 

Note that our presentation of $\K_t$ does \emph{not} yield a new presentation of $\R$.
Indeed, in order to obtain $\R$ from $\K_t$ by specializing $t$ at $1$, one needs to use the integral form $K_t$ defined
in \S\ref{main_iso} below, and the $x_{i,m}^Q$ are not generators of~$K_t$ if $\g \not = \Sl_2$.

\medskip
(b)\ \ For $m\in\Z$, let $\K^{(m,m+1)}$ denote the subalgebra of $\K_t$ generated by $y_{i,m}, y_{i,m+1}\ (i\in I)$.
It follows from Theorem~\ref{presentation} that $\K^{(m,m+1)}$ is isomorphic to the $t$-deformed
boson algebra $\B_t(\g)$ introduced by Kashiwara \cite[\S3.3]{K}.
} 
\end{remark}


\section{Derived Hall algebras}\label{sect_Hall}


\subsection{The Hall algebra $H(Q)$}\label{Hall}

Let $F$ be a finite field, and let $u := |F|^{1/2} \in \RR_{>0}$. 
Let $\md(FQ)$ be the abelian category of representations of $Q$ over $F$.
The twisted Hall algebra $H(Q)$, introduced by Ringel, 
is the $\C$-algebra with basis $\{z_X\}$
labelled by the isoclasses of objects in $\md(FQ)$, with multiplication
\[
z_X z_Y = u^{\<Y,\, X\>}\sum_W g_{X,Y}^W z_W, 
\]
where $g_{X,Y}^W$ is the number of submodules $T$ of $W$ such that
$T\simeq X$ and $W/T\simeq Y$.
Ringel \cite{Ri,Ri2,Ri4} has shown that $H(Q)$ is isomorphic
to the $\C$-algebra $U_u(\n)$ obtained from $U_v(\n)$ by specializing $v$ at $u$. 
In this isomorphism, the basis $\{z_X\}$ is mapped to a PBW-basis of
$U_u(\n)$. In particular, if $S_i$ denotes the 1-dimensional simple
supported on $i\in I$, $z_{S_i}$ is mapped to the Chevalley generator $e_i$.

\subsection{The derived Hall algebra $DH(Q)$}
Let $D^b(\md(FQ))$ be the bounded derived category of $\md(FQ)$.
To\"en \cite[\S7]{T} has associated with this triangulated category
an associative algebra $DH(Q)$ with the following presentation.
The generators $z_X^{[m]}$ are labelled by all pairs $(X,m)$
where $X$ is an isoclass of $\md(FQ)$ and $m\in\Z$. (The pair 
$(X,m)$ corresponds to the stalk complex with $X$ in degree $m$.)
The relations are:
\begin{itemize}
 \item[(D1)] for every $m\in\Z$,
\[
z_X^{[m]} z_Y^{[m]} = u^{\<Y,\, X\>}\sum_W g_{X,Y}^W z_W^{[m]};  
\]
\item[(D2)] for every $m\in\Z$,
\[
z_X^{[m]} z_Y^{[m+1]} = u^{-\<Y,\, X\>}\sum_{W,T} u^{-\<W,T\>}\ga_{X,Y}^{T,W} z_T^{[m+1]}z_W^{[m]};  
\]
\item[(D3)] for $p>m+1$, 
\[
z_X^{[m]}z_Y^{[p]} = u^{(-1)^{p-m}(X,Y)} z_Y^{[p]} z_X^{[m]}.
\]
\end{itemize}
Here, the Hall number $\ga_{X,Y}^{T,W}$ is defined by To\"en as 
\[
\ga_{X,Y}^{T,W} := \frac{|\Ex(W,Y,X,T)|}{|\Aut(X)||\Aut(Y)|},
\]
where $\Ex(W,Y,X,T)$ is the finite subset of $\Hom(W,Y)\times \Hom(Y,X)\times\Hom(X,T)$
consisting of exact sequences $0\to W\to Y \to X \to T \to 0$.
Note that, as in \S\ref{Hall}, we have twisted the multiplication
by inserting in the original Hall product $z_X^{[m]}z_Y^{[p]}$
of \cite{T} a factor $u^{(-1)^{p-m}\<Y,X\>}$, see \cite{S}.

Consider the elements $z_{i,m} := z_{S_i}^{[m]}$ for $i\in I$ and $m\in\Z$.
As in \S\ref{generators}, we see that the $z_{i,m}$ generate $DH(Q)$.
More precisely, we have:
\begin{Prop}\label{presentationDH}
The algebra $DH(Q)$ is generated by the $z_{i,m}\ (i\in I, m\in\Z)$ subject only
to the following relations:
\begin{itemize}
 \item[(H1)] for every $m\in\Z$,
\[
\begin{array}{lll}
&z_{i,m}\,z_{j,m} - z_{j,m}\,z_{i,m} = 0 &\mbox{if $(\a_i,\a_j) = 0$,}\\[2mm]
&z_{i,m}^2\,z_{j,m} - (u+u^{-1})z_{i,m}\,z_{j,m}\,z_{i,m} + z_{j,m}\,z_{i,m}^2= 0\qquad &\mbox{if $(\a_i,\a_j) = -1$;} 
\end{array}
\] 
\item[(H2)]
for every $m\in\Z$ and every $i,j\in I$,
\[
z_{i,m}\, z_{j,\,m+1} = u^{-(\a_i,\a_j)} z_{j,\,m+1}\, z_{i,m} + \de_{ij}\frac{u^{-1}}{u^2-1};
\]
\item[(H3)]
For every $p>m+1$ and every $i,j\in I$,
\[
z_{i,m}\,z_{j,p} = u^{(-1)^{p-m}(\a_i,\a_j)} \,z_{j,p}\,z_{i,m}. 
\]
\end{itemize}
\end{Prop}
\proof
The relations (H1) follow immediately from (D1) and Ringel's theorem.
The relations (H3) follow immediately from (D3).
Let us deduce the relations (H2) from (D2).

If $i\not = j$, the only exact sequences $0\to W\to S_j \to S_i \to T \to 0$
are of the form 
\[
0\to S_j\overset{\f}{\to} S_j \overset{0}{\to} S_i \overset{g}{\to} S_i \to 0
\]
where $0$ means the zero map, and $f$ and $g$ are isomorphisms.
Clearly there are $(|F|-1)^2$ such sequences, and since $|\Aut(S_i)|=|\Aut(S_j)|=|F|-1$,
we get that $\ga_{S_i,S_j}^{S_j,S_i} = 1$.
Hence 
\[
z_{i,m}\, z_{j,\,m+1} = u^{-\<S_j,S_i\>} u^{-\<S_i,S_j\>} z_{j,\,m+1}\, z_{i,m}
= u^{-(\a_i,\a_j)} z_{j,\,m+1}\, z_{i,m}.
\]

If $i=j$, we have two types of exact sequences $0\to W\to S_i \to S_i \to T \to 0$,
namely
\[
0\to S_i\overset{\f}{\to} S_i \overset{0}{\to} S_i \overset{g}{\to} S_i \to 0,
\quad \mbox{and}\quad
0\to 0\overset{0}{\to} S_i \overset{h}{\to} S_i \overset{0}{\to} 0 \to 0,
\]
where $f,g,h$ are isomorphisms.
It follows that
\[
\ga_{S_i,S_i}^{S_i,S_i} = 1,
\qquad\mbox{and}\qquad
\ga_{S_i,S_i}^{0,0} = \frac{1}{|F|-1} = \frac{1}{u^2-1}, 
\]
hence
\[
z_{i,m}\, z_{i,\,m+1} = u^{-(S_i,S_i)} z_{i,\,m+1}\, z_{i,m} + u^{-\<S_i,S_i\>}\frac{1}{u^2-1}. 
\]
This proves (H2). Finally, the proof that relations (H1), (H2), (H3) give a presentation
of $DH(Q)$ is entirely similar to the proof of the analogous statement in 
Theorem~\ref{presentation} 
(the basis $\B^{(m)}$ is replaced by $\{z_X^{[m]} \mid X \mbox{ isoclass of }\md(FQ) \}$), 
and we omit it. \cqfd

\subsection{The isomorphism between $\K_u$ and $DH(Q)$} 

Define the integral form
\[
K_t:= \bigoplus_{L} \C[t^{1/2},\,t^{-1/2}]\ \chi_{q,t}(L) \subset \K_t, 
\]
where the sum runs over all isoclasses $L$ of simple objects in $\CC_\Z$.
By Theorem~\ref{posmult}, this is a subring of $\K_t$. Set
\[
 \K_u := \C \otimes_{\C[t^{1/2},\,t^{-1/2}]} K_t, 
\]
where $\C$ is regarded as a $\C[t^{1/2},\,t^{-1/2}]$-module via the specialization
$t^{1/2}\mapsto u^{1/2}$.

For $\b\in \hDe = \De_+\times\Z$, we denote by $z_\b^{[m]}$ the basis element $z_X^{[m]}$
of $DH(Q)$ with $X\in\md(FQ)$ indecomposable of dimension vector $\b$. 

The following is a slightly more precise formulation of Theorem~\ref{Th_DHall}.

\begin{Thm}\label{main_iso}

There is a $\C$-algebra isomorphism $\iota \colon \K_u \overset{\sim}{\to} DH(Q)$ such that:
\begin{itemize}
\item[(a)] the class of the fundamental $U_q(L\g)$-module $L(Y_{i,p})$ of $\CC_\Z$ is mapped 
by $\iota$ to a scalar multiple of $z_\b^{[m]}$, where $(\b,m) = \varphi(i,p)$. 
\item[(b)] the basis of classes of standard $U_q(L\g)$-modules of $\CC_\Z$ 
is mapped by $\iota$ to a rescaling of the natural basis  of $DH(Q)$ labelled by 
all isoclasses of objects of $D^b(\md(FQ))$.
\end{itemize}
\end{Thm}
\proof
We first assume, as in the proof of Theorem~\ref{presentation}, 
that $Q$ is a sink-source orientation of the Dynkin diagram.
We can rescale the generators $x_{i,m}^Q$ of $\K_t$ by setting
\[
\tx_{i,m}^Q := \frac{1}{u^{1/2}(u-u^{-1})}\ x_{i,m}^Q,\qquad (i\in I,\ m\in\Z). 
\]
Clearly the new generators $\tx_{i,m}^Q$ still satisfy the homogeneous relations (R1) and (R3)
of Theorem~\ref{presentation}, and the relations (R2) become
\[
\tx_{i,m}^Q\, \tx_{j,\,m+1}^Q = t^{-(\a_i,\a_j)} \tx_{j,\,m+1}^Q\, \tx_{i,m}^Q + \de_{ij}\frac{1-t^{-2}}{u(u-u^{-1})^2}.
\] 
Let $\bx_{i,m}^Q = 1\otimes \tx_{i,m}^Q \in \K_u$. 
By Theorem~\ref{presentation}, Proposition~\ref{presentationDH}, 
the assignment $\tx_{i,m}^Q\mapsto z_{i,m}$ extends to an algebra isomorphism $\iota$.
Indeed, in the relations (R2) we have
\[
\frac{1-u^{-2}}{u(u-u^{-1})^2} = \frac{u^{-1}}{u^2-1}
\]
so the generators $\bx_{i,m}^Q$ of $\K_u$ and $z_{i,m}$ of $DH(Q)$ give rise
to identical presentations.

Since the PBW-basis of $U_v(n)$ is orthogonal with respect to the bilinear form
of \S\ref{defAqn}, it only differs from the dual PBW-basis $\bE^*$ by scalar multiples.
Hence by Ringel's theorem,   
it follows from Theorem~\ref{mainth} 
that the classes of fundamental modules 
in $\CC_Q$, which correspond under $\Phi$ to the elements $E^*(\b)\ (\b\in\De_+)$
of $U_v(\n)$, are mapped by $\iota$ to scalar multiples of the $z_{\b}^{[0]}$.
So, if $\varphi(i,p)=(\b,0)$ we have $\iota([L(Y_{i,p})]_u) = \la_{i,p} z_\b^{[0]}$ for some $\la_{i,p}\in\C$.
Therefore, using on one side the automorphism of $\K_u$ given by $[L(Y_{i,p})]_u \mapsto [L(Y_{i,p-2})]_u$,
and on the other side the corresponding automorphism of $DH(Q)$ induced by the Auslander-Reiten translation 
$\tau$ of $D^b(\mod(FQ))$, we get~(a).

Since the classes of standard modules are the ordered products of the $[L(Y_{i,p})]_u$
(up to powers of $u$), and the basis elements of $DH(Q)$ are the ordered products
of the $z_\b^{[m]}$ (up to powers
of $u$), we get (b).

Therefore we have proved Theorem~\ref{main_iso} in the case of a sink-source
orientation.
But the $\C$-algebras $\K_u$ and $DH(Q)$ are both independent of $Q$.
For $\K_u$ this is clear. On the other hand if $Q'$ is another orientation
of the Dynkin diagram, then $D^b(\md(FQ))$ and $D^b(\md(FQ'))$
are equivalent triangulated categories, so $DH(Q)$ and $DH(Q')$
are isomorphic. Thus $\K_u$ is isomorphic to $DH(Q)$ for an arbitrary orientation.
More precisely, recall that the map $\varphi=\varphi_Q\colon \hI \to \hDe$ depends on 
the choice of $Q$. 
There is a triangle equivalence $F_{QQ'}\colon D^b(\md(FQ)) \to D^b(\md(FQ'))$
such that the induced isomorphism $f_{QQ'}\colon DH(Q) \to DH(Q')$
satisfies 
\[
f_{QQ'}\left(z_{\b}^{[m]}\right) = z_{\b'}^{[m']} \quad \mbox{where} 
\quad (\b',m') = \varphi_{Q'}\varphi_Q^{-1}(\b,m). 
\]
Therefore (a) and (b) hold for an arbitrary orientation.
\cqfd

In the proof of Theorem~\ref{presentation}, we have shown
that if $Q$ is a sink-source orientation, the generators $x_{i,m}^Q$ of $\K_t$ 
satisfy the relations (R1), (R2), (R3). We can see now that 
this holds for any orientation $Q$.

\begin{Cor}
The generators $x_{i,m}^Q$ of $\K_t$ satisfy the same relations for
every orientation $Q$ of the Dynkin diagram, namely the relations
(R1), (R2), (R3) of Theorem~\ref{presentation}.
\end{Cor}

\proof
Let $Q$ be any orientation, by Theorem~\ref{main_iso},
the elements $\bx_{i,m}^{Q}$ of $\K_u$
are mapped by $\iota$ to scalar multiples of the 
generators $z_{i,m}$ of $DH(Q)$. 
Now the relations (H1), (H2), (H3) satisfied by the $z_{i,m}$
are independent of $Q$. Moreover, they are all homogeneous
except for (H2) with $i=j$. 
Since scalar multiplication does not affect homogeneous relations,
the elements $1\otimes x_{i,m}^Q$ of $\K_u$ satisfy the 
relations (R1), (R2) $(i \not = j)$, (R3) with $t$ replaced by $u=|F|^{1/2}$. 
Since this is true for every finite field~$F$, it follows that the elements
$x_{i,m}^Q$ of $\K_t$ satisfy the 
relations (R1), (R2) $(i\not = j$), (R3) where $t$ is an indeterminate.

Finally, the relations (R2) $(i=j)$ follow from Lemma~\ref{lemme_commut}~(c) with
$i=\nu(j)$.
\cqfd

\begin{remark}\label{Remark_untwistedDH}
{\rm
Using Remark~\ref{Remark_differ2}, one can modify the presentation of $\K_t$
to obtain a presentation of the deformed 
Grothendieck ring $\R_t$ of \cite{N,VV}. This presentation shows that the specialization of $\R_{t}$
at $t=u^{-1}$ is isomorphic to the 
\emph{non-twisted} derived Hall algebra of $D^b(\md(FQ))$ with the opposite product.
} 
\end{remark}


\section{Quiver varieties}\label{sect7}


In this section we show that the variety $E_\bd$ of representations
of $Q$ with dimension vector $\bd$ can be regarded as 
a Nakajima graded quiver variety $\M^\bullet_0(W^\bd)$ for an appropriate
$\hI$-graded vector space~$W^\bd$. 
Moreover the stratification of $E_\bd$ by $G_\bd$-orbits
coincides with Nakajima's stratification of $\M^\bullet_0(W^\bd)$.
It follows that the set of perverse sheaves used by Lusztig to
define the (dual) canonical basis of $U_v(\n)$ can be identified
with a subset of the set of perverse sheaves used by Nakajima
for describing the classes $[L]_t$ of simple $U_q(L\g)$-modules.
This gives a geometric way of understanding Theorem~\ref{mainth}~(b).

\subsection{The quiver representation space $E_{\bd}$}\label{quiverLusztig}
Let 
$
\bd = (d_i)_{i\in I} \in \N^I
$ 
denote a dimension vector for $Q$.
We will identify $\bd$ with the element $\sum_{i\in I} d_i\a_i$ of the 
root lattice of $\g$.
The variety $E_\bd$ of representations of $Q$ of dimension $\bd$ is by definition
\[
E_\bd := \bigoplus_{i\to j} \Hom_\C(\C^{d_i},\,\C^{d_j}), 
\]
the sum being over all arrows $i\to j$ of $Q$.
This is just a $\C$-vector space of dimension $\sum_{i\to j} d_id_j$, 
but the interesting geometry comes from the following stratification.
Consider the algebraic group
\[
 G_{\bd} := \prod_{i\in I} GL(d_i,\C).
\]
It acts on $E_\bd$ by base change. There are finitely many orbits in one-to-one
correspondence with the isomorphism classes of representations of $Q$ of 
dimension $\bd$. Thus, using Gabriel's theorem, these orbits have a natural
labelling by the set
\[
I_\bd := \left\{\aa = (a_k) \in \N^r \mid \sum_{k=1}^r a_k \b_k = \bd \right\},
\]
where the positive roots $\b_k$ are enumerated as in (\ref{eq_beta}).
Let $\O_\aa$ denote the orbit labelled by the element $\aa$ of $I_\bd$.
Let $IC(\overline{\O_\aa})$ be the intersection cohomology complex
of $\overline{\O_\aa}$, extended by zero on the complement
of $\overline{\O_\aa}$. Let $\H^i(IC(\overline{\O_\aa}))$ be its
ith cohomology sheaf, and $\H^i(IC(\overline{\O_\aa}))_{\cc}$ the
stalk of this sheaf at a point of $\O_\cc$.

Recall from \S\ref{bases} the dual PBW basis $\bE^*$ and the dual canonical basis
$\bB^*$ of $A_v(\n)$. Write
\[
E^*(\cc) = \sum_{\aa\in I_\bd} \kappa_{\aa,\cc}(v)\, B^*(\aa).  
\]
Lusztig has shown:

\begin{Thm}{\rm\cite[\S9, \S10]{Lu}}\label{LusztigTh}
The coefficients $\kappa_{\aa,\cc}(v)$ are given by
\begin{equation}\label{Lusztig_form}
\kappa_{\aa,\cc}(v) = v^{\dim\O_\cc - \dim\O_\aa} \sum_{i\ge 0} v^{i} \dim\H^i(IC(\overline{\O_\aa}))_{\cc}. 
\end{equation}
\end{Thm}
\subsection{Nakajima's variety $\M^\bullet_0(W)$}\label{quiverNakajima}

Let 
\[
W = \bigoplus_{(i,p)\in \hI} W_i(p)
\] 
be a finite-dimensional $\hI$-graded $\C$-vector space.
In his geometric construction of representations of $U_q(L\g)$, Nakajima \cite{N} has associated
with $W$ an affine variety $\M^\bullet_0(W)$ whose definition we shall now recall.

Let $\hJ:=\{(i,p)\in I\times \Z \mid (i,p-1)\in \hI\}$,
and let 
\[
V = \bigoplus_{(i,s)\in \hJ} V_i(s)
\]
be a finite-dimensional $\hJ$-graded $\C$-vector space.
Define
\begin{align*}
&L^\bullet(V,W)= \bigoplus_{(i,s)\in \hJ} \Hom(V_i(s),W_i(s-1)),
\\
&L^\bullet(W,V)= \bigoplus_{(i,p)\in \hI} \Hom(W_i(p),V_i(p-1)), 
\\
&E^\bullet(V)= \bigoplus_{(i,s)\in \hJ;\ j\sim i} \Hom(V_i(s),V_j(s-1)).  
\end{align*}
Put  
$M^\bullet(V,W) = E^\bullet(V) \oplus L^\bullet(W,V) \oplus L^\bullet(V,W)$.
An element of $M^\bullet(V,W)$ is written 
$(B,\a,\b)$, and its components are denoted by:
\begin{align*}
&B_{ij}(s)\in \Hom(V_i(s),V_j(s-1)),
\\
&\a_i(p)\ \in \Hom(W_i(p),V_i(p-1)),
\\
&\b_i(s)\ \in \Hom(V_i(s),W_i(s-1)).
\end{align*}
We denote by $\L^\bullet(V,W)$ the subvariety of the affine space $M^\bullet(V,W)$
defined by the equations
\begin{equation}\label{ADHM}
\a_i(s-1)\b_i(s) + \sum_{j\sim i}  \epsilon(i,j)\, B_{ji}(s-1)B_{ij}(s) = 0,
\qquad ((i,s)\in\hJ),
\end{equation}
where $\epsilon(i,j)=1$ (\resp $\epsilon(i,j)=-1$) if $i\to j$ is an arrow of $Q$
(\resp $i\to j$ is not an arrow of $Q$).
The algebraic group 
\[
G_V := \prod_{(i,s)\in\hJ}GL(V_i(s))
\]
acts on $M^\bullet(V,W)$ by base change in $V$:
\[
 g\cdot (B,\a,\b) = \left((g_j(s-1)B_{ij}(s) g_i(s)^{-1}),\ (g_i(p-1)\a_i(p)),\ (\b_i(s)g_i(s)^{-1})\right).
\]
Note that there is no action on $W$.
This action of $G_V$ preserves the subvariety $\L^\bullet(V,W)$.
One defines the affine quotient
\[
\M^\bullet_0(V,W) := \L^{\bullet}(V,W) \sslash G_V.
\]
By definition, the coordinate ring of $\M^\bullet_0(V,W)$
is the ring of $G_V$-invariant functions on $\L^{\bullet}(V,W)$,
and $\M^\bullet_0(V,W)$ parametrizes the closed $G_V$-orbits.
If $V_i(s) \subseteq V'_i(s)$ for every 
$(i,s)\in \hJ$, then we have a natural
closed embedding $\M^\bullet_0(V,W) \subset \M^\bullet_0(V',W)$.
Finally, one defines
\[
\M^\bullet_0(W) := \bigcup_V \M^\bullet_0(V,W). 
\]
This is an affine variety, acted upon by the algebraic group
\[
G_W := \prod_{(i,p)\in\hI}GL(W_i(p)).
\]
Let $\M^{\bullet\, {\rm reg}}_0(V,W)$ be the open subset of $\M^\bullet_0(V,W)$
parametrizing the closed {\em free} $G_V$-orbits. 
For a given $W$, we have 
$\M^{\bullet\, {\rm reg}}_0(V,W)\not = \emptyset$ only for a finite
number of $V$'s.
Nakajima has shown that this gives a stratification of $\M^\bullet_0(W)$:
\[
\M^\bullet_0(W) = \bigsqcup_V \M^{\bullet\,{\rm reg}}_0(V,W). 
\]
A necessary condition for $\M^{\bullet\, {\rm reg}}_0(V,W)$ to be nonempty is that
\[
\dim W_i(p) - \dim V_i(p+1) - \dim V_i(p-1) + \sum_{j\sim i} \dim V_j(p) \ge 0 
\]
for every $(i,p)\in\hI$. In this case we say that $(V,W)$ is a \emph{dominant pair}.
Equivalently, by (\ref{defA}), the pair $(V,W)$ is dominant if and only if the monomial $Y^W A^V \in \Y$
is dominant, where we use the shorthand notation
\[
Y^W := \prod_{(i,p)\in\hI} Y_{i,p}^{\dim W_i(p)},
\qquad
A^V := \prod_{(i,s)\in\hJ} A_{i,s}^{-\dim V_i(s)}. 
\]
Note that this stratification of $\M^\bullet_0(W)$ is $G_W$-invariant.
Hence each stratum is a union of $G_W$-orbits.

Given a dominant pair $(V,W)$ such that $\M^{\bullet\, {\rm reg}}_0(V,W) \not = \emptyset$,
we denote by $IC_W(V)$ the intersection cohomology complex of the closure of 
the stratum $\M^{\bullet\, {\rm reg}}_0(V,W)$. Let $\H^i(IC_W(V))$ be its $i$th
cohomology sheaf, and $\H^i(IC_W(V))_{V'}$ be the stalk of this sheaf at a point
of $\M^{\bullet\, {\rm reg}}_0(V',W)$.

For a dominant monomial $m\in\Y$, recall from \S\ref{sectqtstandard} and \S\ref{sectqtsimple} 
the $(q,t)$-characters 
$\chi_{q,t}(M(m))$ and $\chi_{q,t}(L(m))$
of the standard and of the simple $U_q(L\g)$-modules labelled by $m$. 
Write
\[
\chi_{q,t}(M(m')) = \sum_{m'} \zeta_{m,m'}(t)\, \chi_{q,t}(L(m)). 
\]
Nakajima has shown:
\begin{Thm}{\rm\cite[\S8]{N}}\label{NakajimaTh}
The coefficients $\zeta_{m,m'}(t)$ are given by
\begin{equation}\label{Nakajima_form}
\zeta_{m,m'}(t) = t^{\dim\M^{\bullet\, {\rm reg}}_0(V',W) - \dim\M^{\bullet\, {\rm reg}}_0(V,W)} 
\sum_{i\ge 0} t^{i} \dim\H^i(IC_W(V))_{V'}, 
\end{equation}
for any pair of strata $\M^{\bullet\, {\rm reg}}_0(V,W)$ and $\M^{\bullet\, {\rm reg}}_0(V',W)$ 
such that $m = Y^WA^V$ and $m'=Y^WA^{V'}$.
\end{Thm}
\begin{remark}
{\rm
(a)\ \ In order to make the comparison between Theorem~\ref{LusztigTh} and Theorem~\ref{NakajimaTh}
easier, we stated Nakajima's formula (\ref{Nakajima_form}) in a different way from
the original one.
In \cite{N}, Nakajima writes
\[
t^{\dim\M^{\bullet\, {\rm reg}}_0(V',W)} 
\sum_{i\ge 0} t^{-i} \dim\H^i(i_{x_{V'}}^{!}IC_W(V)) 
\]
for the right-hand side of (\ref{Nakajima_form}),
but in his degree convention the trivial local system on the open stratum $S = \M^{\bullet\, {\rm reg}}_0(V,W)$
appears in the intersection cohomology complex of $\overline{S}$ in degree $\dim S$, 
while in Lusztig's convention it appears in degree 0. 
Here we follow Lusztig's convention. Moreover Nakajima uses the costalk 
$i_x^{!}$ at a point $x$ instead of the stalk $i_x^*$, which explains the change of $t^i$ into $t^{-i}$. 

\medskip
(b)\ \ A dominant monomial $m$ can be written in several ways as $m = Y^WA^V$.
The fact that the right-hand side of (\ref{Nakajima_form}) depends only on the monomials
$m$ and $m'$, and not on the particular choices of spaces $W, V, V'$, follows from
a transversal slice argument \cite[\S3]{NJAMS}.
} 
\end{remark}

\subsection{An isomorphism}\label{isom}

Let $\bd=(d_i)$ be a dimension vector, as in \S\ref{quiverLusztig}. 
Recall the bijection $\varphi\colon \hI \to \hDe$ of \S\ref{subsect_quivers}.
We define an $\hI$-graded space $W^{\bd}$ by taking
\[
W^{\bd}_j(p) := \C^{d_i} \quad \mbox{if} \quad\varphi(j,p) = (\a_i,0),
\]
and  $W^{\bd}_j(p):= 0$ for all others $(j,p)\in\hI$.
Clearly, the group $G_{W^{\bd}}$ is isomorphic to $G_\bd$ and we may
identify $G_{W^{\bd}} \equiv G_\bd$.

\begin{Prop}\label{propgeo1}
There is a $G_\bd$-equivariant closed immersion of affine varieties 
\[
\Psi\colon \M^\bullet_0(W^{\bd}){\longrightarrow} E_\bd.
\]
\end{Prop}

The proof of Proposition~\ref{propgeo1} will follow from the next two lemmas.
\begin{Lem}\label{Lemdim1}
Let $i\not = i'\in I$ and set $(j,p)=\varphi^{-1}(\a_i,0)$, $(j',p')=\varphi^{-1}(\a_{i'},0)$.
Assume that $p'\le p$, and write $\varphi(j',\xi_j+p'-p+2) = (\b,0)$.
Then, the coefficient of $\a_j$ in the expansion of the root $\b$ 
on the basis of simple roots is equal to $1$ if there is an arrow $i\to i'$
in $Q$, and to $0$ otherwise. 
\end{Lem}

\proof
By definition of $(j,p)$, $(j',p')$, and of $\varphi$ (see \S\ref{subsect_quivers}), we have
\[
 \a_i = \tau^{(\xi_j-p)/2}(\ga_j),\qquad \a_{i'} = \tau^{(\xi_{j'}-p')/2}(\ga_{j'}). 
\]
It follows that
\[
\b = \tau^{-1+(\xi_{j'}-\xi_j-p'+p)/2}(\ga_{j'}) = \tau^{-1+(p-\xi_{j})/2}(\a_{i'}). 
\]
Recall the Ringel bilinear form $\<\cdot,\cdot\>$. It may be characterized by
\[
 \<\a_i,\ga_j\> = \de_{ij},\qquad (i,j\in I).
\]
Hence, the coefficient of $\a_j$ in $\b$ is equal to:
\[
\<\b,\ga_j\> = \<\tau^{-1+(p-\xi_{j})/2}(\a_{i'}),\ga_j\>
=\<\tau^{-1}(\a_{i'}),\tau^{(\xi_j-p)/2}(\ga_j)\> =  
\<\tau^{-1}(\a_{i'}),\a_i\>=-\<\a_i,\a_{i'}\>.
\]
Now
\[
-\<\a_i,\a_{i'}\> = -\dim(\Hom(S_i,S_{i'})) + \dim(\Ext^1(S_i,S_{i'})) = \dim(\Ext^1(S_i,S_{i'})), 
\]
and this is equal to 1 if there is an arrow from $i$ to $i'$ in $Q$, and to $0$ otherwise.
\cqfd

We now introduce an algebra $\tL_Q$ defined by a quiver $\tG_Q$ with relations.
The vertices of $\tG_Q$ are of two types:
\begin{itemize}
 \item $w_j(p)$ for every $(j,p) = \varphi^{-1}(\a_i,0)\ (i\in I)$;
 \item $v_j(p-1)$ for every pair $(j,p)\in \hI_Q$ such that
$(j,p-2)\in \hI_Q$.
\end{itemize}
The arrows of $\tG_Q$ are of three types:
\begin{itemize}
 \item $a_j(p)\colon w_j(p) \to v_j(p-1)$;
 \item $b_j(p)\colon v_j(p) \to w_j(p-1)$;
 \item $\B_{ij}(p)\colon v_i(p)\to v_j(p-1)$ if $j\sim i$.
\end{itemize}
The relations are:
\[
a_i(p-1)b_i(p) = \sum_{j\sim i} \epsilon(i,j) \B_{ji}(p-1)\B_{ij}(p). 
\]
Obviously, as suggested by the notation,
the definition of $\tL_Q$ is so that the affine variety
$\L^\bullet(V,W^{\bd})$ is nothing but the representation variety
of $\tL_Q$ consisting of representations for which the spaces
$W_j^\bd(p) = \C^{d_i}$ are attached to the vertices $w_j(p)$, 
and the spaces $V_j(p-1)$ are attached to the vertices $v_j(p-1)$
(we assume that $V_j(p-1)=0$ if $v_j(p-1)$ is not a vertex of $\tG_Q$).  

For $i\in I$, we denote by $\veps_i$ the idempotent of $\tL_Q$
associated with the vertex $w_j(p)$ such that $(j,p) = \varphi^{-1}(\a_i,0)$.
We endow $I$ with a total ordering such that $i>i'$ if $p>p'$, where
as above $\varphi^{-1}(\a_{i'},0) = (j',p')$. It is well-known that
if there is an arrow $i\to i'$ in $Q$ then $i>i'$. 
\begin{Lem}\label{Lemdiminf1}
For $i\not =i'\in I$, we have 
\[
 \dim\left(\veps_{i'} \tL_Q \veps_i\right) =
\left\{
\begin{matrix}
1 & \mbox{if there is a path from $i$ to $i'$ in $Q$,}\\[2mm]
0 &\mbox{otherwise.} 
\end{matrix}
\right.
\]
\end{Lem}
\proof
If $i<i'$ then $p\le p'$ and clearly $\veps_{i'} \tL_Q \veps_i = 0$.
On the other hand if $i<i'$ there can be no path from $i$ to $i'$
in $Q$. Thus the lemma is clear in this case, and we may assume 
from now on that $i>i'$.

Let $x\in \veps_{i'} \tL_Q \veps_i$, and let $i''\in I$.
Then $x\in \veps_{i'} \tL_Q \veps_{i''} \tL_Q \veps_i$ if and only if
$x$ belongs to the two-sided ideal of $\tL_Q$ generated by
$a_{j''}(p-1)b_{j''}(p)$ where $\varphi(j'',p-1)=(\a_{i''},0)$.
This is because $b_{j''}(p)$ is the only arrow entering in $w_{j''}(p-1)$,
and $a_{j''}(p-1)$ is the only arrow exiting from $w_{j''}(p-1)$.
Note that $\veps_{i'} \tL_Q \veps_{i''} \tL_Q \veps_i \not = 0$
implies that $i>i''>i'$.

Let $\I$ be the the two-sided ideal of $\tL_Q$ generated by
all the degree two paths:
\[
a_{j}(p-1)b_{j}(p), \qquad ((j,p-1) = \varphi^{-1}(\a_i,0),\ (i\in I)).
\]
Then, the algebra $\tL_Q/\I$ is defined by the same relations as the graded 
preprojective algebra $\hL$ of \cite[\S2.8]{L} (but we have the additional
vertices $w_j(p)$ and only a finite set of vertices $v_j(p)$). 
It follows that 
$\dim\left(\left(\veps_{i'} \tL_Q \veps_i\right) / \left(\veps_{i'} \I \veps_i\right)\right)$
is equal to the $v_{j'}(p'+1)$-component of the dimension vector of the
indecomposable projective $\hL$-module with top $v_j(p-1)$. 
Now it is well-known that this dimension vector can be read off from the 
Auslander-Reiten quiver of $Q$. Namely, using our notation, 
the $v_j(p-1)$-component of the dimension vector of the projective
with top $v_i(\xi_i-1)$ is equal to the coefficient of $\a_i$
in the root $\b$ such that $\varphi(j,p) = (\b,0)$. 
The dimension vectors of the remaining indecomposable projectives 
are obtained from these particular ones by translation.
It follows that we can reformulate Lemma~\ref{Lemdim1} as follows:
\begin{equation}\label{dimproj}
 \dim\left(\left(\veps_{i'} \tL_Q \veps_i\right) / \left(\veps_{i'} \I \veps_i\right)\right) =
\left\{
\begin{matrix}
1 & \mbox{if there is an arrow from $i$ to $i'$ in $Q$,}\\[2mm]
0 &\mbox{otherwise.} 
\end{matrix}
\right.
\end{equation}
In particular if $i'$ is the successor of $i$ in the descending total order defined above,
then $\veps_{i'} \I \veps_i = 0$, and we have
\begin{equation}\label{successor}
 \dim\left(\veps_{i'} \tL_Q \veps_i\right)  =
\left\{
\begin{matrix}
1 & \mbox{if there is an arrow from $i$ to $i'$ in $Q$,}\\[2mm]
0 &\mbox{otherwise.} 
\end{matrix}
\right.
\end{equation} 

Assume now that $i>i'$ are such that there is no path from $i$ to $i'$
in $Q$. Then in particular there is no arrow $i\to i'$, so by
(\ref{dimproj}) we have 
\[
\veps_{i'}\tL_Q\veps_i = \veps_{i'} \I \veps_i 
= \sum_{i'<j<i} \left(\veps_{i'}\tL_Q\veps_j\right) \left(\veps_{j}\tL_Q\veps_i\right).
\]
For each summand, we have either no path from $i$ to $j$ or no path from $j$ to $i'$.
So we can iterate the splitting until we obtain an expression of the form
\[
\veps_{i'}\tL_Q\veps_i = 
\sum_{i'<i_1<\cdots<i_k<i} \left(\veps_{i'}\tL_Q\veps_{i_k}\right)\cdots \left(\veps_{i_1}\tL_Q\veps_i\right),
\] 
where in the right-hand side each factor $\veps_{k}\tL_Q\veps_{j}$ is such that either
we have a path from $j$ to $k$ or $k$ is the successor of $j$ and there is no arrow from
$j$ to $k$. Moreover, since we have no path from $i$ to $i'$ each summand contains
at least one factor of the second type, which is equal to $0$ by (\ref{successor}).
Hence we have shown that $\veps_{i'}\tL_Q\veps_i = 0$.

Assume now that there is an arrow $i\to i'$ in $Q$. Then we have as above
\[
\veps_{i'} \I \veps_i 
= \sum_{i'<j<i} \left(\veps_{i'}\tL_Q\veps_j\right) \left(\veps_{j}\tL_Q\veps_i\right),
\]
where for each $j$ we have either no path from $i$ to $j$ or no path from
$j$ to $i'$ (because the Dynkin diagram is a tree). Thus it follows from
above that $\veps_{i'} \I \veps_i = 0$, so by (\ref{dimproj}) we get
$ \dim\left(\veps_{i'} \tL_Q \veps_i\right)  = 1$.

Finally, if there is a path $i\to i_1\to \cdots \to i_k \to i'$ in $Q$, with $k\ge 1$,
then there is no arrow from $i$ to $i'$, and by (\ref{dimproj}) we have
$\veps_{i'}\tL_Q\veps_i = \veps_{i'} \I \veps_i$.
Moreover, this path is unique, and arguing as above we can write
\[
\veps_{i'}\tL_Q\veps_i = 
\left(\veps_{i'}\tL_Q\veps_{i_k}\right)\cdots \left(\veps_{i_1}\tL_Q\veps_i\right), 
\]
where each factor has dimension $1$, so again $ \dim\left(\veps_{i'} \tL_Q \veps_i\right)  = 1$.
\cqfd

\medskip\noindent
{\it Proof of Proposition~\ref{propgeo1} --- \ }
Let $V$ be a $\hJ$-graded space, and pick $(B,\a,\b)\in \L^\bullet(V,W^\bd)$.
As explained above, we can regard $(B,\a,\b)$ as a representation of $\tL_Q$.
Choose two vertices $i$ and $i'$ of $Q$, and set as before
\[
(j, p) = \varphi^{-1}(\a_i,0),
\qquad
(j', p') = \varphi^{-1}(\a_{i'},0).
\]
By Lemma~\ref{Lemdiminf1} we have $\dim(\veps_{i'} \tL_Q \veps_i) \le 1$.
Let $\theta_{ii'}$ be a generator of $\veps_{i'} \tL_Q \veps_i$.
By the proof of Lemma~\ref{Lemdiminf1}, we can normalize
the $\theta_{ii'}$ so that they verify $\theta_{i'i''}\theta_{ii'} = \theta_{ii''}$
for every $i,i',i''\in I$. 
More precisely, if there is a path $i\to i_1\to \cdots \to i_k\to i'$ in $Q$, then 
$\theta_{i,i'} = \theta_{i_ki'}\cdots \theta_{ii_1}$, and if there
is no path from $i$ to $i'$ then $\theta_{ii'}=0$.
Evaluating $\theta_{ii'}$ in the representation $(B,\a,\b)$ we
obtain a linear map 
$\psi_{ii'}\colon W^{\bd}_{j}(p) \to W^{\bd}_{j'}(p')$.
The collection of maps $(\psi_{ii'})$ for all arrows $i\to i'$
of $Q$ can be regarded as a representation $\psi$ of $Q$ of dimension vector~$\bd$.
It follows easily from the definition of the $G_V$-action that
$\psi$ depends only on the $G_V$-orbit of $(B,\a,\b)$.
Hence the assignment $(B,\a,\b) \mapsto \psi$ induces a morphism of varieties
$\Psi_V\colon \M^{\bullet}_0(V,W^\bd) \to E_\bd$.
Moreover, it follows from the known description of the generators
of the coordinate ring of 
$\M^{\bullet}_0(V,W^\bd)$ (see \cite[\S3.1]{Ncl}) that this coordinate ring  
is generated by the matrix coefficients of
the linear maps $\psi_{ii'}$ for all pairs $(i, i')$.
Hence $\Psi_V$ induces a surjective morphism from 
$\C[E_\bd]$ to $\C[\M^{\bullet}_0(V,W^\bd)]$, thus $\Psi_V$ is a closed immersion.
Since for $V$ large enough we have $\M^{\bullet}_0(V,W^\bd) = \M^{\bullet}_0(W^\bd)$,
we obtain a closed immersion $\Psi\colon\M^{\bullet}_0(W^\bd) \to E_\bd$.
By construction, $\Psi$ commutes with the actions of $G_\bd$ on
both varieties. 
\cqfd

\begin{figure}[t]
\[
\xymatrix@-1.0pc{
&\qquad\qquad\qquad\qquad\qquad&\qquad&&
\\
&&& W_3(5)\ar[d]^{\a_3(5)} &
\\
&&& \ar[lld]^{B_{31}(4)} \ar[ld]^{B_{32}(4)}V_3(4)\ar[rd]^{B_{34}(4)} &
\\
&V_1(3) \ar[rrd]^{B_{13}(3)}& V_2(3)\ar[rd]^{B_{23}(3)}&&\ar[ld]^{B_{43}(3)} V_4(3)
\\
&&&\ar[lld]^{B_{31}(2)}\ar[ld]^{B_{32}(2)}V_3(2)\ar[rd]^{B_{34}(2)}&
\\
&V_1(1)\ar[d]^{\b_1(1)} &V_2(1)\ar[d]^{\b_2(1)}&&\ar[d]^{\b_4(1)} V_4(1)
\\
&W_1(0) &W_2(0)&& W_4(0)
}
\]
\caption{\label{FigD0} {\it $(B,\a,\b)$ in type $D_4$}.}
\end{figure}
\begin{example}\label{ExaPropGeo1}
{\rm
Take $\g$ of type $D_4$.
We label the Dynkin diagram so that the central node is numbered~$3$,
and we choose $\xi_1=\xi_2=\xi_4=4$ and $\xi_3=5$. Thus $Q$ has 
a sink-source orientation with source $3$
and sinks $1,2,4$.
Given a dimension vector $\bd = (d_1,d_2,d_3,d_4)$, the corresponding
graded space $W^\bd$ is given by
\[
W_1^\bd(0)=\C^{d_1},\quad
W_2^\bd(0)=\C^{d_2},\quad
W_3^\bd(5)=\C^{d_3},\quad
W_4^\bd(0)=\C^{d_4}, 
\]
and the other $W_i(p)$'s are zero (see the Auslander-Reiten quiver of $Q$
in Figure~\ref{FigD}).
An element $(B,\a,\b)$ of $\L(V,W^\bd)$ is represented in Figure~\ref{FigD0}. 
The defining equations (\ref{ADHM}) read
\begin{eqnarray*}
B_{13}(3)B_{31}(4) + B_{23}(3)B_{32}(4) + B_{43}(3)B_{34}(4)&=&0, \\[2mm]
B_{31}(2)B_{13}(3)\ =\ B_{32}(2)B_{23}(3)\ =\ B_{34}(2)B_{43}(3)
&=&0.
\end{eqnarray*}
Thus, 
\[
\b_1(1)B_{31}(2)B_{43}(3)B_{34}(4)\a_3(5) = 
- \b_1(1)B_{31}(2)B_{23}(3)B_{32}(4)\a_3(5)
\]
and
\[
\b_1(1)B_{31}(2)B_{13}(3)B_{31}(4)\a_3(5)=0. 
\]
Hence we can take 
\[
\psi_{31}:= \b_1(1)B_{31}(2)B_{43}(3)B_{34}(4)\a_3(5),
\]
and similarly
\[
\psi_{32}:= \b_2(1)B_{32}(2)B_{13}(3)B_{31}(4)\a_3(5),
\qquad
\psi_{34}:= \b_4(1)B_{34}(2)B_{23}(3)B_{32}(4)\a_3(5).
\]
We get a representation 
$\psi := (\psi_{31},\,\psi_{32},\,\psi_{34})$
of $Q$ on the space $W^\bd$.
} 
\end{example}

\begin{Prop}\label{propgeo2}
There is a bijection between the set 
of (nonempty) strata $\M^{\bullet\, {\rm reg}}_0(V,W^\bd)$
and the set
$I_\bd$ of $G_\bd$-orbits.
\end{Prop}

\proof
Let us first consider a stratum $\M^{\bullet\, {\rm reg}}_0(V,W^\bd)$.
By \S\ref{quiverNakajima}, the pair $(V,W^\bd)$ is a dominant pair.
This means that we have nonnegative integers $a_j\ (1\le j\le r)$ such
that 
\begin{equation}\label{eqVa}
Y^{W^\bd}A^V = \prod_{j=1}^r Y_{i_j, p_j}^{a_j}. 
\end{equation}
Here for $1\le j \le r$, we have put $(i_j,p_j) = \varphi^{-1}(\b_j,0)$.
Indeed, by definition,
every dominant commutative monomial of the form $Y^{W^\bd}A^V$ 
belongs to $\Y_{t,\,Q}$. 
Moreover we have a natural grading of $\Y_{t,\,Q}$ by the root lattice
of $\g$ given by
\[
 \deg(Y_{i_j,p_j}) = \b_j,\qquad (1\le j \le r).
\]
It is easy to see that 
for every $A_{i,s}\in \Y_{t,\,Q}$ we have $\deg A_{i,s} = 0$.
Therefore 
\[
\sum_{k=1}^r a_j\b_j = \deg\left(Y^{W^\bd}\right) = \bd. 
\]
Hence, to every stratum $\M^{\bullet\, {\rm reg}}_0(V,W^\bd)$
corresponds an element $\aa$ of $I_\bd$ given by
(\ref{eqVa}).  

Conversely, if $\aa\in I_\bd$, we need to show that
$m_\aa:=\prod_{j=1}^r Y_{i_j, p_j}^{a_j}$ 
can be written in the form (\ref{eqVa}) for some 
nonempty stratum $\M^{\bullet\, {\rm reg}}_0(V,W^\bd)$.
By \cite[Th.~14.3.2]{NJAMS}, this is equivalent to the
fact that $m_\aa$ appears 
in the $q$-character of the standard module $M(Y^{W^\bd})$.
For $i\in I$ write as in \S\ref{generators}
$(k_i,p_i)=\varphi^{-1}(\a_i,0)$.
Then we have by definition of $W^\bd$
\[
Y^{W^\bd} = \prod_{i\in I} Y_{k_i,p_i}^{d_i}. 
\]
By \S\ref{generators} we know that the $(q,t)$-characters
of the fundamental modules $L(Y_{k_i,p_i})\ (i\in I)$ 
generate $\K_{t,\,Q}$. Hence the simple module $L(m_\aa)$,
which is an object of $\CC_Q$,
is a composition factor of a standard module of the form
$M(\prod_{i\in I} Y_{k_i,p_i}^{e_i})$ for some nonnegative
integers $e_i$. But, as before, we must have 
\[
\bd=\deg(m_\aa) = \sum_{i\in I}e_i\deg(Y_{k_i,p_i})=\sum_{i\in I}e_i\a_i,
\]
hence $e_i=d_i$ for every $i$. Therefore $m_\aa$ is indeed 
a weight of $M(Y^{W^\bd})$. 
This proves the claim.
\cqfd

\begin{remark}
{\rm 
The proof of Proposition~\ref{propgeo2} shows that 
$\M^{\bullet\, {\rm reg}}_0(V,W^\bd)$ is a 
nonempty stratum of $\M^\bullet_0(W^{\bd})$ 
if and only if
$(V,W^\bd)$ is a dominant pair, a purely combinatorial
condition. In general Nakajima \cite[Th.~14.3.2]{NJAMS} only 
shows that this is a necessary condition.
In representation-theoretic terms, this means that
every dominant monomial of the form $Y^{W^\bd}A^V$ occurs
in the $q$-character of the standard module $M(Y^{W^\bd})$.
}
\end{remark}

\begin{example}
{\rm
\begin{figure}[t]
\[
\def\objectstyle{\scriptstyle}
\xymatrix@-1.0pc{
&&& Y_{3,5} &
\\
&Y_{1,4} &Y_{2,4}&&Y_{4,4}
\\
&&&Y_{3,3}&
\\
&Y_{1,2} &Y_{2,2}&&Y_{4,2}
\\
&&& Y_{3,1}&
\\
&Y_{1,0} &Y_{2,0}& & Y_{4,0} 
\\
}
\qquad
\def\objectstyle{\scriptstyle}
\xymatrix@-1.0pc{
&&& \a_3 &
\\
&\a_1+\a_3 \ar[rru]&\a_2+\a_3\ar[ru]&&\ar[lu]\a_3+\a_4
\\
&&&\ar[llu]\ar[lu]\a_1+\a_2+2\a_3+\a_4\ar[ru]&
\\
&\a_2+\a_3+\a_4\ar[rru] &\a_1+\a_3+\a_4\ar[ru]&&\ar[lu]\a_1+\a_2+\a_3
\\
&&& \ar[llu]\ar[lu]\a_1+\a_2+\a_3+\a_4\ar[ru]&
\\
&\a_1\ar[rru] &\a_2\ar[ru]& & \ar[lu]\a_4 
\\
}
\]
\caption{\label{FigD} {\it The skeleton of $\CC_Q$ and the Auslander-Reiten quiver in type $D_4$}.}
\end{figure}
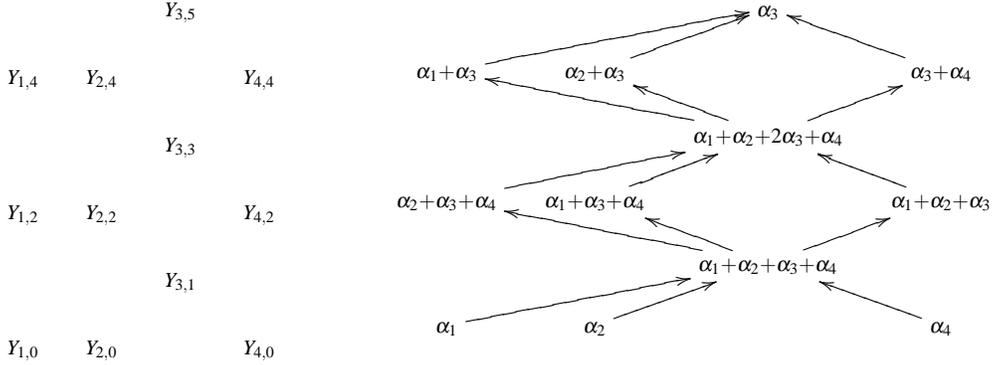
We continue Example~\ref{ExaPropGeo1}.
There are $12$ positive roots $\b_k$, which we 
identify with the vertices of the Auslander-Reiten quiver of $Q$ represented
in Figure~\ref{FigD}. The numbering is obtained by reading this graph 
from top to bottom and left to right:
\[
 \b_1=\a_3,\quad \b_2=\a_1+\a_3,\quad \b_3 = \a_2+\a_3,\quad \b_4 = \a_3+\a_4,
\quad \cdots ,\quad \b_{12}=\a_4.
\]
The corresponding generators $Y_{i_k,p_k}$ of $\Y_{t,\,Q}$ can be read 
at the corresponding place on the left side of Figure~\ref{FigD}. 
Let $\bd=(d_1,d_2,d_3,d_4)$ be a dimension vector for $Q$.
Then
\[
Y^{W^{\bd}} = Y_{1,0}^{d_1}Y_{2,0}^{d_2}Y_{3,5}^{d_3}Y_{4,0}^{d_4}. 
\]
The elements of $I_\bd$ are $12$-tuples $\aa\in\N^{12}$ encoding the decompositions
of $\bd$ into a sum of positive roots. 
By Proposition~\ref{propgeo2},
they are in one-to-one correspondence 
with the dominant monomials of the form $Y^{W^\bd} A^V$.
This bijection 
can be read immediately from Figure~\ref{FigD}.

For example, if
$\bd = (1,1,1,1)\equiv \a_1+\a_2+\a_3+\a_4$, 
the correspondence is:
\[
\begin{array}{lclcl}
(\a_1)+(\a_2)+(\a_3)+(\a_4) &\leftrightarrow & Y_{1,0}Y_{2,0}Y_{3,5}Y_{4,0}
&\leftrightarrow &1
,\\ [1mm]
(\a_1+\a_3)+(\a_2)+(\a_4) &\leftrightarrow & Y_{1,4}Y_{2,0}Y_{4,0}
&\leftrightarrow &A_{1,1}A_{3,2}A_{2,3}A_{4,3}A_{3,4}
,\\ [1mm]
(\a_2+\a_3)+(\a_1)+(\a_4) &\leftrightarrow & Y_{2,4}Y_{1,0}Y_{4,0}
&\leftrightarrow &A_{2,1}A_{3,2}A_{1,3}A_{4,3}A_{3,4}
,\\ [1mm]
(\a_3+\a_4)+(\a_1)+(\a_2) &\leftrightarrow & Y_{4,4}Y_{1,0}Y_{2,0}
&\leftrightarrow &A_{4,1}A_{3,2}A_{1,3}A_{2,3}A_{3,4}
,\\ [1mm]
(\a_1+\a_2+\a_3)+(\a_4) &\leftrightarrow & Y_{4,2}Y_{4,0}
&\leftrightarrow &A_{1,1}A_{2,1}A_{3,2}^2A_{1,3}A_{2,3}A_{4,3}A_{3,4}
,\\ [1mm]
(\a_1+\a_3+\a_4)+(\a_2) &\leftrightarrow & Y_{2,2}Y_{2,0}
&\leftrightarrow &A_{1,1}A_{4,1}A_{3,2}^2A_{1,3}A_{2,3}A_{4,3}A_{3,4}
,\\ [1mm]
(\a_2+\a_3+\a_4)+(\a_1) &\leftrightarrow & Y_{1,2}Y_{1,0}
&\leftrightarrow &A_{2,1}A_{4,1}A_{3,2}^2A_{1,3}A_{2,3}A_{4,3}A_{3,4}
,\\  [1mm]
(\a_1+\a_2+\a_3+\a_4)   &\leftrightarrow & Y_{3,1}
&\leftrightarrow &A_{1,1}A_{2,1}A_{4,1}A_{3,2}^2A_{1,3}A_{2,3}A_{4,3}A_{3,4}
.
\end{array}
\]
It is obtained by replacing each root $\b_k$ in a decomposition of $\a_1+\a_2+\a_3+\a_4$
by the corresponding variable $Y_{i_k,p_k}$.
The third column gives the monomial 
$Y^{W^\bd}\left(\prod_{k=1}^r Y_{i_k, p_k}^{a_k}\right)^{-1}$.
} 
\end{example}

We can now state the main result of this section.

\begin{Thm}\label{thmgeo}
\begin{itemize}
\item[\rm(a)] We have a $G_\bd$-equivariant isomorphism of varieties 
$\Psi \colon \M^\bullet_0(W^\bd)\overset{\sim}{\longrightarrow} E_\bd$.
\item[\rm(b)] $\M^\bullet_0(W^\bd)$ is an affine space of dimension
$\sum_{i\to j} d_id_j$.
\item[\rm(c)] Lusztig's perverse sheaves $IC(\O_\aa)$ on $E_\bd$ are the same
as Nakajima's perverse sheaves $IC_{W^{\bd}}(V)$ on $\M^\bullet_0(W^\bd)$.
\end{itemize}
\end{Thm}

\proof
In Proposition~\ref{propgeo1} we have constructed a 
$G_\bd$-equivariant closed immersion 
$\Psi$ of  $\M^\bullet_0(W^\bd)$ into $E_\bd$.
Since each stratum $\M^{\bullet\, {\rm reg}}_0(V,W^\bd)$
is $G_{\bd}$-invariant,
$\Psi$ maps every stratum to a union of orbits $\O_\aa$.
Since $\Psi$ is injective and 
the number of strata is equal to the number of orbits (Proposition~\ref{propgeo2}),
it follows that $\Psi$ maps each 
stratum of $\M^\bullet_0(W^{\bd})$
to a single $G_\bd$-orbit in $E_\bd$, so every orbit is 
contained in the image of $\Psi$. Thus $\Psi$ is surjective.
Since a surjective closed immersion between reduced schemes is an 
isomorphism, this proves (a). Claim (b) follows immediately from (a), and claim (c)
is again a consequence of Proposition~\ref{propgeo2}, which shows that the stratifications
used for defining the perverse sheaves are the same.
\cqfd

\begin{remark}\label{2remarks}
{\rm
(a)\ \ By the proof of Proposition~\ref{propgeo2}, for every 
dominant monomial $m$ in $\Y_{t,\,Q}$ there is a unique
pair $(V, W^\bd)$ such that $m = Y^{W^\bd}A^V$. Hence, even if
the varieties $\M^\bullet_0(W^{\bd})$ involve very particular spaces $W^\bd$,
the isomorphisms $\Psi\colon \M^\bullet_0(W^{\bd})\overset{\sim}{\longrightarrow} E_\bd$ 
are enough to identify
all the irreducible $(q,t)$-characters of $\CC_{Q}$.

Thus Theorem~\ref{LusztigTh}, Theorem~\ref{NakajimaTh},
and Theorem~\ref{thmgeo}
provide a geometric explanation of part~(b) of Theorem~\ref{mainth}.
By comparing convolution diagrams in \cite{Lu2} and \cite{VV}, it 
should also be possible to understand 
in a geometric manner
part (a) of Theorem~\ref{mainth},
that is, the multiplicative structure (see \cite[\S3.5]{Ncl}). 

\medskip\noindent
(b)\ \
If we take for $Q$
a quiver of type $A$ with all arrows in the same direction,
then 
$\M^\bullet_0(W^{\bd})$ is just a 
space of graded nilpotent endomorphisms
as in the Ginzburg-Vasserot construction \cite{GV} of type $A$
quantum loop algebras
(see \eg\cite[\S2.5.3]{L}). So 
Theorem~\ref{thmgeo} becomes tautological in this case.
}
\end{remark}

\subsection*{Acknowledgements}
We thank H. Nakajima for helpful comments and discussions, in particular
about the proof of Theorem~\ref{thmgeo}.


\bigskip
\small
\noindent
\begin{tabular}{ll}
David {\sc Hernandez} : &Universit\'e Paris Diderot-Paris 7,\\
& Institut de Math\'ematiques de Jussieu - Paris Rive Gauche CNRS UMR 7586,\\
& B\^atiment Sophie Germain, Case 7012,\\ 
& 75205 Paris Cedex 13, France. \\
&email : {\tt hernandez@math.jussieu.fr}\\[5mm]
Bernard {\sc Leclerc} : & Normandie Univ, France\\ 
&UNICAEN, LMNO F-14032 Caen, France\\
&CNRS UMR 6139, F-14032 Caen, France\\
&Institut Universitaire de France,\\
&email : {\tt bernard.leclerc@unicaen.fr}
\end{tabular}
\end{document}